\documentclass [10pt]{article}      
\usepackage[margin=1.5in]{geometry}
\usepackage{amsmath} 
\usepackage{latexsym}
\usepackage{pb-diagram}
\usepackage{amssymb}
\usepackage{color}
\usepackage{verbatim}
\usepackage{graphicx}
\usepackage{hyperref}
\hypersetup{
    colorlinks, 
    citecolor=blue,
    filecolor=blue,
    linkcolor=blue,
    urlcolor=blue
}
\hypersetup{linktocpage}

\usepackage{todonotes}

\newtheorem{theorem}{Theorem}
\newtheorem{corollary}[theorem]{Corollary}

\newtheorem{lemma}[theorem]{Lemma}
\newtheorem{definition}[theorem]{Definition}
\newtheorem{prop}[theorem]{Proposition}
\newtheorem{proposition}[theorem]{Proposition}

\newtheorem{example}[theorem]{Example}
\newtheorem{fact}[theorem]{Fact}
\newtheorem{remark}[theorem]{Remark}

\newtheorem{numreq}{Numerical Requirement}
\newtheorem{sublemma}[theorem]{Sublemma}

\newcommand{\cat}{^\frown}
\newcommand{\rest}{\ensuremath{\upharpoonright}}

\newcommand{\inv}{{^{-1}}}

\newcommand{\rev}[1]{\mathop{\rm rev}({#1})}

\newcommand{\floor}[1]{{\lfloor{#1}\rfloor}}
\newcommand{\ceil}[1]{{\lceil{#1}\rceil}}

\newcommand{\la}{\langle}
\newcommand{\ra}{\rangle}

\newcommand{\Chi}{\ensuremath{\mbox{\Large{$\chi$}}}}

\newcommand{\sig}{\ensuremath{{\Sigma\!\!\!\!_{_\sim}\, }}}

\newcommand{\trees}{\mathcal{T}\!\! rees}

\newcommand{\xbmt}{\ensuremath{(X,\mcb,\mu,T)}}
\newcommand{\ycns}{\ensuremath{(Y,\mcc,\nu,S)}}

\newcommand{\MPT}{\mbox{\bf MPT}}

\newcommand{\zoo}{{[0,1)}}

\newcommand{\mudiff}{\mbox{Diff}^k(M, \mu)}

\newcommand{\difflam}{\mbox{Diff}^\infty(\mathbb T^2, \lambda)}

\newcommand{\sz}{\Sigma^\poZ}

\newcommand{\poZ}{\mathbb Z}
\newcommand{\nn}{{\mathbb N}}

\newcommand{\poz}{\mathbb Z}
\newcommand{\bk}{{\mathbb K}}
\newcommand{\bl}{\mathbb L}
\newcommand{\bt}{\mathbb T}

\newcommand{\mct}{{\ensuremath{\mathcal T}}}

\newcommand{\mci}{\ensuremath{\mathcal I}}
\newcommand{\mcl}{\ensuremath{\mathcal L}}
\newcommand{\mcq}{\ensuremath{\mathcal Q}}

\newcommand{\mcr}{\ensuremath{\mathcal R}}

\newcommand{\mcw}{\mathcal W}
\newcommand{\mce}{{\mathcal E}}

\newcommand{\mcc}{{\mathcal C}}
\newcommand{\mck}{{\mathcal K}}
\newcommand{\mcb}{\mathcal B}
\newcommand{\mcs}{{\mathcal S}}
\newcommand{\mco}{\mathcal O}
\newcommand{\mcf}{\mathcal F}

\newcommand{\dbar}{\ensuremath{\bar{d}}}

\newcommand{\bfni}[1]{\noindent {{\bf{#1}}}}

\newcommand{\qed}{{\nopagebreak \hfill $\dashv$ 
 \par\bigskip}}

\newcommand{\pf}{{\par\noindent{$\vdash$\ \ \ }}}

\newcommand{\qn}{{q_n}}
\newcommand{\qnpo}{{q_{n+1}}}

\newcommand{\kn}{{k_n}}

\newcommand{\misal}{\not\Downarrow}

\newcommand{\hoo}[2]{{[{#1}/{#2}, ({#1}+1)/{#2}) }}

\newcommand{\boundary}{\partial}

\newcommand{\wcuprev}[2]{{\mcw_{#1}/\mcq^{#1}_{#2} \cup \rev{\mcw_{#1}/\mcq^{#1}_{#2}}}}

 \title{Measure Preserving Diffeomorphisms of the Torus\\ are Unclassifiable}
 \author{Matthew Foreman and Benjamin Weiss}
 
\begin{document}
 \maketitle

  \begin{abstract}
In 1932 von Neumann proposed classifying the statistical behavior of differentiable systems. In modern language this is interpreted as classifying diffeomorphisms of compact manifolds up to measure isomorphism.  This paper proves that this is impossible in a rigorous sense.

\end{abstract}

 \tableofcontents

\pagebreak

 \section{Introduction}
  The isomorphism problem in ergodic theory was formulated by von Neumann in 1932 in his pioneering paper \cite{vN}.\footnote{Two measure  preserving transformations (abbreviated to `MPTs' in the paper) $T$ and $S$ are isomorphic if there is an invertible measurable mapping between
the corresponding measure spaces which commutes with the actions of $T$ and $S$.} The problem
 has been solved for some  classes of  transformations that have special properties. Halmos and von Neumann \cite{HvN} used the unitary operators defined by Koopman to completely characterize ergodic
measure preserving transformations  with pure point spectrum. They showed that these are exactly the transformations that can be realized  as translations on compact groups. Another
 notable success in solving this problem was the classification of Bernoulli shifts using the notion of entropy introduced by Kolmogorov. 

Starting in the late 1990's a different type of result began to appear: \emph{anti-classification} results that demonstrate in a rigorous way that classification is not possible. This type of theorem requires a precise definition of what a classification is. Informally, a classification is a {method} of determining isomorphism between transformations  by computing (in a liberal sense) other invariants for which equivalence is easy to determine.

The key words here are \emph{method} and \emph{computing}. For negative 
theorems, the more liberal a notion one takes for these words, the stronger the 
theorem. One natural way of what a computation is  uses  the Borel/non-Borel 
distinction. Saying a set $X$ or function $f$ is Borel is a loose way of saying that
 membership in $X$ or the computation of $f$ can be done using a countable 
 (possibly transfinite) protocol whose basic input is membership in open sets. Saying 
 that $X$ or $f$ is \emph{not} Borel is saying that determining membership in $X$ 
 or computing $f$ cannot be done with any countable amount of resources. (See 
 \cite{what_is} for an elementary discussion and a comparison with the more strict
notion of \emph{recursive computation}, which requires inherently finite 
resources.)

In the context of classification problems, saying that an equivalence relation $E$ on a space $X$ is \emph{not} Borel is
saying that there is no countable amount of initial  information and no countable, potentially transfinite, protocol based on this information for determining, for arbitrary
$x,y\in X$ whether $xEy$. {Any} such method must inherently use uncountable
resources.\footnote{Many well known classification theorems have as  immediate corollaries that the resulting equivalence
relation is Borel. An example of this is the Spectral Theorem, which has a consequence that the relation of Unitary
Conjugacy for normal operators is a Borel equivalence relation.} 

An example of a positive theorem in the context of ergodic theory is due to Halmos (\cite{Halmos}) who showed that the collection of ergodic measure preserving transformations is a dense $\mathcal G_\delta$ set in the space of all measure preserving transformations of $([0,1],\lambda)$ endowed with the weak topology. Moreover he showed that the set of weakly mixing transformations is also a dense $\mathcal G_\delta$.\footnote{Relatively
straightforward arguments show that the set of strongly mixing transformation is a first category $\Pi_0^3$ set. See \cite{descview}.}

The first anti-classification result in the area is due to Beleznay and Foreman \cite{BF} who showed that the class of \emph{measure distal} transformations used in early ergodic theoretic proofs of Szemeredi's theorem is not a Borel set. Later Hjorth \cite{hj} introduced the notion of \emph{turbulence} and showed that there is no Borel way of attaching algebraic invariants to ergodic transformations that completely determine isomorphism. Foreman and Weiss \cite{FW} improved this result by showing that the conjugacy action of the measure preserving transformations is turbulent--hence no generic class  can have a complete set of algebraic invariants.

In considering the isomorphism relation as a collection $\mathcal I$ of pairs $(S,T)$ of measure preserving transformations, Hjorth (\cite{hjorth}) showed that $\mathcal I$ is not a Borel set. However the pairs of transformations he used to demonstrate this were inherently non-ergodic,  leaving open the essential problem:

\begin{center}
Is isomorphism of ergodic measure preserving transformations Borel?
\end{center}

This question was answered in the negative by Foreman, Rudolph and Weiss in \cite{FRW}. This answer can be interpreted as saying that determining isomorphism between ergodic transformations is inaccessible to countable methods that use countable amounts of information.

In the same foundational paper from 1932 von Neumann expressed the likelihood that any abstract MPT is isomorphic to a continuous MPT and
perhaps even to a {differentiable} one. This brief remark eventually gave rise to one of the yet outstanding problems in smooth dynamics, namely:

\begin{center}
 Does every
ergodic MPT with finite entropy have a smooth model?
\footnote{In \cite{vN} on page 590,
``Vermutlich kann sogar zu jeder allgemeinen Str\"{o}mung eine isomorphe stetige Str\"{o}mung gefunden werden [footnote 13], vielleicht sogar eine stetig-differentiierbare, oder gar eine mechanische.
Footnote 13: Der Verfasser hofft, hierf\"{u}r demn\"{a}chst einen Beweis anzugeben."}
\end{center}

  By a smooth model it is meant an isomorphic copy of the MPT which is given by smooth diffeomorphism of a compact manifold preserving a measure equivalent to the volume element.
Soon after entropy was introduced, A. G. Kushnirenko showed that such a diffeomorphism must have finite entropy, and up to now this is the only
restriction that is known. The current paper is the culmination of  a series whose purpose is to show that the variety of ergodic transformations that have  smooth models is rich enough so
that the abstract isomorphism relation, when restricted to these smooth systems, is as complicated as the general isomorphism problem for ergodic measure preserving systems. We show that even when restricting to
diffeomorphisms of the 2-torus that preserve Lebesgue measure this is the case. The formal statement of our solution to the isomorphism problem is:

\begin{theorem}\label{big tuna} If $M$ is either the torus $\bt^2$, the disk $D$ or the annulus then the measure-isomorphism relation among pairs $(S,T)$ of measure preserving $C^\infty$-diffeomorphisms of $M$ is not a Borel set with respect to the $C^\infty$-topology.
\end{theorem}
Thus the isomorphism problem is impossible even for diffeomorphisms of compact surfaces.
\bigskip

How does one prove a result such as Theorem \ref{big tuna}? The main tool is the idea of a \emph{reduction} (see \cite{what_is} and Section \ref{dst basics}). A function $f:X\to Y$ reduces $A$ to $B$ if and only for all $x\in X$:
\[x\in A \mbox{ if and only if }f(x)\in B.\]
If $X$ and $Y$ are completely metrizable spaces and $f$ is a Borel function then $f$ is a method of reducing the question of membership in $A$ to membership in $B$. Thus if $A$ is not Borel then $B$ cannot be either.

In the current context, the $C^\infty$-topology on the smooth transformations refines the weak topology. Thus, by Halmos' result quoted earlier, on the torus (disc etc.), the ergodic transformations are still a $\mathcal G_\delta$-set. (However the famous KAM theory shows that the ergodic transformations are no longer dense.) In particular the $C^\infty$-topology induces a  metrizable complete and perfect topology on the measure preserving diffeomorphisms of $\bt^2$. If $M$ is a manifold with supporting a measure $\mu$ we denote the space of $C^\infty$, $\mu$-measure preserving diffeomorphisms of $M$ with the notation $\mbox{Diff}^\infty(M, \mu)$. {Elements of 
$\mbox{Diff}^\infty(M, \mu)$ are also members of the group MPT of $\mu$-measure preserving transformations. 
For $T\in \mbox{Diff}^\infty(M, \mu)$ the centralizer of $T$ in MPT is denoted $C(T)$.}

If $X$ is perfect and completely metrizable, a set $A\subseteq X$ is \emph{analytic} if and only $A$ is the continuous image of a Borel set. A is \emph{complete analytic} if and only if every analytic set can be reduced to $A$. It is a classical fact that complete analytic sets are not Borel.

The proof of Theorem \ref{big tuna} uses a well-known example of a complete analytic set. The underlying space $X$ is the space $\trees$ and $A$ is the collection of ill-founded trees; those that have infinite branches.  A precise  statement of the main result of the paper: 
\begin{theorem}\label{lele}
There is a continuous function $F^s:\trees\to \difflam$, taking values among the ergodic transformations, such that for $\mct\in \trees$, if $T=F^s(\mct)$:
\begin{enumerate}
\item $\mct$ has an infinite branch if and only if $T\cong T^{-1}$, and
\item $\mct$ has two distinct infinite branches if and only if 
\[C(T)\ne \overline{\{T^n:n\in\poZ\}}.\]
\end{enumerate}
\end{theorem}
\begin{corollary}\ \label{cor3}
\begin{itemize}
\item $\{T\in \difflam:T\mbox{ is ergodic and }T\cong T\inv\}$ is complete analytic.
\item $\{T\in \difflam:T\mbox{ is ergodic and }C(T)\ne \overline{\{T^n:n\in\poZ\}}$ is complete analytic.
\end{itemize}
\end{corollary}
Since the map $\iota(T)=(T, T\inv)$ is a continuous mapping of $\difflam$ to $\difflam\times \difflam$ and reduces $\{T:T\cong T\inv\}$ to $\{(S,T):S\cong T\}$, it follows that:
\begin{corollary}
$\{(S,T):S $ and $T$ are ergodic diffeomorphisms of $\bt^2$ and are isomorphic$\}$ is a complete analytic set and hence not Borel.
\end{corollary}

{We note that the problem of finding even one measure preserving transformation not isomorphic to its inverse is difficult. This was not done until Anzai in \cite{anzai}. In Math Review MR0047742, Halmos said, ``By constructing an example of the type described in the title the author solves (negatively) a problem proposed by the reviewer and von Neumann [Ann. of Math. (2) 43, 332?350 (1942): MR0006617]".}

More fine-grained information is now known and will be published elsewhere. For example, Foreman, in unpublished work, showed that the problem of ``isomorphism of countable graphs" is Borel reducible to the isomorphism problem for ergodic measure preserving transformations. 

The techniques of this paper also have foundational interest. {A close analysis of our construction shows that the problem of whether $T$ is isomorphic to its inverse is ``$\Pi^0_1$-hard."  (See \cite{GD}).}  This enables one to prove 
that truth or falsity of various open problems like the Riemann hypothesis is equivalent to the question of is $T_{RH}$ isomorphic 
or not to its inverse for a specific measure preserving diffeomorphism $T_{RH}$ of the torus  given by our construction.  Another consequence is the existence of a different diffeomorphism $T_{ZFC}$ such that the question of whether $T_{ZFC}$ is isomorphic to its inverse is independent of ZFC, the usual axioms for mathematics.
\medskip

Here are two problems that remain open: 

\begin{description}
\item[Problem 1] In contrast to  \cite{FW}, where the authors were able to show that the equivalence relation of isomorphism on abstract ergodic measure preserving transformations is \emph{turbulent}, this remains open for ergodic diffeomorphisms of a compact manifold. 

\item[Problem 2] The problem of classifying diffeomorphisms of compact surfaces up to topological conjugacy remains largely open. {Work of the first author with A Gorodetski shows that the isomorphism relation itself is not Borel, but for a very specific type of diffeomorphisms of manifolds of dimension 5 and above. It is not know, for example for topologically minimal transformations.}
\end{description}

We owe a substantial debt to everyone who has helped us with this project. Jean-Paul Thouvenot brought the Anosov-Katok technique to our attention and suggested using it to solve the von Neumann problem. Philipp Kunde aided us by reading the paper and providing comments and corrections. Others include  Eli Glasner, Anton Gorodetski,  Alekos Kechris, and Anatole Katok.  

We particularly want to acknowledge the contribution of the late Dan Rudolph, who helped pioneer these ideas and was a co-author in \cite{FRW}, contributing techniques fundamental to this paper.

\section{An Outline of the Argument}\label{outline}

This section gives an outline of the argument for Theorem \ref{lele}. It uses the main results from our earlier papers:
\emph{A Symbolic Representation of {A}nosov-{K}atok Systems}(\cite{prequel}) and \emph{From Odometers to Circular Systems: A Global Structure Theorem} (\cite{global_structure})
 which we briefly summarize.
   In \cite{prequel}, the Anosov-Katok technique of Approximation by Conjugacy is used to  give a new symbolic representation for a class of measure preserving diffeomorphisms that are extensions of the rotations
by certain Liouvillean $\alpha$. These are called
strongly uniform \textbf{Circular Systems}.\footnote{In a forthcoming paper we show how to drop the ``strongly uniform" assumption.}

   In \cite{global_structure} two classes of symbolic systems are defined.  The first, called \emph{Odometer Based} systems, contains representatives of every finite entropy measure preserving transformation with an odometer factor. The second class is the collection of \emph{Circular Systems}. These classes are made into categories by taking as morphisms synchronous and anti-synchronous factor maps.  The main result is that there is a functorial isomorphism between $\mcf$ between these categories that takes strongly uniform systems to strongly uniform  systems.

   Since the main construction in \cite{FRW} uses  Odometer Based systems this map enables us to adapt that construction
   to the smooth setting. However in order to prove our main result we still have to take into account potential isomorphisms
   of Circular Systems that are neither synchronous nor anti-synchronous. It is to deal with this difficulty that we analyze what we call
   the \emph{displacement function}.

To each $\alpha$ arising as a rotation factor of a circular system $T$ one can associate a \emph{displacement function} (Section \ref{def of delta}) and use it to associate the set of \emph{central values}, a subgroup of the unit circle. Its significance is 
the following:
\begin{itemize}
\item[1.]  (Theorem \ref{central values in closure}) If $\beta$ is central then there is an $\phi\in \overline{\{T^n:n\in \poZ\}}$ such that the rotation factor of $\phi$ is rotation by $\beta$.

\item[2.] (Theorem \ref{centralizer central values}) If $T$ is built from sufficiently random words,\footnote{i.e. $T$ satisfies the \emph{Timing Assumptions}.} and $\phi\in C(T)$, then the canonical rotation factor of $\phi$ is rotation by a central value.

\item[3.] It follows  that if there is a $\phi\in C(T)$ and $\phi\notin \overline{\{T^n:n\in \poZ\}}$, then there is a synchronous $\psi\in C(T)$ such that $\psi\notin \overline{\{T^n:n\in \poZ\}}$.

\item[4.] (Theorem \ref{conjugacy central values}) The analogous results {relating} isomorphisms $\phi$ between $T$ and $T^{-1}$ {with central values} 	is proved, allowing us to conclude that if $T$ is isomorphic to $T^{-1}$ then there is an anti-synchronous isomorphism between $T$ and $T^{-1}$. 

\item[5.] The previous two items are the content of Theorem \ref{perfect match}, which says that for $T$ satisfying the Timing Assumptions,   to decide whether  
$T\cong T^{-1}$ or $C(T)\ne \overline{\{T^n:n\in \poZ\}}$ it suffices to consider anti-synchronous and synchronous isomorphisms.

\end{itemize}
In \cite{FRW} a continuous function $F$ from the space of Trees to the strongly uniform odometer based transformations is constructed that:
\begin{itemize} 
\item reduces the set of ill-founded trees to the transformations $T$ that are isomorphic to their inverses (and \emph{if} $T\cong T^{-1}$ then this is witnessed by an anti-synchronous isomorphism) and
\item reduces the set of trees with two infinite branches to the transformations $T$ whose centralizer is different from the powers of $T$ (and if the centralizer contains an exotic element, it contains a synchronous exotic element).
\end{itemize}
Moreover, in the second case, there is a synchronous element of the centralizer with a specific piece of evidence that it is not the identity (it moves a $\mcq^1_1$-equivalence class).

Composing one concludes that $\mcf\circ F$:

\begin{itemize} 
\item reduces the set of ill-founded trees to collection of circular systems that are isomorphic to their inverses and
\item reduces the set of trees with two infinite branches to the circular systems whose centralizer is different from the closure of the powers of $T$. 
\end{itemize}

\begin{figure}[!h]
\centering
\includegraphics[height=.40\textheight]{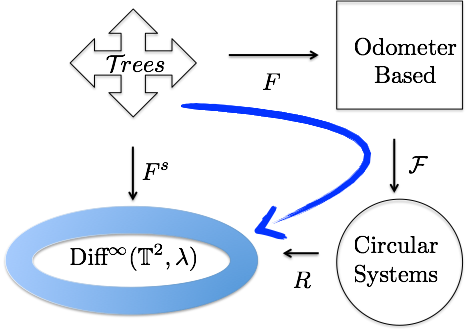}
\caption{The reduction $F^s$.}
\label{square diagram}
\end{figure}

Continuously realizing the circular systems by $R$ (as in \cite{prequel}) completes the proof that:

\begin{itemize} 
\item The collection of ergodic measure preserving diffeomorphisms $T$  of the torus that are isomorphic to their inverses is complete analytic. Consequently the set of pairs $(S,T)$ of ergodic conjugate  measure preserving diffeomorphisms is
a complete analytic set.
\item The collection of  ergodic measure preserving diffeomorphisms $T$ whose centralizer is different from the closure of the powers of $T$ is complete analytic.
\end{itemize}
Figure \ref{square diagram} illustrates $F^s=R\circ \mcf\circ G$.

\bigskip

 The next two sections  review   basic facts in ergodic theory and descriptive set theory,  define \emph{odometer based} and \emph{circular} systems and review their
 properties and the facts shown in \cite{prequel} and \cite{global_structure}.

   The analysis of the displacement function and the associated central values, which are a subgroup of the circle canonically
   associated to the Liouvillean $\alpha$, is carried out in sections 5-7. Finally the proof of the main theorems are given in section 8
   modulo certain properties which impose some additional conditions on the parameters of the construction in \cite{FRW}. These are
   verified in section 9 and in section 10  we spell out the dependencies between the various parameters and show that they can be
   realized.  
   
   \section{Numerical Requirements}\label{OMG}

The proof of Theorem \ref{lele} uses a construction with many interconnecting pieces, most of which are built by taking limits.  This results in a large number of related sequences of variables, each having their own requirements and the estimates for the different pieces must be compatible. 

The least interesting part of this paper is is verifying the consistency of the numerical requirements. Sorting these requirements out is completely independent of the rest of the paper.  For this reason, we list the numerical requirements in section \ref{alltogethernow}, and then give an argument for their consistency. We also note the specific requirement by number in the text as they are posited and used.

Contributing to the complexity of the situation is that many of the relationships between the variables come from internal arguments of the general form ``taking $\delta$ small enough you can guarantee that $x<\epsilon$", with various variables in place of $\epsilon, \delta$ and $x$.   The exact relationship between 
$\epsilon$ and $\delta$  is not clear from the argument, but there is a requirement of the form  \emph{``$\delta$ is small as a function of $\epsilon$."} 
A typical example of this is  Sublemma \ref{fortheref} which says that, as a function of $Q^n_1$, if $\epsilon_n$ is take sufficiently small then an involved inequality involving $I^*, u_i', v_i'$ and $Q^n_1$ holds.

Complicating this task further is the fact that the construction in this paper depends on the construction in \cite{FRW}, which has its own numerical requirements.   For a reader tracking the correspondence, in the appendix,  we include a  table for translating between the notation in this paper  and the notation in \cite{FRW}

\paragraph{The variables} Here is a list of variable sequences that have to be chosen during the construction:

\[k_n, l_n, q_n, s_n, e(n), p_n, q_n, \alpha_n,  \epsilon_n, \varepsilon_n, \mu_n, Q^n_1.\]

Some of these variables have clear relationships that are externally determined.  The main construction is of a function that has a tree as in input.  That tree \emph{directly} determines a sequence of parameters, such as $G^n_1$ and $\la M(s):s\le n\ra$ that are not \emph{chosen} during the construction.  (In section \ref{numreqs}, we call these \emph{exogenous} variables.) These parameters determine some of the numerical requirements. 

\begin{example} \label{to confuse the reader} The words in the collection $\mcw_{n+1}$ are built by a sequence of $M$ substitutions into equivalence classes of the relations $\mcq^{n+1}_i$, where 
$M=\sup_SM(s)$ for $S$ the collection of heights on nodes in the given tree at stage $n$. These substitution instances 
are closed under a sequence of $\poZ_2$ actions of the groups $\la G^n_i:i\le M\ra$.  The number $M$ and $\poZ_2$ dimension of the 
$\poZ_2$ actions is also determined by the tree.  Thus $s_{n+1}$ is determined by the exogenous variables $G^n_i$, 
$M(s)$, and the internally chosen variable $e(n+1)$. In this particular example, It is possible to give a completely explicit formula for $e(n+1)$ in terms $s_{n+1}$ and vice versa.\footnote{{$s_{n+1}=(2^{Me(n+1)})G$ for numbers $M$ and $G$ determined exogenously.}}

However that would be uninformative. What we need to see is that if $e(n+1)$ is large then $s_{n+1}$ is and vice versa and that each determines the other.   This is the only relevant information for determining the consistency of the numerical requirements. We have thus eliminated one variable.
\end{example}

\medskip

It would perhaps be more conventional to define all of the variables in advance, write down the list of inequalities and then show they are consistent.  However the examples above illustrate the difficulties with this. The inequalities are intimately intertwined with the details of the construction and are completely enigmatic without that context.  For this reason we note the Numerical Requirements one by one as they accumulate and collect them in section \ref{alltogethernow}. We then proceed to show that they are consistent by the method we describe next.  A reader with a preference for the conventional presentation is advised to skip directly to section \ref{numreqs}, read the reconciliation and then return to read the rest of the paper.

\paragraph{What could possibly go wrong?} The only potential issue is that that there may be a situation where the requirements are circular: for example, $\delta$ might have to be small as a function of $\epsilon$, $\epsilon$ small as a function of $\mu$ and $\mu$ small as a function of $\delta$.  In symbols 

\[\begin{diagram}%
\node{\epsilon}\arrow{e}\node{\delta}\arrow{e}\node{\mu}\arrow{e}\node{\epsilon}
\end{diagram}\]
So if you choose $\epsilon$ first, then $\delta$ then $\mu$ you might find that your choice of $\epsilon$ was inadequate. Indeed, because there is a cycle in the dependency diagram there is no variable you can choose first and be certain of consistency.

\paragraph{Method for showing consistency} In section \ref{numreqs} we analyze the dependencies and draw a dependency diagram giving the order of choice. Since that diagram is cycle free, all of the variables can be chosen to satisfy the accumulated requirements.

\section{Preliminaries}
The reader is referred to standard texts such as \cite{Shields}, \cite{walters} or \cite{Peterson}. Facts that are not standard and are simply cited here are proved in \cite{global_structure}, \cite{prequel} and \cite{FRW}.

\subsection{Measure Spaces}\label{abstract measure spaces} We will call separable non-atomic probability spaces \emph{standard measure spaces} and denote them $(X,\mcb, \mu)$ where $\mcb$ is the Boolean algebra of measurable subsets of $X$ and $\mu$ is a countably additive, non-atomic measure defined on $\mcb$. 
Maharam and von Neumann proved that every standard measure space is isomorphic to $([0,1],\mcb,\lambda)$ where $\lambda$ is Lebesgue measure and $\mcb$ is the algebra of Lebesgue measurable sets.

If $(X, \mcb, \mu)$ and $(Y, \mcc, \nu)$ are measure spaces, an isomorphism between $X$ and $Y$ is a bijection $\phi:X\to Y$ such that $\phi$ is measure preserving and both $\phi$ and $\phi^{-1}$ are measurable. We will ignore sets of measure zero when discussing isomorphisms; i.e.  we allow  the domain and range of $\phi$ to be  subsets of $X$  and $Y$  of measure one.

A measure preserving system is an object $\xbmt$ where $T:X\to X$ is a measure isomorphism. A \emph{factor map} 
between two measure preserving systems $\xbmt$ and $\ycns$ is a measurable, measure preserving  function 
$\phi:X\to Y$ such that $S\circ\phi=\phi\circ T$. A factor map is an isomorphism between systems iff $\phi$ is a 
measure isomorphism.

Let $T:\xbmt\to \xbmt$ be measure preserving, $(Y,\mcc)$ be a measurable space, $S:Y\to Y$ a measurable map and $\phi:X\to Y$ be a measurable map such that $\phi T=S\phi$.  Then we can define a measure $\nu=\phi^*\mu$ by setting $\nu(A)=\mu(\phi^{-1}(A))$. This measure makes $\phi$ a factor map from $\xbmt$ to $\ycns$.
\medskip

\subsection{Presentations of Measure Preserving Systems}\label{presentations}
Measure preserving systems occur naturally in many guises with diverse topologies. As far as is known, the Borel/non-Borel distinction for dynamical properties is the same in each of these presentations and many of the presentations have the same generic classes. (See the forthcoming paper \cite{models} which gives a precise condition for this.)

Here is a  review the properties of the types of presentations relevant to this paper, which are:
abstract invertible preserving systems, smooth transformations preserving volume elements and symbolic 
systems. 

\subsubsection{Abstract Measure Preserving systems}
Since every standard measure space is isomorphic to the unit interval with Lebesgue measure, every invertible measure preserving transformation of a standard measure space is isomorphic to an invertible Lebesgue measure preserving transformation on the unit interval. 

In accordance with the conventions of \cite{descview} we denote the group of measure preserving transformations of $\zoo$ by \MPT.\footnote{Recently several authors have adopted the notation $Aut(\mu)$ for the same space.} Two measure preserving transformations are identified if they are equal on sets of full measure. 

Two measure preserving transformations are isomorphic if and only if they are conjugate in the group \MPT\ and we will use \emph{isomorphic} and \emph{conjugate} as synonyms. However some caution is order. If $(M,\mu)$ is a manifold, $T:M\to M$ is a smooth measure preserving transformation and $\phi$ is an arbitrary measure preserving transformation from $M$ to $M$, then $\phi T\phi^{-1}$ is unlikely to be smooth.  Thus, the equivalence relation of isomorphism of diffeomorphisms is not given by an action of the group of measure preserving transformations in an obvious way.

Given a measure space $(X,\mu)$ and a measure preserving transformation $T:X\to X$,  define the \emph{centralizer of $T$} to be 
the collection of measure preserving $S:X\to X$ such that $ST=TS$. This group is denoted $C(T)$. Note that this is the centralizer \emph{in the group of measure preserving transformations}.  In the case that $X$ is a manifold and $T$ is a diffeomorphism, $C(T)$ differs from  the centralizer of $T$ inside the group of diffeomorphisms.

To each invertible measure preserving transformation $T\in \MPT$,  associate a unitary operator $U_T:L^2([0,1])\to L^2([0,1])$ by defining $U(f)=f\circ T$.  In this way \MPT\ can be identified with a closed subgroup of the unitary operators on $L^2([0,1])$ with respect to the weak operator 
topology\footnote{Which coincides with the strong operator topology in this case.} on the space of  unitary 
transformations. This makes \MPT\ into a Polish group. We will call this the \emph{weak topology} on \MPT. 
Halmos (\cite{Halmos}) showed that the ergodic transformations, which we denote
 $\mce$, is a dense $\mathcal G_\delta$ set in $\MPT$. In particular the weak topology makes $\mce$ into a Polish subspace of $\MPT$.

 There is another topology on the collection of measure preserving transformations of $X$ to $Y$  for measure spaces $X$ and $Y$. If $S, T:X\to Y$ are measure preserving transformations, the \emph{uniform distance} between $S$ and $T$ is defined to be:
\[d_U(S, T)=\mu\{x:Sx\ne Tx\}.\]
This topology refines the weak topology and is a complete, but  not a separable topology.

\subsubsection{Diffeomorphisms}
Let $M$ be a {$C^m$}-smooth compact finite dimensional manifold and $\mu$ be a standard measure on $M$ determined by a smooth volume element. 
For each {$k\le m$} there is a Polish topology on the $k$-times differentiable homeomorphisms of $M$, the $C^k$-topology. {If $M$ is $C^\infty$, then }the $C^\infty$-topology is the coarsest topology refining the $C^k$-topology for each $k\in \nn$. 
It is also a Polish topology and a sequence of $C^\infty$-diffeomorphisms converges in the $C^\infty$-topology if and only if it converges in the $C^k$-topology for each $k\in \nn$. 

The collection of $\mu$-preserving  diffeomorphisms forms a closed nowhere dense set in the $C^k$-topology on the $C^k$-diffeomorphisms, and as such, inherits a Polish topology.\footnote{One 
can also consider the space of measure preserving homeomorphisms with the $\|\ \|_\infty$ topology, 
which behaves in some ways similarly.}
We will denote this space by $\mudiff$.

Viewing $M$ as an abstract measure space one can also consider the space of abstract $\mu$-preserving transformations on $M$ with the weak topology. In 
\cite{bjorkthesis} it is shown that the collection of a.e.-equivalence classes of smooth transformations form  a 
$\Pi^0_3$-set in \MPT(M), and hence the collection has the  Property of Baire.

\subsubsection{Symbolic Systems}
\label{symbolic shifts}
Let $\Sigma$ be a countable or finite alphabet endowed with the discrete topology. Then $\Sigma^\poZ$ can be given the product topology, which makes it into a separable, totally disconnected space that is compact if $\Sigma$ is finite.

\medskip

\bfni{Notation:}  If $u=\la \sigma_0, \dots \sigma_{n-1}\ra\in \Sigma^{<\infty}$ is a finite sequence of elements of $\Sigma$, then we denote the cylinder set based at $k$ in $\Sigma^\poZ$ by writing $\la u\ra_k$. If $k=0$ we abbreviate this and write $\la u\ra$. 
Explicitly: $\la u\ra_k=\{f\in \Sigma^\poZ: f\rest[k,k+n)=u\}$. The collection of cylinder sets form a base for the product topology on $\Sigma^\poZ$.

  Let $u, v$ be finite sequences of elements of $\Sigma$ having length $q$.  Given intervals $I$ 
  and $J$ in $\poZ$ of length $q$ we can view $u$ and $v$ as functions having domain $I$ and $J$ 
  respectively. We will say that $u$ and $v$ are located at $I$ and $J$.  We will say that $u$ is \emph{shifted by $k$} relative to $v$ iff $I$ is the shift of the interval 
  $J$ by $k$.  We say that $u$ is the \emph{$k$-shift} of $v$ iff $u$ and $v$ are the same words and $I$ is 
  the shift of the interval $j$ by $k$.

\medskip

\noindent The shift map:
\[sh:\Sigma^\poZ\to \Sigma^\poZ\]
defined by setting $sh(f)(n)=f(n+1)$ is a homeomorphism. If $\mu$ is a shift-invariant Borel measure then the resulting measure preserving system $(\Sigma^\poZ, \mcb,\mu, sh)$ is called a \emph{symbolic system}. The closed support of $\mu$ is a shift-invariant closed subset of $\Sigma^\poZ$ called a \emph{symbolic shift} or \emph{sub-shift}.

Symbolic shifts are often described intrinsically by giving a collection of words that constitute a clopen basis for the support of an invariant measure.  
Fix a language $\Sigma$,  and a sequence of collections of words $\la\mcw_n:n\in\nn\ra$ with the properties that:
\begin{enumerate}
\item for each $n$ all of the words in $\mcw_n$ have the same length $q_n$,
\item each $w\in\mcw_{n}$ occurs  at least once as a subword of every $w'\in \mcw_{n+1}$,
\item \label{not too much space} there is a summable sequence $\la \epsilon_n:n\in\nn\ra$ of positive numbers such that for each $n$, every word $w\in \mcw_{n+1}$ can be uniquely parsed into segments 
\begin{equation}u_0w_0u_1w_1\dots  w_lu_{l+1}\label{words and spacers}
\end{equation}
 such that each $w_i\in \mcw_n$, {$u_i\in \Sigma^{<q_n}$} and for this parsing
\begin{equation}\label{small boundary numeric} {\sum_i|u_i|\over q_{n+1}}<\epsilon_{n+1}.
\end{equation}
\end{enumerate}

\noindent The segments $u_i$ in condition \ref{words and spacers} are called the \emph{spacer} or \emph{boundary} portions of $w$.
\begin{definition}A sequence $\la \mcw_n:n\in\nn\ra$ satisfying properties 1.)-3.) will be called a \emph{construction sequence}.
\end{definition}

If $\mcw$ is a collection of words in an alphabet $\Sigma$, we will say that 
$\mcw$ is \emph{uniquely readable} if and only if whenever $u, v, w\in \mcw$ and $uv=pws$ then either:
\begin{itemize}
\item $p=\emptyset$ and $u=w$ {or}
\item $s=\emptyset$ and $v=w$.
\end{itemize}
Equation 1 of clause 3 implies that each $\mcw_n$ is uniquely readable. {We will need unique readability to parse elements of $\bk$,  the symbolic shift associated with the construction sequence}.

{\begin{definition} Let $\bk$  be the collection of 
 $x\in \Sigma^\poZ$ such that every finite contiguous subword of $x$ occurs inside some $w$ belonging to some 
 $\mcw_n$. Then $\bk$ is a 
 closed shift-invariant subset of $\sz$ that is compact if $\Sigma$ is finite. 
 \end{definition}}
  
 The symbolic shifts built from construction sequences coincide with transformations built by \emph{cut-and-stack} constructions.
 \medskip
 
{ \bfni{Notation:} For a word $w\in \Sigma^{<\nn}$ we will write  $|w|$ for the length of $w$.}
\medskip

 Here is a natural set of measure one for the relevant measures:

\begin{definition}\label{def of S} Suppose that $\la \mcw_n:n\in\nn\ra$ is a construction sequence  for a symbolic system $\bk$ with each $\mcw_n$ uniquely readable. Let
$S$ be the collection $x\in \bk$ such that there are  sequences of natural numbers 
$\la a_m: m\in\nn\ra$, $\la b_m: m\in\nn\ra$ going to infinity such that for all $m$ there is an 
$n, x\rest [-a_m, b_m)\in \mcw_n$.
\end{definition}

Note that  $S$ is a dense shift-invariant $\mathcal G_\delta$ set.
\begin{lemma}\label{catch all} \cite{prequel}
Fix a construction sequence $\la\mcw_n:n\in\nn\ra$ for a symbolic system $\bk$ in a finite language. Then:
\begin{enumerate}
\item $\bk$ is the smallest shift-invariant closed subset of $\Sigma^\poZ$ such that for all $n$,  and $w\in\mcw_n$, $\bk$ has non-empty intersection with the basic open interval $\la w\ra\subset \Sigma^\poZ$.

\item Suppose that there is a unique invariant measure $\nu$ on $S\subseteq \bk$, then $\nu$ is ergodic.
\item (See \cite{global_structure}) If $\nu$ is an invariant measure on $\bk$ concentrating on $S$, then for  $\nu$-almost every $s$ there is an $N$ for all $n>N$, there are $a_n\le 0< b_n$ such that $s\rest[a_n, b_n)\in \mcw_n$.
\end{enumerate}
\end{lemma}

\begin{example}
\label{definition of uniform}
Let $\la\mcw_n:n\in\nn\ra$ be a construction sequence.  Then $\la \mcw_n:n\in\nn\ra$ is  \hypertarget{uniform}{\emph{uniform}} 
if there is a summable sequence of positive numbers $\la \epsilon_n:n\in\nn\ra$ and $\la d_n:n\in\nn\ra$, where $d_n:\mcw_n\to (0,1)$ such that for 
each $n$ all  words $w\in \mcw_n$ and $w'\in \mcw_{n+1}$ if $f(w,w')$ is the number of $i$ such that  $w=w_i$
\begin{equation}\label{uniform occurrance}
\left|{f(w,w')\over q_{n+1}/q_n}-d_n(w)\right|<{\epsilon_{n+1}\over q_n}.
\end{equation}
It is shown in \cite{prequel} that uniform construction sequences are uniquely ergodic. A  special case of uniformity is \emph{strong uniformity:}  when each $w\in \mcw_n$ 
occurs exactly the same number of times in each $w'\in \mcw_{n+1}$. This property holds for  the circular systems considered in \cite{prequel} and that are used for the proof of the main theorem of this paper (Theorem \ref{lele}).
\end{example}

\subsubsection{Locations}

Let $\la \mcw_n:n\in\nn\ra$ be a uniquely readable construction sequence and  let $\nu$ be a shift invariant measure on $S$. For 
$s\in S$  and  each $n$   either $s(0)$ lies in a well-defined subword  of $s$ belonging to $\mcw_n$ or in 
a spacer of a subword of $s$ belonging to some  $\mcw_{n+k}$.  By Lemma \ref{catch all} for $\nu$-almost all $x$ and  for all large enough $n$ there is a unique $k$ with $0\le k<q_n$ such that 
$s\rest[-k, q_n-k)\in \mcw_n$.

\begin{definition}\label{def of rn} Let $s\in S$ and suppose that for some $0\le k<q_n, s\rest[-k,q_n-k)\in \mcw_n$.
Define $r_n(s)$ to be the unique $k$ with with this property.   We will call the interval $[-k, q_n-k)$ the \emph{principal $n$-block} of $s$, and $s\rest [-k, q_n-k)$ its \emph{principal $n$-subword}. The sequence of $r_n$'s will be called the \emph{location sequence of $s$}.
\end{definition}
Thus $r_n(s)=k$ is saying that \emph{$s(0)$ is the $k^{th}$ symbol in the principal $n$-subword  of $s$ containing $0$.} We can view the principal $n$-subword of $s$ as being located on an interval $I$ inside the principal $n+1$-subword. Counting from the beginning of the principal $n+1$-subword, the $r_{n+1}(s)$ position is located at the $r_n(s)$ position in $I$. 

\medskip

\begin{remark} \label{interval coherence}  It follows immediately from  the definitions that  if $r_n(s)$ is well-defined and $n\le m$,  the $r_m(s)^{th}$ position of the word 
  occurring in the principal $m$-block of $s$  is in the $r_n(s)^{th}$ position inside the principal $n$-block of $s$. 

\end{remark}

\begin{lemma}\cite{global_structure}\label{specifying elements} Suppose that $s, s'\in S$ and $\la r_n(s):n\ge N\ra=\la r_n(s'):n\ge N\ra$ and for all $n\ge N$, $s$ and $s'$ have the same principal $n$-subwords. Then $s=s'$.
\end{lemma}

Thus an element of $s$ is determined by   knowing  any tail of the  sequence  $\la r_n(s):n\ge N\ra$ together with  a tail of the principal subwords of $s$.

\begin{remark}Here are some consequences of Lemma \ref{specifying elements}: \label{rebuilding}
\begin{enumerate}
\item Given a sequence $\la u_n:M\le n\ra$ with $u_n\in\mcw_n$,   if we specify which occurrence of $u_n$ in $u_{n+1}$ is the principal occurrence, then $\la u_n:M\le n\ra$ determines an $s\in \bk$ completely up to a shift $k$ with $|k|\le q_M$.

\item  
A sequence $\la r_n:N\le n\ra$ and sequence of words $w_n\in \mcw_n$ comes from an infinite word $s\in S$ if  both $r_n$ and $q_n-r_n$  go to infinity and 
that the $r_{n+1}$ position in $w_{n+1}$ is in the $r_n$ position in  a subword of $w_{n+1}$ identical to $w_n$.

\emph{Caveat}: just because $\la r_n:N\le n\ra$ is the location sequence of some $s\in S$ and $\la w_n:N\le n\ra$ is the sequence of principal subwords of some $s'\in S$, it does not follow that there is an $x\in S$ with location sequence $\la r_n:N\le n\ra$ and sequence of subwords $\la w_n:N\le n\ra$. 

\item If
 $x, y\in S$ have the same principal $n$-subwords and $r_n(y)=r_n(x)+1$ for all large enough $n$, then $y=sh(x)$. 

\end{enumerate}
\end{remark}

\subsubsection{A note on inverses of symbolic shifts}\label{note on inverses}
We define operators we label $\rev{}$, and apply them in several contexts
\begin{definition} If $x$ is in $\bk$, define the 
reverse  of $x$ by setting $\rev{x}(k)=x(-k)$.  For $A\subseteq \bk$, define:
\hypertarget{reverse}{
\[\rev{A}=\{\rev{x}:x\in A\}.\]}
If $w$ is a word, let $\rev{w}$ to be the reverse of $w$ sitting on the same interval. Explicitly,  if $w:[a_n, b_n)\to \Sigma$ is the word then $\rev{w}:[a_n,b_n)\to \Sigma$ and $\rev{w}(i)=w((a_n+b_n)-(i+1))$.
If $\mcw$ is a collection of words, $\rev{\mcw}$ is the collection of reverses of the words in $\mcw$.
\end{definition}

If $(\bk, sh)$ is an arbitrary symbolic shift then its inverse is $(\bk, sh^{-1})$.   It will be convenient to have all of the shifts go in the same direction, thus:
\begin{prop}\label{spinning}
The map $\phi$ sending $x$ to $\rev{x}$ is a canonical isomorphism between $(\bk, sh\inv)$ and $(\rev{\bk},sh)$. \end{prop}
 The notation $\bl^{-1}$ stands for the system $(\bl,sh\inv)$ and $\rev{\bl}$ for the system $(\rev{\bl},sh)$.

 \subsection{Generic Points} 
Let $T$ be a measure preserving transformation from $(X,\tau, \mu)$ to $(X,\tau, \mu)$, where $\tau$ is a compact  separable topology, and $\mu$ is a standard measure. Then a point $x\in X$ is \emph{generic} for $T$ if and only if for all $f\in C(X)$, 
 \begin{equation}\lim_{N\to \infty}\left({1\over N}\right)\sum_0^{N-1} f(T^n(x))=\int_X f(x)d\mu(x).\label{def of gen}
 \end{equation}
 
 The Ergodic Theorem tells us that for a given $f$ and ergodic $T$ equation (\ref{def of gen}) holds for a set of $\mu$-measure one. Intersecting over a countable dense set of $f\in C(X)$ gives a set of $\mu$-measure one of generic points. For symbolic systems $\bk\subseteq \Sigma^\poZ$ the generic points are those $x$ such that the $\mu$-measure of all basic open intervals $\la u\ra_0$ is equal to the density of $k$ such that $u$ occurs in $x$ at $k$.

 \subsection{Stationary Codes and $\dbar$-Distance}\label{random label}
In this section we briefly review a standard idea, that of a \emph{stationary code}.  A reader unfamiliar with this material who is interested in the  proofs of the facts cited here should see \cite{Shields} 
 \begin{definition}
Suppose that $\Sigma$ is a countable language. A \emph{code} of length $2N+1$ is a function $\Lambda:\Sigma^{[-N, N]}\to \Sigma$ (where $[-N, N]$ is the interval of integers starting at $-N$ and ending at $N$).

Given a code $\Lambda$, the  \emph{stationary code} determined by $\Lambda$ is the function $\bar{\Lambda}:\Sigma^\poZ \to \Sigma^\poZ$ where, given $s$ 
\[\bar{\Lambda}(s)(k)=\Lambda(s\rest[k-N, k+N]).\]
\end{definition} 
Let $(\Sigma^\poZ, \mcb, \nu, sh)$ be a symbolic system.
Given two codes $\Lambda_0$ and $\Lambda_1$ (not necessarily of the same length),
define $D=\{s\in \Sigma^\poZ:\overline{\Lambda}_0(s)(0)\ne \bar{\Lambda}_1(s)(0)\}$ and $d(\Lambda_0, \Lambda_1)=\nu(D)$. Then $d$ is a semi-metric on the collection of codes.  The following is a consequence of the Borel-Cantelli lemma.

\begin{lemma}\label{cauchy codes} 
Suppose that $\la \Lambda_i:i\in\nn\ra$ is a sequence of codes such that $\sum_{i}d(\Lambda_i, \Lambda_{i+1})<\infty$. Then there is a shift-invariant Borel map $S:\Sigma^\poZ\to \Sigma^\poZ$
such that for $\nu$-almost all $s$, $\lim_{i\to \infty}\overline\Lambda_i(s)=S(s)$.

\end{lemma}

A shift-invariant Borel map $S:\Sigma^\poZ\to \Sigma^\poZ$, determines a factor  $(\Sigma^\poZ, \mcb, \mu, sh)$ of $(\Sigma^\poZ, \mcb, \nu, sh)$ by setting 
$\mu=S^*\nu$. Hence a convergent sequence of stationary codes determines a factor of $(\Sigma^\poZ, \mcb, \nu, sh)$. 

Let $\Lambda_0$ and $\Lambda_1$ be codes. Define $\dbar(\bar\Lambda_0(s), \bar\Lambda_1(s))$  to be
\[
	\lim_{N\to \infty}{|\{k\in [-N, N]:	\bar\Lambda_0(s)(k)\ne 	\bar\Lambda_1(s)(k)\}|\over 2N+1}			
\]
More generally  define the $\dbar$ metric on $\Sigma^{[a,b]}$ by setting 
\[\dbar_{[a,b]}(x,y)={|\{k\in [a, b):	x(k)\ne 	y(k)\}|\over b-a}.\]  For $x, y\in \Sigma^\poZ$, we set
\[\dbar(x,y)=\lim_{N\to \infty}\dbar_{[-N,N]}(x\rest[-N,N], y\rest[-N,N]),\] 
provided this limit exists.

To compute distances between codes we will use the following  application of the Ergodic Theorem.
\begin{lemma}\label{computing code distances} \label{d-bar}{Suppose that $\nu$ is ergodic. Let $\Lambda_0$ and $\Lambda_1$  be codes. Then for almost all $s\in \Sigma^\poZ$:}
\[
d(\Lambda_0, \Lambda_1)=\dbar(\bar\Lambda_0(s), \bar\Lambda_1(s))\]
\end{lemma}
The next proposition is used to study alleged isomorphisms between measure preserving transformations. We again refer the reader to \cite{Shields} for a proof. 
\begin{prop}\label{codes exist}
Suppose that $\bk$ and $\bl$ are symbolic systems and $\phi:\bk\to \bl$ is a factor map. Let $\epsilon>0$. Then there is a code $\Lambda$ such that for almost all $s\in \bk$, 
\begin{equation}\label{dbar approx}
\dbar(\bar{\Lambda}(s), \phi(s))<\epsilon.
\end{equation}\end{prop} 
\bigskip

To show that equation \ref{dbar approx} cannot hold (and hence show that  $\bl$ is not a factor of $\bk$) we will want to view 
$\bar\Lambda(s)$ as limits of  $\Lambda$-images of large blocks of the form $s\rest[a, b]$ with 
$a<0<b$.  There is an ambiguity in doing this: if the code $\Lambda$ has length $2N+1$ it does not make sense to 
apply it to  
$s\rest[k-N,k+N]$ for $k\in [a, a+2N]$ or $k\in [b-2N, b]$.  However if $b-a$ is quite large with respect to $N$, then 
filling in the values for $\Lambda(s\rest[k-N, k+N])$ arbitrarily as $k$ ranges over these initial and final intervals  
makes a negligible difference to the $\dbar$-distances of the result. In particular if 
$\dbar(\bar\Lambda(s), \phi(s))<\epsilon$ then for all large enough $a, b\in \nn$, we have
\[\dbar_{[-a, b]}(\bar{\Lambda}(s\rest[-a,b]),\phi(s)\rest[-a,b])<\epsilon,\]
no matter how we fill in the ambiguous portion.

The general phenomenon of ambiguity or disagreement at the beginning and end of large intervals is referred to by the phrase \emph{end effects}. Because the end effects are usually negligible on large intervals we will often neglect them when computing $\dbar$ distances.

The next proposition is standard:
\begin{prop}\label{weak and dbar} Suppose that $(\Sigma^\poZ, \mcb, \nu, sh)$ is an ergodic symbolic system and $\la T_n:n\in\nn\ra$ is a sequence of functions  from $\Sigma^\poZ\to \Sigma^\poZ$ that commute with the shift.  Then the following are equivalent:
\begin{enumerate}
\item The sequence $\la T_n\ra$ converges to $S$ in the weak topology.
\item $\nu(\{s:T_n(s)(0)\ne S(s)(0)\})\to 0$.
\item For $\nu$-almost all $s, \dbar({T_n(s), S(s)})\to 0$.
\item For some $\nu$-generic $s$, for all $\gamma>0$ we can find an {$N$ for all $n\ge N$}, for all large enough $a,b$, the distance $\dbar({T_n(s)\rest[-a,b), S(s)\rest[-a,b)})<\gamma$.

\end{enumerate}
\end{prop}

\smallskip
\noindent We finish with a remark  that we will use in several places:

\begin{remark}\label{cheating on dbar}
If $w_1$ and $w_2$ are words in a language $\Sigma$ defined on an interval $I$ and $J\subset I$ with ${|J|\over |I|}\ge \delta$, then $\dbar_I(w_1, w_2)\ge \delta\dbar_J(w_1, w_2)$.
 \end{remark}
 
 \subsection{Rotations of the circle}
\label{add or multiply} {Many of the arguments in this paper are based on an understanding of rational approximations to rotations of the circle.  It is usually convenient to adopt additive notation and work on the unit interval $[0,1)$, but this introduces ambiguities.  Fix an $\alpha\in \mathbb R$.  We use the symbol $\mcr_\alpha$ in two ways.  The first way is that 
\[\mcr_\alpha:S^1\to S^1\]
by rotating the circle by $\alpha*2\pi$ radians. The second, equivalent, way is that 
\[\mcr_\alpha:[0,1)\to [0,1)\] and is given by the formula
\[x\mapsto x+\alpha \mod 1.\]
We note in both cases that we are really concerned with $[\alpha] (\mbox{mod} 1)$.}
 
 \subsection{Descriptive Set Theory Basics}\label{dst basics}
 
Let $X$ and $Y$ be Polish spaces and $A\subseteq X, B\subseteq Y$.\footnote{{The ideas in section are just summaries, they are exposited in  \cite{descview} and \cite{kechris}.} }  A function $f:X\to Y$ \emph{reduces} $A$ to $B$ if and only if for all $x\in X$:
 \[x\in A \mbox{ if and only if } f(x)\in B.\]
 For this definition to have content there must be some definability restriction on $f$. The relevant restrictions for this paper are {either that}
 \emph{$f$ is a Borel function} (i.e. the inverse image of an open set is Borel) {or that} \emph{$f$ is a continuous function} (i.e. the inverse image of an open set is open). {The latter is clearly a stronger condition.}  If $B$ is Borel and $f$ is a Borel reduction, 
 then $A$ is clearly Borel. Taking the contrapositive, if $A$ is \emph{not} Borel then $B$ is not. If $A$ is Borel (resp. continously) reducible to $B$ we will write $A\preceq_B B$ (resp. $A\preceq_c B$). Both  $\preceq_B$ and $\preceq_c$ are clearly pre-partial-orderings.\footnote{The reader should be aware that this is a different notion than the notion of a reduction of \emph{equivalence relations}.}
  
 If $\mathcal S$ is a collection of  {pairs $(A,X)$ and $(B,Y)\in \mathcal S$, then $B$ is 
 \emph{$\mathcal S$-complete} for Borel 
reductions  (resp. continuous reductions) if and only if every $(A,X)\in \mathcal S$ is Borel  reducible (resp. continuously reducible) to 
 $(B,Y)$.} Being complete is interpreted as being at least as complicated as each set in $\mathcal S$.
 
 For this to be useful there must be examples of sets that are not Borel. If $X$ is a Polish space and $B\subseteq X$, then $B$ is \emph{analytic} ($\sig^1_1$) if and only if it the continuous image of a Borel subset of a Polish space. This is equivalent to   there being  a Polish space $Y$ and a Borel set $C\subseteq X\times Y$ such that $B$ is the projection to the $X$-axis of $C$.
 
  Correcting a famous mistake of Lebesgue, Suslin proved that there are analytic sets that are not Borel. {It follows immediately that complete analytic sets are not Borel.} This paper uses a canonical example of such a set.
  
   Let $\la \sigma_n:n\in\nn\ra$ be an enumeration of $\nn^{<\nn}$, the finite sequences of natural numbers. Using this enumeration subsets $S\subseteq \nn^{<\nn}$ can be identified with functions $\Chi_S:\nn\to \{0,1\}$.
  
  A \emph{tree} is a set {$\mathcal T\subseteq \nn^{<\nn}$ such that if $\tau\in \mathcal T$ and $\sigma$ is an initial segment of $\tau$, then $\sigma\in \mathcal T$. The set $\{\Chi_{\mathcal T}:{\mathcal T}$ is a tree$\}$ is a closed subset of 
  $\{0,1\}^\nn$, hence a Polish space with the induced topology.  We call the resulting space $\trees$. (In the sequel we will not always distinguish between $\mathcal T$ and $\Chi_{\mathcal T}$.)}
  
  Because the topology on the space of trees is the ``finite information" topology, inherited from the product topology on $\{0,1\}^\nn$, the  following characterizes  continuous maps defined on $\trees$.

\begin{proposition}
Let $Y$ be a topological space and $f:\trees\to Y$. Then $f$ is continuous if and only if for all open $O\subseteq Y$ and all $\mct$ with $f(\mct)\in O$ there is an $M\in \nn$ for all $\mct'\in \trees$:
\begin{quotation}
if $\mct\cap \{\sigma_n:n\le M\}=\mct'\cap\{\sigma_n:n\le M\}$, then $f(\mct')\in O$.
\end{quotation}
\end{proposition}
  
  An infinite branch through $T$ is a function $f:\nn\to \nn$ such that for all $n\in \nn, f\rest \{0, 1, 2, \dots n-1\}\in T$. A tree $T$ is \emph{ill-founded} if and only if it has an infinite branch.
 
 The following theorem is classical; proofs can be found in \cite{kechris}, \cite{marker}. 
\begin{fact}\label{comp anal}
Let $\trees$ be the space of trees. Then:
\begin{enumerate}
\item The collection of ill-founded trees is a complete analytic subset of $\trees$.
\item The collection of trees that have at least two distinct infinite branches is a complete analytic subset of $\trees$.
\end{enumerate}
\end{fact} 
The main results of this paper (Theorem \ref{lele} and Corollary \ref{cor3}) are proved by reducing the sets mentioned in Theorem \ref{comp anal} to conjugate pairs of diffeomorphisms and concluding that the sets of conjugate pairs is complete analytic--so not Borel.

\section{Odometer and Circular Systems}\label{odo circ review}

Two types of symbolic shifts play central roles for the proofs of the main theorem, the \emph{odometer based} and the  \emph{circular} systems. Most of the material in this section appears in \cite{global_structure} in  more detail and is reviewed here without proof.
\subsection{Odometer Based Systems}
We now define the class of \emph{Odometer Based Systems}. In a sequel to this paper (\cite{part_IVa}), we prove that these are exactly the finite entropy transformations that have non-trivial odometer factors. We recall the definition of an odometer transformation.  Let $\la k_n:n\in\nn\ra$ be a sequence of natural numbers greater than or equal to 2. Let 
\[O=\prod_{n=0}^\infty \poZ/k_n\poZ\]
be the $\la k_n\ra$-adic integers. 
Then $O$ naturally has a compact abelian  group structure and hence carries a Haar measure $\mu$. 
{The set $O$ becomes a measure preserving system 
$\mco$ by defining $T:O\to O$ to be addition by 1 in the $\la k_n\ra$-adic integers. Concretely, this is the map that ``adds one to $\poZ/k_0\poZ$ and carries right".} Then $T$ is  an invertible transformation that preserves the Haar measure $\mu$ on $\mco$. Let $K_n=k_0*k_1*k_2\dots k_{n-1}$.

The following results are standard:
\begin{lemma} \label{odometer basics}Let $\mco$ be an odometer system. Then:
\begin{enumerate}
\item $\mco$ is ergodic.
\item The map $x\mapsto -x$ is an isomorphism between $(O, \mcb, \mu, T)$ and $(O, \mcb, \mu, T^{-1})$.
\item {Odometer maps are transformations with discrete spectrum and the eigenvalues of the associated linear operator are the $K_n^{th}$ roots of unity ($n>0$).}
\end{enumerate}
\end{lemma}

Any natural number {$a<K_j$} can be uniquely written as:
\[a=a_0+a_1k_0+a_2(k_0k_1)+ \dots +a_j(k_0k_1k_2\dots k_{j-1})\] 
for some sequence of natural numbers $a_0, a_1, \dots a_j$ with $0\le a_j<k_j$.

\begin{lemma}\label{specifying elements of the odometer} Suppose that 
$\la r_n:n\in\nn\ra$ is a sequence of natural numbers with $0\le r_n<k_0k_1\dots k_{n}$ and 
{$r_n\equiv r_{n+1} \mod(k_0k_1\dots k_{n})$}.  Then there is a unique element $x\in O$ such that $r_n=x(0)+x(1)k_0+\dots x(n)(k_0k_1\dots k_{n-1})$ for each $n$.

\end{lemma}
We now define the collection of symbolic systems that have odometer maps as their timing mechanism. This timing mechanism can be used to parse typical elements of the symbolic system.

\begin{definition}Let $\la \mcw_n:n\in\nn\ra$ be a uniquely readable construction sequence with the properties that {$\mcw_0=\Sigma$} and for all $n, \mcw_{n+1}\subseteq (\mcw_n)^{k_n}$ for some $k_n$. The  associated symbolic system will be called an \emph{odometer based system}.
\end{definition}
Thus odometer
based systems are those built from  construction sequences $\la \mcw_n:n\in\nn\ra$ such that the words in 
$\mcw_{n+1}$ are concatenations of words in $\mcw_n$ of a fixed length $k_n$. The words in $\mcw_{n}$ all have length $K_n$ and the words $u_i$  in  equation \ref{words and spacers} are all the empty words. 

Equivalently, an odometer based transformation is one that can be built by a cut-and-stack construction using no spacers. An easy consequence of the definition is that for odometer based systems, for all $s\in S$ and for all 
$n\in \nn$,
 $r_n(s)$ exists.\footnote{$S$ is defined in  Definition \ref{def of S}.}

 The next lemma justifies the terminology.
 \begin{lemma}\label{odometer factor}
 Let $\bk$ be an odometer based system with each $\mcw_{n+1}\subseteq (\mcw_n)^{k_n}$. Then there is a canonical factor map
 \begin{equation*}
 \pi:S\to \mco
 \end{equation*}
 where $\mco$ is the odometer system determined by $\la k_n:n\in\nn\ra$.
 \end{lemma}
\pf For each $s\in S$, for all $n, r_n(s)$ is defined and both $r_n$ and $k_n-r_n$ go to infinity. By 
Lemma \ref{specifying elements of the odometer}, the sequence $\la r_n(s):n\in\nn\ra$ defines a unique element $\pi(s)$ in 
$\mco$. It is easily checked that $\pi$ intertwines $sh$ and $T$.\qed

Heuristically, the odometer transformation $\mco$ parses the sequences $s$ in $S\subseteq \bk$ by indicating where 
the words constituting $s$ begin and end. Shifting $s$ by one unit shifts this parsing by one. We can understand 
elements of $S$ as being an element of the odometer with words in $\mcw_n$ filled in inductively.

The following remark is useful when studying the canonical factor of  the inverse of an odometer based system.
\begin{remark}
If $\pi:\bl\to\mco$ is the canonical factor map, then the function $\pi:L\to O$ is  also factor map from $(\bl,sh\inv)$ to $\mco\inv$ (i.e. $O$ with the operation ``$-1$").  If $\la\mcw_n:n\in\nn\ra$ is the construction sequence for $\bl$, then $\la \rev{\mcw_n}:n\in\nn\ra$ is a construction sequence for $\rev{\bl}$. If $\phi:\bl\inv \to \rev{\bl}$ is the canonical isomorphism given by Proposition \ref{spinning}, then Lemma  \ref{odometer basics} tells us that the projection of $\phi$ to a map $\phi^\pi:\mco\to \mco$ is given by  $x\mapsto -x$.
\end{remark}

The following is proved in \cite{global_structure}:
\begin{prop} \label{kiss your S goodbye}
Let $\bk$ be an odometer based system and suppose that $\nu$ is a shift invariant measure. Then 
$\nu$ concentrates on $S$.
\end{prop}

\subsection{Circular Systems} \label{circular systems 1} We now define \emph{circular systems}.  In \cite{prequel} it is shown that the strongly uniform circular systems give symbolic  characterizations of certain smooth diffeomorphisms defined by the Anosov-Katok method of conjugacies.

 These systems are called \emph{circular}  because they are related to the behavior of rotations by a convergent sequence of rationals $\alpha_n=p_n/q_n$.  
 The rational rotation by $p/q$ permutes the $1/q$ intervals of the circle cyclically in a manner that the interval $[i/q, (i+1)/q)$ occurs in position 
 $j_i=_{def}p^{-1}i$ (mod $q$).\footnote{We assume that $p$ and $q$ are relatively prime and the exponent $-1$ indicates the multiplicative inverse  modulo $q$.} 
 The operation 
 $\mcc$ which we are about to describe models the relationship between rotations by $p/q$ and $p'/q'$ when $p'/q'$ is very close to $p/q$.

Let $k, l, p, q$ be positive natural numbers with $p<q$ relatively prime. {For $0\le i<q$, setting 
$j_i\equiv_q(p)^{-1}i \label{j sub i}$
with $j_i<q$. It is easy to verify that:
\begin{equation}\label{reverse numerology} 
q-j_i=j_{q-i}. 
\end{equation}
For notational convenience later we set $j_q=q$.}

Let $\Sigma$ be a non-empty set such that neither $b$ nor $e$ belongs to 
$\Sigma$ and $w_0, \dots w_{k-1}$ be words in $\Sigma\cup \{b,e\}$. 
Define:\footnote{We use $\prod$  and powers for repeated concatenation of
 words.}

\begin{equation}\mcc(w_0,w_1,w_2,\dots w_{k-1})=\prod_{i=0}^{q-1}\prod_{j=0}^{k-1}(b^{q-j_i}w_j^{l-1}e^{j_i}). \label{definition of C}
\end{equation} 

{We note that the product symbol $\Pi$ is repeated concatenation as is the exponent.  If $w$ is a word then $w^0$ is the empty string, $w^1=w$, $w^2=ww$ and so forth. The formula in equation \ref{definition of C} is a concatenation  of $q$ words, each of which is itself, a concatenation of $k$ words.  The words inside the parenthesis in equation \ref{definition of C} start with $q-j_i$ $b$'s, followed by concatenating $l-1$ many $w$'s, followed by concatenating $j_i$ many $e$'s.  Written with parenthesis:
\begin{equation}\mcc(w_0,w_1,w_2,\dots w_{k-1})=\prod_{i=0}^{q-1}\left(\prod_{j=0}^{k-1}((b^{q-j_i})(w_j^{l-1})(e^{j_i}))\right). 
\end{equation} 
Informally, the $i^{th}$ term, $\prod_{j=0}^{k-1}(b^{q-j_i}w_j^{l-1}e^{j_i})$  can be written as a block of $q-j_i$ $b$'s followed by $w_0$ concatenated with itself $l-1$ times, followed by a block of $j_i$ many $e$'s, followed by a block of $q-j_i$ $b$'s followed by $w_1$ concatenated with itself $l-1$ times followed by a block of $j_i$ $e$'s and so forth, ending with a block of $w_{k-1}$ repeated $l-1$ times followed by $e$ repeated $j_i$ many times:
\footnotesize{\[(bbb\dots) (w_0w_0\dots )(ee\ .. e)(bb\dots b) (w_1w_1\dots w_1)(ee\dots e)\dots \dots(bb\dots b)(w_{k-1}w_{k-1}w_{k-1}\dots w_{k-1})(ee\dots e)\]}}

{\begin{remark}\label{word length} \
\begin{itemize}
\item Suppose that each $w_i$ has length $q$, then the length of  $\mcc(w_0, w_1, \dots w_{k-1})$ is $klq^2$. 
\item For each occurrence of an $e$ in $\mcc(w_0, \dots w_{k-1})$  there is an occurrence of $b$ to the left of it.
\item {Suppose that $n<m$ and $b$ occurs at $n$ and $e$ occurs at $m$ and neither occurrence is in a $w_i$. Then there must be some $w_i$ occurring between $n$ and $m$.}
\item Words constructed with $\mcc$ are uniquely readable. 
\end{itemize}
\end{remark}}

The $\mcc$ operation is used to build a collection of symbolic shifts. \emph{Circular systems} will be defined using a sequence of natural number parameters $ k_n$ and $l_n$ that is fundamental to the version of the Anosov-Katok construction  presented in \cite{katoksbook}.

 Fix an arbitrary sequence of positive natural numbers $\la k_n:n\in\nn\ra$. Let {$\la l_n:n\in\nn\ra$} be an increasing sequence of natural numbers such that 
\begin{numreq}\label{l_n summability}
$l_0>20$ and $\sum_{k\ge n} 1/l_k< 1/l_{n-1}$. \end{numreq}
 From the $k_n$ and $l_n$ we define  sequences of numbers: $\la p_n, q_n, \alpha_n:n\in\nn\ra$. 
 Begin by letting $p_0=0$ and $q_0=1$  and inductively set 
\begin{eqnarray}\label{qn}
\qnpo=\kn l_n\qn^2
\end{eqnarray}
{(thus $q_1=k_0l_0$)} and take 
\begin{eqnarray}p_{n+1}=p_nq_nk_nl_n+1.\label{pn}
\end{eqnarray}
Then clearly $p_{n+1}$ is relatively prime to $q_{n+1}$.\footnote{{$p_n$ and $q_n$ being relatively prime for $n\ge 1$, allows us to define the integer $j_i$ in equation \ref{j sub i} . For $q_0=1$, $\poZ/q_0\poZ$ has one element, $[0]$, so we set $p_0\inv=p_0=0$.}}

Setting $\alpha_n=p_n/q_n$,  then it is easy to check that there is an irrational $\alpha$ such that the sequence $\alpha_n$ converges rapidly to $\alpha$.

\begin{definition}
A sequence of integers $\la k_n, l_n:n\in\nn\ra\ra$ such that $k_n\ge 2$, $\sum 1/l_n<\infty$  will be called a \hypertarget{circ coef}{\emph{circular coefficient sequence}}.
\end{definition}

 Let $\Sigma$ be a non-empty finite or countable alphabet. Build collections of words $\mcw_n$ in $\Sigma\cup \{b, e\}$ by induction as follows:
 \begin{itemize}
 \item Fix a  \hyperlink{circ coef}{circular coefficient sequence} $\la k_n, l_n:n\in\nn\ra\ra$.
 \item Set $\mcw_0=\Sigma\cup \{b, e\}$.  
 \item Having built $\mcw_n$ choose a set $P_{n+1}\subseteq (\mcw_{n})^{k_n}$ and  form $\mcw_{n+1}$ by taking all words of the form $\mcc(w_0,w_1\dots w_{k_n-1})$  with $(w_0, \dots w_{k_n-1})\in P_{n+1}$.\footnote{Passing from $\mcw_n$ to $\mcw_{n+1}$, use $\mcc$ with parameters $k=k_n, l=l_n, p=p_n$ and $q=q_n$ and take $j_i=(p_n)^{-1}i$ modulo $q_n$.  By Remark \ref{word length}, the length of each of the words in $\mcw_{n+1}$ is $q_{n+1}$.} 
 \end{itemize}
We will call the elements of $P_{n+1}$ \hypertarget{pwords}{\emph{prewords}.} The $\mcc$ operator automatically creates uniquely readable words, however we will need a  stronger unique readability assumption for our definition of circular systems.

 \medskip 
 \hypertarget{strong unique readability}{\bfni{Strong Unique Readability Assumption:}}
  Let $n\in \nn$, and view $\mcw_n$ as a collection $\Lambda_n$ of letters. Then each element of $P_{n+1}$ can be viewed as a word  with letters in   $\Lambda_n$. In the alphabet $\Lambda_n$,  each $w\in P_{n+1}$ is uniquely readable.  
  \begin{definition}A construction sequence $\la \mcw_n:n\in\nn\ra$  will be called \emph{circular} if 
  it is built in this manner using the $\mcc$-operators, a circular coefficient sequence and each 
  $P_{n+1}$ satisfies the strong unique readability assumption.   \end{definition}

\begin{definition}\label{circular definition}
A symbolic shift $\bk$ built from a  circular construction sequence  will be called a \emph{circular system}. 
\end{definition}
\bfni{Notation:} we will often write $\bk^c$ and $\la \mcw_n^c:n\in\nn\ra$ to emphasize that we are building circular systems and circular construction sequences. Circular words will often be denoted $w^c$ for emphasis.
\smallskip

\begin{definition}
Suppose that $w=\mcc(w_0,w_1,\dots w_{k-1})$.  Then $w$ consists of blocks of $w_i$ repeated $l-1$ times, together with some $b$'s and $e$'s that are not in the $w_i$'s. The \emph{interior} of $w$ is the portion of $w$ in the $w_i$'s.
{The remainder of  $w$ consists of blocks of the form $b^{q-j_i}$ and $e^{j_i}$. We call this portion the  \emph{boundary} of $w$.} 

In a block of the form $w_j^{l-1}$ the first and last occurrences  of $w_j$ will be called the  \emph{boundary} occurrences of the block $w_j^{l-1}$. The other occurrences will be the \emph{interior} occurrences. 
\end{definition}
{While the boundary consists of sections of $w$ made up of $b$'s and $e$'s, not all $b$'s and 
$e$'s occurring in $w$ are in the boundary, as they may be part of a power $w_i^{l-1}$.}

 The boundary of $w$ constitutes a small portion of the word:

\begin{lemma}\label{stabilization of names 1} Suppose that $w=\mcc(w_0, w_1,\dots ,w_{k-1})$ and each $w_i$ has length $q$.  Then the proportion of the word $w$ that belongs to its boundary is $1/l$. Moreover the proportion of the word that is within $q$ letters of boundary of $w$ is $3/l$. 
\end{lemma}

\pf The length of $w$ is $klq^2$.  The boundary portions are $q*k*q$ long. The number of letters within $q$ letters of the boundary is $q*k*3*q$.\qed

\begin{remark}
Let $v_0, \dots v_{k-1}$ and $w_0,\dots w_{k-1}$ be sequences of words of length $q$. The boundary portions of $\mcc(v_0,\dots v_{k-1})$ and 
$\mcc(w_0,\dots w_{k-1})$ occur in the same positions and by Lemma \ref{stabilization of names 1} have proportion $1/l$ of the length. Since all of the $v_i$'s and $w_i$'s have the same length and the same multiplicity in the circular words we see:
\begin{eqnarray*}
\dbar(\mcc(v_0,\dots v_{k-1}),\mcc(w_0,\dots w_{k-1}))\ge& \\ &(1-1/l)\dbar(v_0v_1v_2\dots v_{k-1}, w_0w_1\dots w_{k-1})
\end{eqnarray*}
where $v_0v_1v_2\dots v_{k-1}$ and  $w_0w_1\dots w_{k-1}$ are the concatenations of the various words.\footnote{Equality holds, a fact we won't use.}
\end{remark}
For proofs of the next lemma see \cite{prequel} (Lemma 20) and  \cite{global_structure}.

\begin{lemma}\label{dealing with S} Let $\bk^c$ be a circular system and
 $\nu$ be a shift-invariant measure on $\bk^c$. Then the following are equivalent: 
	{\begin{enumerate}
	\item $\nu$ has no atoms.   
	\item $\nu$ concentrates on  the collection of $s\in \bk^c$ such that $\{i:s(i)\notin \{b, e\}\}$ is 
unbounded in both $\poZ^-$ and $\poZ^+$.
	\item  $\nu$ concentrates on $S$.
	\end{enumerate}}
	\noindent If $\bk^c$ is a uniform circular system (Example \ref{definition of uniform}), then there is a unique invariant measure concentrating on $S$. 
	
\end{lemma}
Moreover there are only two ergodic invariant measures with atoms:  the one concentrating on the constant sequence $\vec{b}$ and the one concentrating on $\vec{e}$.

\begin{remark}If $\bk^c$ is circular and $s\in \bk^c$ has a principal $n$-subword and $m>n$, then $s$ has a principal $m$-subword.
\end{remark}
\subsection{An Explicit Description of {$\rev{\bk^c}$.}}
The symbolic system ${\bk^c}$ is built by an operation $\mcc$ applied to collections of words. The system {$\rev{\bk^c}$} is built by a similar operation applied to the reverse collections of words. In analogy to equation \ref{definition of C},  we define $\mcc^r$ as follows:
\begin{definition}Suppose that $w_0, w_1, \dots w_{k-1}$ are words in a language $\Sigma$. Given coefficients $p, q, k, l$ with $p$ and $q$ relatively prime, let $j_i\equiv_q (p^{-1})i$ with $0\le j_i<q$. Define
\begin{equation}\mcc^r(w_0,w_1,w_2,\dots w_{k-1})=
\prod_{i=0}^{q-1}\prod_{j=0}^{k-1}(e^{q-j_{i+1}}({w_{k-j-1}}^{l-1})b^{j_{i+1}}). \label{definition of Cr}
\end{equation} 
\end{definition}

From equation \ref{definition of C}, a $w\in \mcw^c_{n+1}$  is of the form $\mcc(w_0, \dots w_{k_n-1})$:
\begin{equation}\label{C again}
w=\prod_{i=0}^{q-1}\prod_{j=0}^{k-1}(b^{q-j_i}w_j^{l-1}e^{j_i})
\end{equation}
where $q=q_n,k=k_n,  l=l_n$ and  $j_i\equiv_{q_n}(p_n)^{-1}i$ with $0\le j_i<q_n$. By examining this formula we see that 
\begin{equation*}
\rev{w}=\prod_{i=1}^{q}\prod_{j=1}^ke^{j_{q-i}}\rev{w_{k-j}}^{l-1}b^{q-j_{q-i}}.
\end{equation*}
Applying the identity in formula \ref{reverse numerology}, we see that this can be rewritten as\footnote{Recall that we take $j_q=q$, so $q-j_q=0$.}
\begin{equation}\label{reverse mcc}
\rev{w}=\prod_{i=1}^{q}\prod_{j=1}^k(e^{q-j_i}\rev{w_{k-j}}^{l-1}b^{j_i}).
\end{equation}
Thus
\begin{equation} \label{reversing words}
\rev{w}=\mcc^r(\rev{w_0},\rev{w_1},\dots \rev{w_{k-1}}).
\end{equation}
In particular if $\la \mcw^c_n:n\in\nn\ra$ is a construction sequence of a circular system $\bk^c$, then 
$\rev{\mcw^c_{n+1}}$ is the collection:
\begin{equation*}
\{\mcc^r(\rev{w_0}, \rev{w_1}, \dots, \rev{w_{k_n-1}}):w_0w_1 \dots w_{k_n-1}\in P_n\} 
\end{equation*}
and {$\la \rev{\mcw_n^c}:n\in \nn\ra$ is}  a construction sequence for $\rev{\bk^c}$.

\subsection{Understanding the Words}\label{understanding the words}
The words used to form circular transformations have quite specific combinatorial properties. 
Fix a sequence $\la \mcw^c_n:n\in \nn\ra$ defining a circular system.
Each $u\in \mcw^c_{n+1}$ has three \emph{subscales}: 
\begin{enumerate}
\item[]{\bf Subscale 0,} the scale of the individual powers of $w\in \mcw^c_n$ of the form $w^{l-1}$; We call each such occurrence of a $w^{l-1}$ a \emph{0-subsection}.
\item[]{\bf Subscale 1,} the scale of each term in the product $\prod_{j=0}^{k-1}(b^{q-j_i}w_j^{l-1}e^{j_i})$ that has the form $(b^{q-j_i}w_j^{l-1}e^{j_i})$; We call these terms \emph{1-subsections}.
\item[]{\bf Subscale 2,} the scale of each term of $\prod_{i=0}^{q-1}\left(\prod_{j=0}^{k-1}(b^{q-j_i}w_j^{l-1}e^{j_i})\right)$ that has the form 
$\prod_{j=0}^{k-1}(b^{q-j_i}w_j^{l-1}e^{j_i})$; We call these terms \emph{2-subsections}.
\end{enumerate}
\begin{center}
{\bf Summary}
\end{center}

\begin{center}
  \begin{tabular}{| l | c |}
    \hline

  {\bf Whole Word:}  &$\prod_{i=0}^{q-1}\prod_{j=0}^{k-1}(b^{q-j_i}w_j^{l-1}e^{j_i})$\\ \hline

{\bf 2-subsection:} &$ \prod_{j=0}^{k-1}(b^{q-j_i}w_j^{l-1}e^{j_i})$ \\ \hline

{\bf 1-subsection:} &$(b^{q-j_i}w_j^{l-1}e^{j_i})$\\ \hline

{\bf 0-subsection:} & $w_j^{l-1}$\\
    \hline
  \end{tabular}
\end{center}
For $m\le n$, we will discuss ``$m$-subwords"  of a word $w$. These will be subwords that lie in $\mcw^c_{m}$, the ${m}^{th}$ stage of the construction sequence. We will use ``$m$-block" to mean the location of the $m$-subword. \bigskip

\begin{lemma}\label{gap calculation}Let $w=\mcc(w_0, \dots w_{k_n-1})$ 
for some $n$ and $q=q_n, k=k_n, l=l_n$. View $w:\{0, 1, 2\dots ,klq^2-1\}\to \Sigma\cup\{b, e\}$. 
\begin{enumerate}
\item If $m_0$ and $m_1$ are such that $w(m_0)$ and $w(m_1)$ are at the beginning of $n$-subwords in the same 2-subsection, then $m_0\equiv_qm_1$.
\item If $m_0$ and $m_1$ are such that $w(m_0)$ is the beginning of an $n$-subword  occurring in a $2$-subsection 
$\prod_{j=0}^{k-1}(b^{q-j_i}w_j^{l-1}e^{j_i})$ and $w(m_1)$ 
is the beginning of an $n$-subword  occurring in the next 2-subsection $\prod_{j=0}^{k-1}(b^{q-j_{i+1}}w_j^{l-1}e^{j_{i+1}})$ then $m_1-m_0\equiv_q -j_1$.
\end{enumerate}
\end{lemma}
\pf To see the first point, the indices of the beginnings of $n$-subwords in the same $2$-subsection  differ by multiples of $q$ coming from powers of a $w_j$ and intervals of $w$ of the form $b^{q-j_i}e^{j_i}$.

To see the second point, let $u$ and $v$ be consecutive $2$-subsections. In view of the first point it suffices to consider the last $n$-subword of $u$ and the first $n$-subword of $v$.  These sit on either side of an interval of the form $e^{j_i}b^{q-j_{i+1}}$. Since
$j_i+q-j_{i+1}\equiv_q (p)^{-1}i-p^{-1}(i+1)\equiv_q-p^{-1}\equiv_q-j_1$, we see that $m_0-m_1=q+j_i+q-j_{i+1}\equiv_q-j_1$.
\qed

Assume that $u\in \mcw^c_{n+1}$ and $v\in \mcw^c_{n+1}\cup\rev{\mcw^c_{n+1}}$ and $v$ is shifted with respect to $u$.
 On the overlap of $u$ and $v$, the 2-subsections of $u$ split each 2-subsection of $v$ into either one or two pieces. Since the 2-subsections all have the same length, the number of pieces in the splitting  and the size of each piece is constant across the overlap except perhaps at the two ends of the overlap. If $u$ splits a 2-subsection of $v$ into two pieces, then we call the leftmost  piece of the  pair  the even piece and the rightmost the odd piece.

If $v$ is shifted only slightly, it can happen that either the even piece or the odd piece does not contain {even one entire} $1$-subsection. In this case we will say that the split is \emph{trivial on the left} or \emph{trivial on the right}

\begin{lemma}\label{numerology lemma} {Assume that $u\in \mcw^c_{n+1}$ and $v\in \mcw^c_{n+1}\cup\rev{\mcw^c_{n+1}}$ and $v$ is shifted with respect to $u$.}
Suppose that the $2$-subsections of $u$ divide the $2$-subsections of $v$ into two non-trivial pieces. Then
\begin{enumerate}
\item the boundary portion of $u$ occurring between each consecutive pair of 2-subsections of 
$u$ completely overlaps at most one $n$-subword of $v$ 
 \item  there are two numbers $s$ and $t$ such that the positions of the $0$-subsections of $v$ in even pieces are  shifted relative to the  $0$-subsections of $u$ by $s$ and the positions of the  $0$-subsections of $v$ in  odd pieces are shifted relative to the $0$-subsections of $u$ by $t$.  Moreover $s\equiv_q t -j_1$.
\end{enumerate}
\end{lemma}

\pf This follows easily from Lemma \ref{gap calculation}\qed
In the case where  the split is trivial  Lemma \ref{numerology lemma} holds with just one coefficient, $s$ or $t$.
A special case {of} Lemma \ref{numerology lemma} that we will use is:

\begin{lemma}\label{weak numer}{Assume that $u\in \mcw^c_{n+1}$ and $v\in \mcw^c_{n+1}\cup\rev{\mcw^c_{n+1}}$ and $v$ is shifted with respect to $u$.}
Suppose that the $2$-subsections of $u$ divide the $2$-subsections of $v$ into two pieces and 
 that for some 
 occurrence of a $n$-subword 

 in an even (resp. odd) piece is lined up with an occurrence of some  $n$-subword 
 in $u$.  Then 
 every occurrence of a $n$-subword 
 in an even (resp. odd) piece of $v$ is either:
 \begin{enumerate}
 \item[a.)] lined up with some $n$-subword of $u$ or
 \item[b.)]lined up with a section of a $2$-subsection that has the form $e^{j_i}b^{q-j_i}$.  
 \end{enumerate}
 Moreover, no $n$-subword in an odd (resp. even) piece of $v$ is lined up with a $n$-subword in $u$.
 \end{lemma}

\subsection{Full Measure Sets for Circular Systems}\label{full measure for ccs}

Fix a sequence $\la \varepsilon_n:n\in\nn\ra$ such that 
\begin{numreq}\label{varepsilons}
$\la \varepsilon_n:n\in\nn\ra$ is a decreasing sequence of numbers in $\zoo$ such that 
$6\sum_{n>N}\varepsilon_n<\varepsilon_N$. 
\end{numreq}

{From Lemma \ref{stabilization of names 1}, the boundary of a word $w_n\in \mcw_n$ has proportion $1/l_n$.  Hence Numerical Requirement 2 implies that for all choices $\la w_n:n\in \nn\ra$ with $w_n\in \mcw_n$, the sum of the proportion of the boundary sections of $w_n$ is finite.}

\begin{definition} Let:
\begin{enumerate}
\item  $E_n$ be the collection of $s\in S$ such that either $s$ does not have a principal $n$-block or $s(0)$ is in the boundary of the principal $n$-block of $s$, 
\item  $E^0_n=\{s:s(0)$ is in the first or last $\varepsilon_nl_n$ copies of $w$ in a power of the form $w^{l_n-1}$ where $w\in\mcw^c_n\}$, 
\item $E^1_n=\{s:s(0)$ is in the first or last $\varepsilon_nk_n$  1-subsections of the 2-subsection in which $s(0)$ is located.$\}$, 
\item $E^2_n=\{s:s(0)$ is in the first or last $\varepsilon_nq_n$ 2-subsections of its principal $n+1$-block$\}$.
\end{enumerate}
\end{definition}

\begin{lemma}\label{bc1} {Assume numerical requirements 1 and 2.}
   Let $\nu$ be a shift-invariant measure on $S\subseteq \bk^c$, where $\bk^c$ 
   is a circular system. Then:
\begin{enumerate} 
\item \[\sum_n\nu(E_n)<\infty.\] 
\smallskip

For $i= 0, 1, 2$:

\item
\[\sum_n\nu(E^i_n)<\infty.\]
\end{enumerate}
\end{lemma}
\pf {By the Ergodic Theorem we have $\nu(E_n)<1/l_n$,  and for $i=0,1, 2, \nu(E_n^i)<\varepsilon_n$. The result then follows by the summability of $1/l_n$ and $1/\varepsilon_n$}\qed
In particular we see:
\begin{corollary}\label{bc2}
For $\nu$-almost all $s$ there is an $N=N(s)$ such that for all $n>N$, 
\begin{enumerate}
\item $s(0)$ is in the interior of its principal $n$-block,
\item For $i=0, 1, 2,$ $s\notin E^i_n$. 

In particular, for almost all $s$ and all large enough $n$:
\item if $s\rest [-r_n(s), -r_n(s)+q_n)=w$, then 
    \begin{equation*}
    s\rest[-r_n(s)-q_n, -r_n(s))=s\rest [-r_n(s)+q_n, -r_n+2q_n)=w.
    \end{equation*}

\item $s(0)$ is not in a string of the form $w_0^{l_n-1}$ or $w_{k_n-1}^{l_n-1}$.
\end{enumerate}
\end{corollary}
\pf {Apply the Borel-Cantelli Lemma using the previous lemma.}\qed

{The elements 
$s$ of $S$ such that   some shift $sh^k(s)$ fails one of the conclusions 1.)-4.) of Corollary \ref{bc2}  form a measure zero set.} Consequently we work on those elements of $S$ whose whole orbit satisfies the conclusions of Corollary \ref{bc2}. {Note however that for $t=sh^k(s)$, the $N(t)$ in Corollary \ref{bc2}, depends on  $k$.}

\begin{definition}\label{mature}
We will call $n$ \emph{mature} for $s$ (or say that \emph{$s$ is mature at stage $n$}) iff $n$ is so large  that 
 $s\notin E_m \cup \bigcup_{0\le i\le 2}E^i_m$
 for all $m\ge n$.
\end{definition}

If $s$ is mature at stage $n$ then {$s$} is mature at stage $n+1$. Moreover, if $sh^k(s)$ has the same principal $n$-block as $s$ does then $sh^k(s)$ is mature if and only if $s(k)$ is not in {the boundary portion of the principal $n$-block.}

\begin{numreq} \label{varepsilons and q's} 
The following hold:
\begin{eqnarray*}
\varepsilon_nk_n&\to \infty\\
\varepsilon_nl_n&\to \infty\\
\varepsilon_nq_n&\to\infty.
\end{eqnarray*}
\end{numreq}

{\begin{definition}\label{def bound} We will use the symbol $\boundary_n$ in multiple equivalent ways. If $s\in S$ or $s\in \mcw^c_m$  define $\boundary_n=\boundary_n(s)\subseteq \poZ$ to be the collection of $i\in\poZ$ such that $sh^i(s)(0)$ is in the boundary portion of an $n$-subword of $s$.  In the spatial context define $s\in \boundary_{n}\subseteq \bk^c$ by putting $s\in \boundary_n$ if  $s(0)$ is the boundary of an $n$-subword of $s$. 
\end{definition}}
\noindent For $s\in S$
\[\boundary_n(s)\subseteq\bigcup\{[l,l+q_n):s\rest[l,l+q_n)\in \mcw^c_n\}.\] 
The relationship between $\boundary_n(s)\subseteq \poZ$ and $\boundary_n\subseteq \bk^c$ is that for $s\in \bk^c$:
\[i\in \boundary_n(s)\subseteq \poZ\mbox{ iff }sh^i(s)\in\boundary_n\subseteq\bk^c.\]

The next lemma says 
 that if $s$ is mature at stage $n$, then we can detect  locally  those $i$ for which the $i$-shifts of $s$ 
 are mature.

\begin{lemma}\label{getting old}
Suppose that $s\in S$, $n$ is mature for $s$ and $n<m$.

\begin{enumerate}
 
    \item Assume the first three numerical requirements. Suppose that 
    $i\in [-r_m(s), q_m-r_m(s))$. Then $n$ is mature for $sh^i(s)$ iff 
        \begin{enumerate}
        \item $i\notin \bigcup_{n\le k\le m} \boundary_k(s)$ and
        \item $sh^i(s)\notin \bigcup_{n\le k< m}(E^0_k\cup E^1_k\cup E^2_k)$.
        \end{enumerate}

   \item\label{simple revision} {For all but at most $(\sum_{n< k\le m}1/l_k) +   (\sum_{n\le k< m}
    6\varepsilon_k)$ proportion  of the $i\in [-r_m(s), q_m-r_m(s))$, 
    the point $sh^i(s)$ is mature for $n$.}
\end{enumerate}
Hence by numerical requirement \ref{varepsilons},  the proportion of $i\in [-r_m(s), q_m-r_m(s))$ for which the $i$-shift of $s$ is not  mature for $n$ is less than $1/l_{n-1}+\varepsilon_{n-1}$.
\end{lemma}
\pf The first item is immediate from the definition of \emph{mature}.  
For the second item, first note that 
\[\bigcup_{n\le k\le m}\boundary_k(s)\cup  
\bigcup_{n\le k< m}(E^0_k\cup E^1_k\cup E^2_k)
=\boundary_m(s)\cup \bigcup_{n\le k< m}\left(\boundary_k(s)\cup E^0_k\cup E^1_k\cup E^2_k\right) .\]
Let $I=[-r_m(s), q_m-r_m(s))$. Since $\boundary_m$ has proportion $1/l_m$ of $I$, it suffices to show that for a fixed $k\in[n,m)$, the proportion of $i\in I$ such that $sh^i(s)\in \boundary_k\cup E^0_k\cup E^1_k\cup E^2_k$ is less than $1/l_k+6\epsilon_k$. 

There are at most $q_m/q_k$ $k$-words appearing in $s\rest I$.  There are at most $1/l_k$ many $i$ in the boundary of  each of these $k$-words.  So total number of $i$ in $\boundary_k(s)\cap I$ is less than or equal to 
$({q_m\over q_k})(q_k/l_k)$, hence has proportion less than or equal to $1/l_k$ of $I$. 

Similarly for $j=0,1,2$ the number of $i$ with $sh^i(s)\in E^j_k$ and $i$ is in the 
block corresponding to a $k$-subword of $s\rest I$ is at most 
$(q_m/q_k)2\varepsilon_k q_k$, and hence those $i$ have proportion bounded 
by 
$({(q_m/q_k)2\varepsilon_k q_k\over q_m}) = 2\varepsilon_k$ in $I$. It follows 
that the collection of $i\in I$ such that $sh^i(s)\in E^0_k\cup E^1_k\cup E^2_k$ 
is bounded by $3*2\varepsilon_k$.

Numerical requirements \ref{l_n summability} and  \ref{varepsilons} imply that the sum in item \ref{simple revision} is bounded by $1/l_{n-1}+\varepsilon_{n-1}$.\qed

A very similar statement is the following:

\begin{lemma}\label{whole $n$-blocks 1} 
Suppose that $s\in S$ and $s$ has a principal $n$-block. Then $n$ is mature provided that $s\notin \bigcup_{n\le m}E^0_m\cup E^1_m\cup E^2_m$. In particular, if $n$ is mature for $s$ and $s$ is not in a boundary portion of its principal $n-1$-block or in $E^0_{n-1}\cup E^1_{n-1}\cup E^2_{n-1}$, then $n-1$ is mature for $s$.
\end{lemma}

\subsection{The Circle Factor}\label{tcf}
Let $\la k_n,l_n:n\in\nn\ra$ be a circular coefficient sequence and $\la p_n, q_n:n\in\nn\ra$ be the associated sequence defined by formulas \ref{qn} and \ref{pn}. Let $\alpha_n=p_n/q_n$ and $\alpha=\lim \alpha_n$. 

For  a natural number $q\ge 1$, let \hypertarget{IQ}{$\mci_q$} be the partition of $[0,1)$ with atoms $\la \hoo{i}{q}:0\le i<q\ra$, and refer to $\hoo{i}{q}$ as $I^q_i$.\footnote{If $i>q$ then $I^q_i$ refers to $I^q_{i'}$ where $i'<q$ and $i'\equiv i \mod{q}$.} 
Since $p_n$ and $q_n$ are relatively prime, the rotation $\mcr_{\alpha_n}$ enumerates the partition $\mci_{q_n}$ starting with $I^{q_n}_0$. Thus $\mci_{q_n}$ has two natural orderings--the usual geometric ordering and the \emph{dynamical} ordering given by the order that  $\mcr_{\alpha_n}$ enumerates $\mci_{q_n}$. Since $j_i=p^{-1}i$ (mod $q$), $I^q_i$ is the $j_i^{th}$ interval in the dynamical ordering. 

\begin{definition}\label{Dn}
For $x\in [0,1)$ we will write $D_n(x)=j$ if $x$ belongs to the $j^{th}$ interval in the dynamical ordering of $\mci_{q_n}$. Equivalently $D_n(x)=j$ if  $x\in I^{q_n}_{jp_n}$.
 \end{definition}

\bfni{Informal description:} Following \cite{prequel},  for each stage $n$,  we have a periodic approximation $\tau_n$ to $\bk^c$ consisting of towers $\mct$ of height $q_{n}$ whose levels correspond to subintervals of $[0,1)$. This approximation refines the periodic permutation of $\mci_{q_n}$ determined by $\mcr_{\alpha_n}$. If $s$ is mature then $s$ lies is the $r_n^{th}(s)$ level of $\mci_{q_n}$ in the dynamical ordering. Passing from $\tau_n$ to $\tau_{n+1}$ the mature points remain in the same levels of the $n$-towers as they are spread into the 
$n+1$-towers in $\tau_{n+1}$. The towers of $\tau_{n+1}$ can be viewed as cut-and-stack constructions--filling in boundary points between cut $n$-towers.
The fillers are taken from portions of the $n$-towers.

With this view each mature point remains in the same interval of $\mci_{q_n}$ when viewed in $\tau_{n+1}$. Moreover if $s\in J\in \mci_{q_{n+1}}$ and $J\subseteq I\in \mci_{q_n}$, then $\mcr_{\alpha_{n+1}}J\subseteq \mcr_{\alpha_n}I$.

Thus the $n+1$-tower for $\mcr_{\alpha_{n+1}}$ has {multiple} contiguous sequences of levels of length $q_n$ that are sublevels of the $n$-tower and the action of $\mcr_{\alpha_n}$ and $\mcr_{\alpha_{n+1}}$ agree on these levels.

\begin{definition}\label{first appearance of circle factor} Let $\Sigma_0=\{*\}$. We define a  circular construction sequence such that each 
$\mcw^c_n$ has a unique element as follows:
\begin{enumerate}
\item $\mcw^c_0=\{*\}$ and
\item If $\mcw^c_n=\{w_n\}$ then $\mcw^c_{n+1}=\{\mcc(w_n, w_n, \dots w_n)\}$.

\end{enumerate}
Let $\mck$ be the resulting circular system. \end{definition}
It is easy to check that $\mck$ has unique non-atomic measure since the unique $n$-word, $w_n$, occurs exactly $k_n(l_n-1)q_n$ many times in $w_{n+1}$. This measure is ergodic.

\medskip

Let $\bk^c$ be an arbitrary circular system with  coefficients $\la k_n, l_n:n\in\nn \ra$. Then $\bk^c$ has a canonical factor isomorphic to $\mck$. This canonical factor plays a role for circular systems analogous to the role  odometer transformations play for odometer based systems.

To see $\mck$ is a factor of $\bk^c$,   define the following function:

\begin{equation}\label{definition of factor map}
\pi(x)(i) = \left\{ \begin{array}{ll}x(i) &
\mbox{if $x(i)\in \{b, e\}$} \\
* &\mbox{otherwise}
\end{array}\right.
\end{equation}

\bfni{Notation:} Write $w_n^\alpha$ for the unique element of $\mcw^c_n$ in the construction sequence for $\mck$. Then $w_n^\alpha$ lies in the principal $n$-block of the projection to $\mck$ of any $s\in \bk^c$ for which $n$ is mature.

\begin{theorem}\label{rank one description}(\cite{prequel}, Theorem 43.)
Let $\nu$ be the unique non-atomic shift-invariant measure on $\mck$. Then 
\[(\mck, \mcb, \nu, sh)\cong (S^1, \mathcal D, \lambda, \mcr_\alpha)\]
 where $\mcr_\alpha$ is the rotation of the unit circle by $\alpha*2\pi$ radians and $\mathcal{B, D}$ are the $\sigma$-algebras of measurable sets.\end{theorem}

 The isomorphism $\phi_0:\mck\to S^1$ asserted to exist  in Theorem \ref{rank one description} is constructed as a limit of functions $\rho_n$, where $\rho_n$ is defined  by setting
{\begin{equation}\label{def of rhon}
{\rho}_n(s)={i\over q_n}
\end{equation}}
iff 
$I^{q_n}_i$ is the $r_n(s)^{th}$ interval in the dynamical ordering.\footnote{Thus $r_n$ and $\rho_n$ both have the same subset of  $S$ as their domain and contain the same information. They map to different places $r_n:S\to \nn$, whereas $\rho_n:S\to\zoo$ and is the left endpoint of the $r_n^{th}$ interval in the dynamical ordering.} Equivalently, since 
the $r_n^{th}$ interval in the geometric ordering is $I^{q_n}_{p_nr_n(s)}$:
	\begin{equation}
	i\equiv p_nr_n(s)\mod{q_n}
	\end{equation}
The following follows from Proposition 44 in \cite{prequel}.
\begin{prop}\label{Dn and rn} Suppose that $n$ is mature for $s$, then
\begin{equation*}
r_n(s)=D_n(\phi_0(s))
\end{equation*}
\end{prop}
{The proof of Theorem \ref{lele} requires understanding the correspondence between the geometric construction and its symbolic representation. The words in $\mcw_n$ correspond to cut-and-stack constructions, passing from stage $n$ to $n+1$ via the $\mcc$ operator corresponds to basing the cut and stack construction on $\mcr_{\alpha_{n+1}}$ which agrees with the $\mcr_{\alpha_n}$ for most consecutive intervals of length $q_n$. A first step in understanding this correspondence is the next remark and lemma. }

\begin{remark}\label{coherence remark}
It will be helpful to understand $\phi_0^{-1}$ explicitly. To each point $x$ in the range of $\phi_0$, 
$s=\phi_0^{-1}(x)$ belongs to $S$.  By Lemma \ref{specifying elements}, to determine $s$ it suffices to 
 know $\la r_n(s):n\ge N\ra$ for some $N$ as well as the sequence $\la w_n:n\ge N\ra$ of principal subwords of $s$. Since 
 we are working with $\mck$, the only choice for $w_n$ is $w_n^\alpha$. For mature $n$, Proposition \ref{Dn and rn} 
 tells us that $r_n(s)=D_n(x)$.  Thus $s$ is the unique element of $S$ with the property that $\la r_n(s):n\in\nn\ra$ 
 agrees with $\la D_n(x):n\in\nn\ra$ for all large $n$.
\end{remark}

We isolate the following fact for later use:

\begin{lemma}\label{coherence of rn}
Suppose that $\phi_0(s)=x$ and $n<m$ are mature for $s$. Then if $I$ and $J$ are  the $D_n(x)^{th}$ and $D_m(x)^{th}$  intervals in the dynamical orderings of $\mci^{q_n}$ and $\mci^{q_m}$, then $J\subseteq I$.
\end{lemma}
The natural way of representing the complex unit circle as an abelian group is multiplicatively: the rotation by $2\pi \alpha$ radians is multiplication by $e^{2\pi i\alpha}$.  It is often  convenient to 
identify the unit circle with $[0,1)$. In doing so, multiplication by $e^{2\pi i\alpha}$ corresponds to ``mod one" addition and the complex conjugate $\bar{z}$ corresponds to $-z$.

The following result is standard:
\begin{prop}\label{nonsense again} Let $\alpha\in [0,1)$ be irrational.
Suppose that $T:S^1\to S^1$ is an invertible measure preserving 
transformation that commutes with $\mcr_\alpha$. Then for some $\beta$, 
$T=\mcr_\beta$ almost everywhere. Identifying $S^1$ with $[0,1)$  there is a 
$\beta$ such that for almost all $x\in S^1$
\begin{equation}
\label{commutator of rotation}
T(x)=x+\beta\mod 1.
\end{equation}
 It follows that if $T$ is an isomorphism between $\mcr_\alpha$ and $\mcr_{\alpha}^{-1}$, then $T(x)=-{x}+\beta\mod 1$.
\end{prop}

{
\begin{definition} \label{identity crisis} Using the identification of $S^1$ with $[0,1)$ we view $\phi_0:\mck\to [0,1)$. Given a rotation $\mcr_\beta$, we get a map $\mcs_\beta:\mck\to \mck$ such that 
\[\mcs_\beta(s)=\phi_0^{-1}\mcr_\beta\phi_0(s).\]
\end{definition}
We will occasionally abuse notation and write $s+\beta$ for $\mcs_\beta(s)$.}

\subsection{Points of view}\label{roshoman}
Circular systems can be viewed from multiple perspectives: geometrically, as limits of periodic processes\footnote{See section 5 of \cite{prequel} for the formal definition.} and as symbolic shifts.

The $n^{th}$ periodic process consists of a collection of $s_n$ periodic towers with each tower having one level designated as a base.  To pass from $\tau_n$ to $\tau_{n+1}$ the bulk of the $\tau_n$-towers are repeated $q_n(k_{n})(l_n-1)$ many times in blocks of length $l_n-1$ in each $\tau_{n+1}$-tower. In between these blocks there are filler levels.

The words $w\in \mcw_n^c$ are in one-to-one correspondence with the towers in $\tau_n$.
The ``$\mcc$" operation encodes the transition from $\tau_n$ to $\tau_{n+1}$.  The towers in $\tau_{n+1}$ correspond to words $\mcc(w_0,\dots w_{k_n-1})$. {Each $\tau_n$-tower $T_j$ has a corresponding word $w_j\in \mcw_n$.  Repeating stacking of $T_j$ corresponds to the powers of $w_j$ in 
$\mcc(w_0,\dots w_{k_n-1})$.}
The levels of a tower in $\tau_{n+1}$ are either contained in levels of $\tau_n$-tower or are filler blocks labelled ``$b$" or ``$e$."
The repetitions of each $w_i$ in $0$-subsections correspond to stacking parts of the levels of the corresponding tower in $\tau_n$ periodically $l_n-1$ times.

The circle factor $\mck_\alpha$ captures exactly the structure of the \emph{levels} of the towers and how they interact as one moves from $\tau_n$ to $\tau_{n+1}$. This is the idea behind for the  construction of the isomorphism between $(\mck_\alpha, \nu,sh)$ and $(S^1,\lambda, \mcr_\alpha)$ and made explicit in Proposition \ref{Dn and rn}.
  
Given an $s\in \bk^c$ that is mature for $n\le m$ we can view its restriction to its principal $m$-subword as a particular tower in $\tau_m$. Since $s$ is mature for $m$, the principal subword is repeated many times on either side of $s(0)$. In particular we see:
    \begin{remark} \label{life is easy in the interior}
    Suppose that $n$ is mature for  $s\in S\subseteq \bk^c$, $n\le m$  and $0\le d<q_m$.  Then 
    \begin{eqnarray}\label{what goes around}
    r_n(sh^d(s))\equiv_{q_n}d+r_n(s)
    \end{eqnarray}
    \end{remark} 
The circle factor $\mck_\alpha$ of $\bk^c$ punctuates the elements of $S\subseteq \bk^c$. Since there is only one word in each element of the construction sequence for $\mck_\alpha$, we can view the levels of its tower as being of the form $[i/q_n, (i+1)/q_n)$ in the dynamical ordering. Then the cyclic permutation of these levels given by $\mcr_{p_n/q_n}$.  This permutation preserves the dynamical ordering and, for $s$ that are mature at stage $n$, reflect the behavior of $r_n(s)$.

\subsection{The Natural Map}\label{natural definition}
A specific isomorphism $\natural:(\mck, sh)\to (\rev{\mck}, sh)$ will serve as a benchmark for 
understanding of potential  maps $\phi:\bk^c\to\rev{\bk^c}$. Viewing $\mcr_\alpha$ as a rotation of the unit circle by $\alpha*2\pi$ radians one can view the transformation $\natural$  as a symbolic analogue of complex 
conjugation $z\mapsto \bar{z}$ on the unit circle, which is an isomorphism between $\mcr_\alpha$ and $\mcr_{-\alpha}$.  Indeed, by 
Theorem \ref{rank one description}, $\mck\cong\mcr_\alpha$ and so $\rev{\mck}\cong \mcr_{-{\alpha}}$. Copying 
$\natural$ over to a map on the unit circle will give an isomorphism $\phi$ between $\mcr_\alpha$ and $\mcr_{{-\alpha}}$.  If we view $z$ and $\alpha$ as elements of the unit interval and the rotation as addition modulo 1, 
Proposition \ref{nonsense again} says that
such an isomorphism must be of the form 
\[\phi(z)=-z+\beta\]
for some $\beta$. It follows immediately from this characterization that $\natural$ is an involution.\footnote{The particular $\beta$  given by $\natural$ is determined by the specific variation of the definition one uses--indeed any \emph{central} value can occur as a $\beta$. (See section \ref{centralizer and central} for the definition and use of central values.)} 

The map $\natural$ is defined as the limit of a sequence of codes 
$\la\Lambda_n:n\in\nn\ra$  that converge  to an isomorphism from 
$\mck$ to $\rev{\mck}$ (see \cite{global_structure} for more details). The $\Lambda_n$ will be shifting and reversing words. The amount of shift is determined by the Anosov-Katok coefficients $p_n, q_n$ defined in equations  \ref{pn} and \ref{qn}.

Let $A_0=0$ and inductively
\begin{equation}
A_{n+1}=A_n-(p_n)^{-1}. \label{code coefficients}
\end{equation}
It is easy to check that 
\begin{equation}|A_{n+1}|<2q_{n}\label{An is small}
\end{equation}

Define a stationary code $\overline{\Lambda}_n$ with domain $S$ that approximates elements of $\rev{\mck}$
 by  defining
\begin{equation}\label{definition of Lambdan}
\Lambda_n(s)=\left\{\begin{array}{ll} sh^{A_n+2r_n(s)-(q_n-1)}(\rev{s})(0) & \mbox{if $r_n(s)$ is defined}\\
				b&\mbox{otherwise}
				\end{array}\right.
\end{equation}
The following result appears in \cite{global_structure}:
\begin{theorem}\label{mr natural}
 The sequence of stationary codes $\la \overline\Lambda_n:n\in\nn\ra$  converges to a 
 {shift invariant function} $\overline{\natural}:\mck\to (\{*\}\cup \{b, e\})^{{\poZ}}$  that induces an 
 isomorphism  $\natural$ from $\mck$ to $\rev{\mck}$.
\end{theorem}
Remark 78 of \cite{global_structure} implies that the convergence is prompt: for a typical $s$ and all large enough $n$, $\natural(s)$ agrees with $\bar{\Lambda}_n(s)$ on the principal $n$-block of $s$.

\paragraph{Caveat} Since $(\bk^c)^{-1}=(\bk^c, sh^{-1})$ is trivially isomorphic to $(\rev{bk^c}, sh)$ we often don't distinguish them.  However, as in  Definition \ref{synch and anti-synch} of the \emph{synchronous} and \emph{anti-synchronous} joinings, the notational distinction becomes important. 

When viewing $(\bk^c)^{-1}$ and $\bk^c$ with the backwards shift and considering the action on the circle factor instead of using $\natural$,  one must use
\begin{equation}\label{getting things backwards}
\rev{}\circ\natural
\end{equation}
instead of simply $\natural$.

\subsection{Categories and the Functor $\mcf$.}
Fix a \hyperlink{circ coef}{circular coefficient sequence} $\la k_n, l_n:n\in\nn\ra$.
Let $\Sigma$ be a language and $\la \mcw_n:n\in\nn\ra$ be a   construction sequence  for an 
odometer based system with coefficients $\la k_n:n\in\nn\ra$.    Then for each $n$ the operation $\mcc_n$ is well-defined. Define a  construction sequence 
$\la \mcw_n^c:n\in\nn\ra$ and bijections 
$c_n:\mcw_n\to \mcw_n^c$ by induction as follows:

\begin{enumerate}
\item Let $\mcw^c_0=\Sigma$ and $c_0$ be the identity map.
\item Suppose that  $\mcw_n, \mcw_n^c$ and $c_n$ have already been defined. 
\[\mcw_{n+1}^c=\{\mcc_n(c_n(w_0),c_n(w_1), \dots c_n(w_{k_n-1})):w_i\in \mcw_n \mbox{ and }w_0w_1\dots w_{k_n-1}\in \mcw_{n+1}\}.\]
{(Words in $\mcw_{n+1}$ are concatenations of $k_n$ words in $\mcw_n$ and so can be written in the required form: as $w_0w_1\dots w_{k_n-1}$ with $w_j\in \mcw_n$.)}

Define the map $c_{n+1}$ by setting
 \[c_{n+1}(w_0w_1\dots w_{k_n-1})=\mcc_n(c_n(w_0),c_n(w_1), \dots c_n(w_{k_n-1})).\] 
 
 \end{enumerate}
Note in case 2 the \hyperlink{pwords}{{prewords}} are:
\[P_{n+1}=\{(c_n(w_0),c_n(w_1),\dots c_n(w_{k_n-1})): w_0w_1\dots w_{k_n-1}\in \mcw_{n+1}\}.\] 

\begin{remark}\label{stick on} Some useful facts are:
\begin{itemize}
\item It follows from Lemma \ref{stabilization of names 1} and Numerical Requirement \ref{l_n summability}  that if $\la \mcw_n:n\in\nn\ra$ is an odometer based construction sequence, then $\la \mcw_n^c:n\in\nn\ra$ is a construction sequence; i.e. the spacer proportions are summable. 
\item If each $w\in \mcw_n$ occurs exactly the same number of times in every element of $\mcw_{n+1}$, then $\la \mcw_n^c:n\in\nn\ra$ is strongly uniform. 
\item Odometer words in $\mcw_n$ have length $K_n$.  The length of the circular words in $\mcw_n^c$ is $q_n$.
\end{itemize}
\end{remark}

\begin{definition}\label{def of functor}
Define a map $\mcf$ from the set  of odometer based subshifts to  circular subshifts  as follows. Suppose that 
$\bk$ is an odometer based shift built from a construction sequence $\la \mcw_n:n\in\nn\ra$. Define 
\[\mcf(\bk)=\bk^c\]
where $\bk^c$ has construction sequence $\la \mcw_n^c:n\in\nn\ra$.
\end{definition}

The map $\mcf$ is one to one by the unique readability of words in $\mcw$. Suppose that $\bk^c$ is a circular system with coefficients 
$\la k_n, l_n:n\in\nn\ra$.  We can  recursively build 
functions $c_n\inv$ from words in $\Sigma\cup \{b,e\}$ to words in $\Sigma$. The 
result is a  
odometer based system $\la \mcw_n:n\in\nn\ra$ with coefficients $\la k_n:n\in\nn\ra$.
If $\bk$ is the resulting odometer based system then 
$\mcf(\bk)=\bk^c$. Thus $\mcf$ is a bijection. 
\bigskip

If $\bk$ is an odometer based system, denote the odometer base by $\bk^\pi$ and let $\pi:\bk\to \bk^\pi$ be the canonical factor map. If $\bk^c$ is a circular system, let $(\bk^c)^\pi$ be the rotation factor  $\mck$ and $\pi:\bk^c\to \mck$ be the canonical factor map. For both odometer based and circular systems the underlying canonical factors serve as  timing mechanisms. This motives the following.

\begin{definition}\label{synch and anti-synch}  \emph{Synchronous} and \emph{anti-synchronous} joinings are defined as follows:\footnote{{We use $\mcl$ for the notation for the rotation factor of a circular system $\bl^c$. In this context, when taking inverses of symbolic systems we keep the same orientation for the symbolic system and use $sh^{-1}$}}

\begin{enumerate}
\item Let $\bk$ and $\bl$ be odometer based systems with the same coefficient sequence, and $\rho$ a joining between $\bk$ and $\bl^{\pm1}$. Then $\rho$ is 
\emph{synchronous} if $\rho$ joins $\bk$ and $\bl$ and the projection of $\rho$ to a joining on $\bk^\pi\times \bl^\pi$ is the graph joining determined by the identity map (the diagonal joining of the odometer factors); $\rho$ is \emph{anti-synchronous} if $\rho$ is a joining of $\bk$ with $\bl^{-1}$ and its projection to 
$\bk^\pi\times (\bl^{-1})^\pi$ is the graph joining determined by the  map $x\mapsto -x$.
\item Let $\bk^c$ and $\bl^c$ be circular systems with the same coefficient sequence and $\rho$ a joining between $\bk^c$ and $(\bl^c)^{\pm 1}$. Then $\rho$ is 
\emph{synchronous} if $\rho$ joins $\bk^c$ and $\bl^c$ and the projection to a joining of $(\bk^c)^\pi$ with $(\bl^c)^\pi$ is the graph joining determined by
the identity map of $\mathcal K$ with $\mathcal L$, the underlying rotations; $\rho$ is \emph{anti-synchronous} if it is a 
joining of $\bk^c$ with $(\bl^c)^{-1}$ and projects to the graph joining determined by $\rev{}\circ\natural$ on $\mck\times \mathcal L^{-1}$.
\end{enumerate}
\end{definition}

\bfni{The Categories} Let ${\mathcal OB}$ be the category whose 
objects are ergodic 
odometer based systems with coefficients $\la k_n:n\in \nn\ra$. The morphisms between objects $\bk$ and $\bl$ will be
 synchronous graph joinings of $\bk$ and $\bl$ or  anti-synchronous graph joinings of $\bk$ and $\bl^{-1}$. We call this the 
\emph{category of odometer based systems.}
\smallskip

Let $\mcc B$ be the category whose objects consists of all ergodic circular systems with coefficients $\la k_n,l_n:n\in\nn\ra$.  The morphisms between objects $\bk^c$ and $\bl^c$ will be synchronous graph joinings of $\bk^c$ and $\bl^c$ or anti-synchronous graph joinings of $\bk^c$ and $(\bl^c)^{-1}$.
We call this the \emph{category of circular systems.}
\medskip

The main theorem of \cite{global_structure} is the following:
\begin{theorem}\label{grand finale}
For a fixed circular coefficient sequence $\la k_n, l_n: n\in\nn\ra$ the categories $\mco B$ and $\mcc B$ are isomorphic by a function $\mcf$ that takes synchronous joinings to synchronous joinings, anti-synchronous joinings to anti-synchronous joinings,  
isomorphisms to isomorphisms and {weakly mixing extensions to weakly mixing extensions.}\footnote{Glasner showed that it takes compact extensions to compact extensions.}
\end{theorem}

It is also easy to verify that the map $\la \mcw_n:n\in\nn\ra\mapsto \la \mcw_n^c:n\in\nn\ra$ takes uniform construction sequences to uniform construction sequences and strongly uniform construction sequences to strongly uniform construction sequences.

\begin{remark}
Were we to be completely precise we would take  objects in $\mco B$ to be \emph{presentations} of odometer based systems by construction
  sequences $\la \mcw_n:n\in\nn\ra$ without spacers and the objects in $\mcc B$ to be \emph{presentations} by circular construction sequences. 
  This subtlety does not cause problems in the sequel so we ignore it.
\end{remark}

\subsection{Propagating Equivalence Relations and Actions}\label{propagano}

{In \cite{FRW}, the number $M(s)$ is the first stage in the tree for which $\sigma_m$ has length $s$.  It is the first stage that the equivalence relation $\mcq^m_s$ is defined.}

 The main result of \cite{FRW} is the existence of a  continuous function from the space of trees to odometer based transformations that reduces ill-founded trees to ergodic transformations isomorphic to their inverses. Components of the construction include equivalence relations $\la \mcq^n_s:M(s)\le n, s\in \nn\ra$ and groups $\la G^n_s:M(s)\le n, s\in \nn\ra$. Some of their properties are:
\begin{enumerate}
\item $M$ is a monotone, strictly increasing function from $\nn$ to $\nn$,
\item $\mcq^0_0$ is the trivial equivalence relation with one equivalence class on $\mcw_0=\Sigma$.
\item $\mcq^n_s$ is an equivalence relation on $\mcw_n$

\item For $n\ge M(s)+1$, viewing elements of $\mcw_n$ as concatenations of words in $\mcw_{M(s)}$,   $\mcq^{n}_s$ is the product equivalence
relation of
$\mcq^{M(s)}_s$. Hence we can view
$\mcw_n/\mcq^n_s$ as sequences of elements of $\mcw_{M(s)}/\mcq^{M(s)}_s$ and similarly for
$\rev{\mcw_n/\mcq^n_s}$. {These sequences have length $K_n$ and are made of $K_n/K_{M(s)}$ many constant blocks of length $K_{M(s)}$.}
\item The groups $\la G^n_s:M(s)\le n, s\in\nn\ra$ are direct sums of copies of 
$\poZ_2$ that have a designated canonical collection of free generators.\footnote{These groups are described in detail in Section \ref{old specs}.}
Each $G^{n+1}_s=G^n_s\oplus H$, where $H$ is either $\poZ/2\poZ$ or $H$ is trivial.
\item Each group  $G^n_s$ acts freely on $\mcw_n/\mcq^n_s\cup \rev{\mcw_n/\mcq^n_s}$ in a manner that even parity group elements preserve the sets $\mcw_n/\mcq^n_s$ and  $\rev{\mcw_n/\mcq^n_s}$ and the odd parity group elements send elements of $\mcw_n/\mcq^n_s$ to  $\rev{\mcw_n/\mcq^n_s}$.

\item The action of $G^{n}_s\subseteq G^{n+1}_s$ on $\mcw_{n+1}\cup \rev{\mcw_{n+1}}$ is propagated from 
$\mcw_n\cup\rev{\mcw_n}$ by the \emph{skew-diagonal} action:
if $g\in G^n_s$ is a canonical generator and 
$w\in \mcw_{n+1}\cup\rev{\mcw_{n+1}}$ is of the form $w_0w_1\dots w_{k_n-1}$ then 
{\[gw=gw_{k_{n-1}}\dots gw_1gw_0.\]}
\end{enumerate}

\medskip

We now define corresponding equivalence relations and group actions on 
 $\la \mcw_n^c:n\in\nn\ra$. They will be used in section \ref{intro timing} to state the timing assumptions and in section \ref{old specs} which gives the construction specifications from \cite{FRW}.\footnote{If $\mcq$ is an equivalence relation on 
 $\mcw^c$ 
define $\rev{\mcq}$ by $(\rev{w_0},\rev{w_1})\in\rev{\mcq}$ if and only if $(w_0,w_1)\in \mcq$. 
In abuse of notation we will not distinguish between 
$(\mcq^n_s)^c$ as a relation on $\mcw^c_n$,  $(\mcq^n_s)^c\cup \rev{(\mcq^n_s)^c}$ as a relation on $\mcw^c_n\cup \rev{\mcw^c_n}$ or $\mcw_n^c/(\mcq^n_s)^c \cup\rev{\mcw_n^c/(\mcq^n_s)^c}$. }

An inductive understanding of $(\mcq^n_s)^s$ and the $G^n_s$-actions is quite useful.

\medskip

\bfni{Inductive definition of $(\mcq^n_s)^c$:}
Define 
\begin{itemize}
\item 
{$(\mcq^n_0)^c$} to  have exactly one class in each $\mcw^c_n$,
\item For $w_0,w_1\in \mcw_{M(s)}$ put $(c_{M(s)}(w_0), c_{M(s)}(w_1)\in (\mcq^{M(s)}_s)^c$ if and only if $(w_0,w_1)\in \mcq^{M(s)}_s$.

\end{itemize}
Suppose we are given $(\mcq^n_s)^c$ on 
$\mcw^c_n$. Define an equivalence 
relation 
$\mcq$ on $\mcw^c_{n+1}$ by setting 
 $\mcc(w_0, \dots w_{k_n-1})$ equivalent to 
$\mcc(w_0', \dots w_{k_n-1}')$ if and only if for all $i, w_i$ is $(\mcq^n_s)^c$-equivalent to $w_i'$. 
\medskip

Rather than a full definition of the action of $G^{n+1}_s$ on $\mcw^c_{n+1}/(\mcq^{n+1}_s)^c\cup \rev{\mcw^c_{n+1}/(\mcq^{n+1}_s})^c$, we describe the how the action of $G^n_s$ propagates: via the \emph{circular skew diagonal action}:
\smallskip

Identify $\rev{\mcw^c_{n+1}/(\mcq^{n+1}_s)^c}$ with the collection of sequences of the form
\[\mcc^{r}(\rev{[w_0]_{(\mcq^n_s)^c}}, \rev{[w_1]_{(\mcq^n_s)^c}}, \dots, \rev{[w_{k_n-1}]_{(\mcq^n_s)^c})}\] 
as $w_0w_1 \dots w_{k_n-1}$ ranges over the elements of  $P_n$.

To define the skew-diagonal action of $G^n_s$ on classes of circular words it suffices to specify it on the canonical 
generators, This is done by setting\footnote{We use $[w_i]$ to denote 
 $[w_i]/(\mcq^n_s)^c$.} 
 \begin{equation*}
g\mcc([w_0], [w_1] \dots [w_{k-1}])=_{def}{\mcc^r([gw_0], [gw_1], \dots [gw_{k-1}])}
\end{equation*}
whenever $g$ is a canonical generator of $G^n_s$.  Note that  the skew-diagonal action has the property that the canonical generators take elements of $\mcw^c_{n+1}/(\mcq^{n+1}_s)^c$ to elements of $\rev{\mcw^c_{n+1}/(\mcq^{n+1}_s)^c}$. It follows that the even parity elements of $G$  leave the sets $\mcw_{n+1}^c/(\mcq^{n+1}_s)^c$  and  $\rev{\mcw_{n+1}^c/(\mcq^{n+1}_s)^c}$ invariant  and odd parity elements of $G$ take $\mcw_{n+1}^c/(\mcq^{n+1}_s)^c$ to elements of $\rev{\mcw_{n+1}^c/(\mcq^{n+1}_s)^c}$ and vice versa.
\medskip

As in \cite{FRW} the equivalence relations $\la \mcq^n_s:n\in\nn\ra$ define factors $\bk_s$ of $\bk$ and similarly $\la (\mcq^n_s)^c:n\in\nn\ra$ define factors $\bk_s$ of $\bk^c$
  The equivariant definitions given here imply that $\mcf$ takes each $\bk_s$ to $\bk^c_s$ and respects the actions of the $G^n_s$.

\section{Understanding Rotations}\label{rotation section}

Let $\mck$ be a rotation factor of a circular system with coefficient sequence 
$\la k_n, l_n:n\in\nn\ra$. This section analyzes how automorphisms of $\mck$ affect the parsing 
of elements of $\mck$. \medskip

Let $(\bk^c,\mu^c)$ and $(\bl^c,\nu^c)$ be two circular systems with that share a given circular coefficient sequence and let $\alpha=\lim \alpha_n$.  Any  isomorphism between ${\bk^c}$ and $(\bl^c)^{\pm 1}$
 induces a unitary isomorphism $U_\phi$ from $L^2((\bl^c)^{\pm 1})$ to  $L^2({\bk^c})$, and this isomorphism sends eigenfunctions for $n \alpha$ to eigenfunctions for 
 $n\alpha$. Thus 
  every  isomorphism  has to send the canonical factor $\mck_\alpha$  of ${\bk^c}$ to the canonical factor
$\mck_\alpha^{\pm 1}$ of $(\bl^c)^{\pm1}$. 
Explicitly:
suppose  that $\phi:{\bk^c}\to (\bl^c)^{\pm 1}$ is an isomorphism. Then $U_\phi:  L^2((\bl^c)^{\pm 1})\to L^2({\bk^c})$, 
and $U_\phi$ takes the space generated by eigenfunctions of $U_{sh}$ in $L^2((\bl^c)^{\pm 1})$ with eigenvalues 
$\{\alpha^n:n\in \poZ\}$ to the  space generated by corresponding eigenfunctions in $L^2({\bk^c})$.  
Consequently there is a measure preserving transformation $\phi^\pi$ making the following diagram 
commute:
\begin{equation}
\label{phi pi}
\begin{diagram}
\node{{\bk^c}}\arrow[1]{e,r}{\phi}\arrow[1]{s,r}{\pi}\node{(\bl^c)^{\pm 1}}\arrow{s,r}\pi\\
\node{\mck_\alpha}\arrow{e,t}{\phi^\pi}\node{\mck_{\alpha}^{\pm 1}}
\end{diagram}
\end{equation}

 { By Theorem \ref{rank one description},  $\mck_\alpha$ is conjugate to the rotation $\mcr_\alpha$ 
 of the unit 
 circle by a map $\phi_0$. Hence (using additive notation) $\phi^\pi$ must be conjugate to a 
 transformation defined 
 on the unit interval of the form $x\mapsto z+\beta$ for some $\beta\in [0,1)$, where $z$ is either $x$ or 
 $-{x}$, 
 depending on whether $\phi^\pi$ maps to $\mck_\alpha$ or 
$\mck_\alpha^{-1}$.  Since $\rev{}\circ\natural:\mck_\alpha\to \mck_\alpha^{-1}$ is an isomorphism, if 
$\phi$ maps to 
$({\bl^c})^{-1}$, $\rev{}\circ\natural(x)$ can serve as an alternative to the  benchmark to the map $x\mapsto -{x}$. Explicitly:  
the $\beta $ associated to $\phi$ is the number making $\phi^\pi(s)=\rev{}\circ\natural(\mathcal S_\beta(s))$; {equivalently, 
$\rev{}\circ\natural^{-1}\circ\phi^\pi(s)=\mathcal S_\beta(s)$.}\footnote{{The reader is referred to the Caveat at the end of section \ref{natural definition}, for the reason $\rev{}\circ \natural$ is used.}}}

\smallskip

\noindent {Summarizing,
	\begin{itemize}
	\item[A.)] If $\phi:\bk^c\to \bl^c$ is an isomorphism, then viewed as a map from $[0,1)$ to $[0,1)$, there is a 
	unique $\beta\in [0,1)$ for almost every $x$,  $\phi^\pi(s)=\mcs_\beta(s)$ .
	\item[B.)] If $\phi:\bk^c\to (\bl^c)^{-1}$ then there is a unique $\beta$ for almost every $s$, 
	$\phi^\pi(x)=\rev{}\circ\natural(\mathcal S_\beta(s))$.
	\end{itemize}
\begin{definition}
In  cases A.) and B.), we call the map $\mcs_\beta$  the \emph{rotation associated with $\phi$}.
\end{definition}
}

We record the following facts:
\begin{lemma}\label{homos}Let ${\bk^c}$ be a  circular system. Then
\begin{enumerate} 
\item The set of $\beta$ associated with automorphisms of $\bk^c$ form a group.

\item If $\phi:{\bk^c}\to {(\bk^c)}\inv$ and $\psi:{\bk^c}\to {\bk^c}$ are isomorphisms where $\phi^\pi=\rev{}\circ\natural\circ\mcs_\beta$ and $\psi^\pi=\mcs_\gamma$, then $(\phi\circ\psi)^\pi=\rev{}\circ\natural \circ\mcs_\delta$ where 
 $\delta=\beta+\gamma$.
\end{enumerate}
\end{lemma}
\pf It is easy to check that 
\begin{itemize}
\item If $\phi, \psi$ are isomorphisms from ${\bk^c} $ to ${\bk^c}$ with $\phi^\pi=\mcs_\beta$ and 
$\psi^\pi=\mcs_\gamma$, then  $(\phi\circ\psi)$ is also an isomorphism from $\bk^c$ to $\bk^c$ and $(\phi\circ\psi)^\pi=\mcs_\delta$, where $\delta=\beta+\gamma$.

\item If $\phi$ is an isomorphism from $\bk^c$ to $\bk^c$, and $\phi^\pi=\mcs_\beta$, then 
$(\phi^{-1})^\pi=\mcs_{-{\beta}}$.

\end{itemize}
The second assertion is similar.
\qed

Given a rotation $\mcr_\beta$,  set
 \begin{equation*}S(\beta)=\bigcap_{n\in\poZ}\mcs_\beta^n(S)
 \end{equation*}
 This can be described independently of $\mcs_\beta$ as:
\begin{equation*} 
\{s\in S:\mbox{for all }n\in \mathbb Z, \phi_0(s)\in (\phi_0[S]+n\beta)\}. 
 \end{equation*}
 It is clear that  $\nu(S(\beta))=1$.

  Define a sequence of functions $\la d^n:n\in\nn\ra$. Each 
  \[d^n:S(\beta)\to \{0, 1, 2, \dots q_n-1\}.\]
  For  $s\in S(\beta)$ and $t=\mcs_\beta(s)$ we have $t\in S(\beta)$ and $\phi_0(t)=\mcr_\beta\phi_0(s)$.  All large enough $n$ are mature for 
  $t$, and 
 $t$ is determined by a tail segment of $\la r_n(t):n\in\nn\ra$. 
 
 \begin{definition}
 If $n$ is mature for both $s$ and  $t=\mcs_\beta(s)$, let 
 \begin{equation}\label{def of dn}
 d^n(s)\equiv_{q_n}r_n(t)-r_n(s),
 \end{equation}
 and $d^n(s)=0$ otherwise. {(We could have made a more general definition $d^n(s,t)$ for arbitrary $t$ and take
$t=\mcs_\beta(s)$ when we want to use $d^n(s)$.)}
 \end{definition}

 Explicitly:     from the definition of $r_n$, 
  $\phi_0(s)+\beta$ belongs to the $(r_n(s)+d^n(s))^{th}$ interval in the 
  dynamical ordering of $\mci_{q_n}$.\footnote{More accurately: if $j<q_n$ and $j\equiv_{q_n}r_n(s)+d^n(s)$, 
  then $\phi_0(s)+\beta$ belongs to the $j^{th}$ interval in the dynamical ordering of $\mci_{q_n}$. {Recall  the relationship between symbolic shifts and the towers  of intervals in the dynamical ordering given in  
Section \ref{roshoman}.}}

Fix an $n$ and suppose that $\beta$ is not a multiple of $1/q_n$. Then the interval $[\beta, \beta+1/q_n)$ 
intersects two geometrically consecutive intervals of the form $[i/q_n, (i+1)/q_n)$.  

\begin{lemma}\label{two values} Suppose that $n$ is mature for $s$ and $\mathcal S_\beta(s)$. Then $d^n(s)$ belongs to 
$ \{D_n(\beta),D_n(\beta+1/q_n)\}$. Thus there are only two possible values for $d^n(s)$ and these values differ by $j_1$. 
\end{lemma}

\begin{figure}[!h]
\centering
\includegraphics[height=.40\textheight]{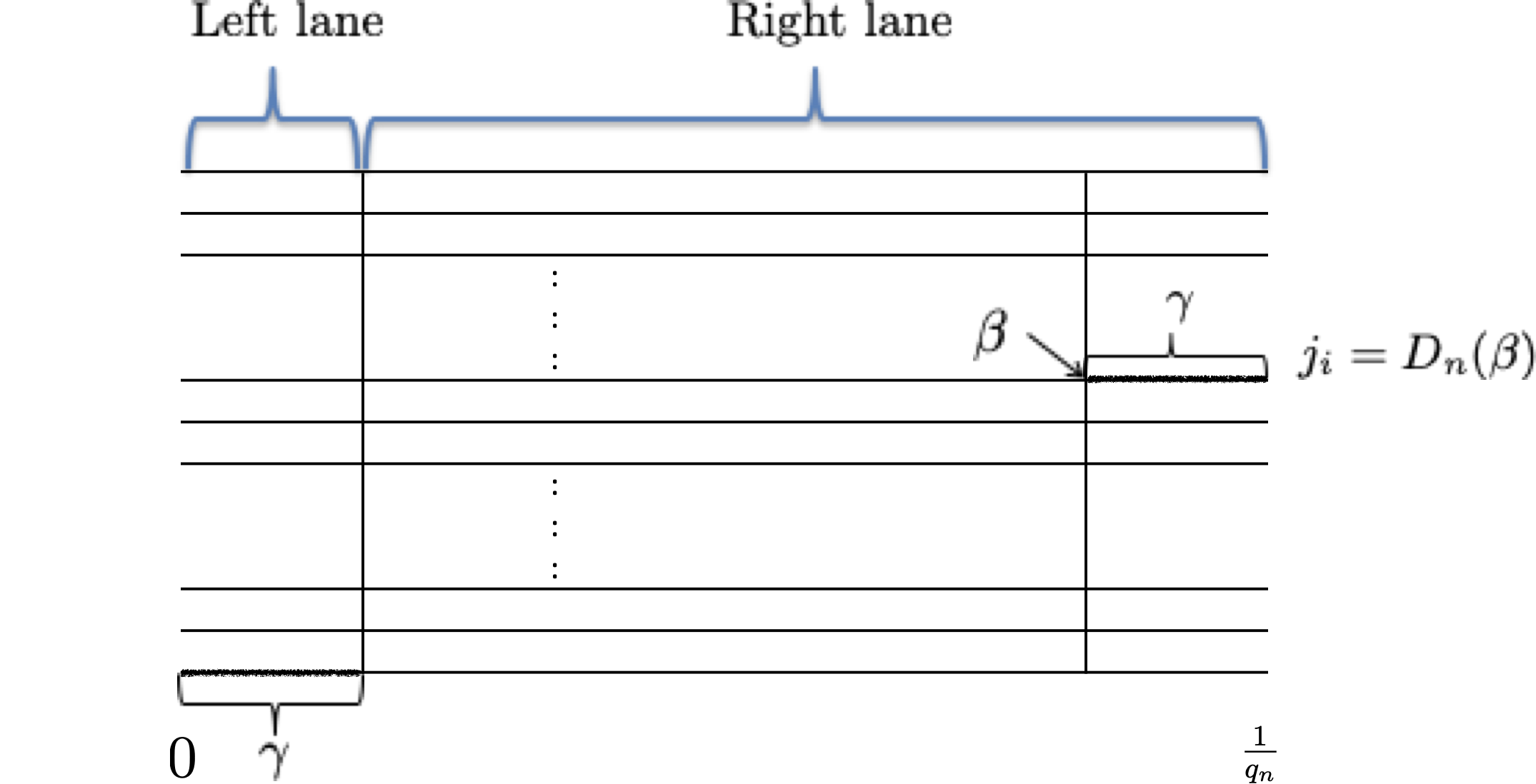}
\caption{Left lane and Right lane of the $q_n$-tower}
\label{fig left right}
\end{figure}

\pf Suppose that $\beta\in [i/q_n, (i+1)/q_n)$ and $\gamma=(i+1)/q_n-\beta$. Then $D_n(\beta)=j_i$. We 
claim that, relative to those $s$ for which $n$ is mature for both $s$ and $\mcs_\beta(s)$, $d^n$ is constant on 
{$\phi_0^{-1}(\bigcup_{j<q_n}[j/q_n, (j+1)/q_n-\gamma))$} and  on 
$\phi_0^{-1}(\bigcup_{j<q_n}[(j+1)/q_n-\gamma, (j+1)/q_n))$, where it takes values $D_n(\beta)$ and 
$D_n(\beta+{1\over q_n})$ respectively (see figure \ref{fig left right}).

We show that $d^n$ is constant on the first set. Suppose that $n$ is mature for $s, \mcs_\beta(s)$ and  
$\phi_0(s)=x$ belongs to the interval $[0, \gamma)$. Then $x+\beta\in \hoo{i}{q_n}$. Hence
 $r_n(\mcs_\beta(s))=j_i=D_n(\beta)$. Since $r_n(s)=0$ we know that $d^n(s)=j_i$. Now suppose that
{$s^*\in \phi_0^{-1}(\bigcup_{j<q_n}[j/q_n, (j+1)/q_n-\gamma))$} and $n$ is mature for {$s^*$} and $\mcs_\beta({s^*})$. 
Let $k=r_n(s^*)$. Then 
  $\phi_0(t)=x+kp_n/q_n$  for some $x\in [0,\gamma)$. So $\phi_0(s^*)+\beta\in[(i+1+kp_n)/q)-\gamma, (i+1+kp_n)/q)$. Hence 
	  \begin{eqnarray*}r_n(\mcs_\beta(s^*))&=&(p_n)^{-1}(i+kp_n)\\
	  &=&j_i+k.
	\end{eqnarray*}
Thus
	\begin{eqnarray*}d^n(s^*)&=&r_n(\mcs_\beta(s^*))-r_n(s^*)\\
					&=&j_i+k-k\\
					&=&j_i.	
	\end{eqnarray*}
If $s^*\in \phi_0^{-1}(\bigcup_{j<q_n}[(j+1)/q_n-\gamma, (j+1)/q_n))$ the proof is parallel.

Finally $\beta$ and $\beta+{1\over q_n}$ fall into consecutive intervals of $\mci^{q_n}$ in the geometric ordering, and hence $D_n(\beta+{1\over q_n})=D_n(\beta)+j_1$.
\qed

Define $d^n_L$ and $d^n_R$ by setting $d^n_L=D_n(\beta)$ and $d^n_R= D_n(\beta+{1\over q_n})$. Let
\begin{equation*}
L_n=\{s:s \mbox{ is mature at stage $n$ and } r_n(s)+d^n_L\equiv_{q_n}r_n(\mcs_\beta(s))\}
\end{equation*}
and 
\begin{equation*}
R_n=\{s:s \mbox{ is mature at stage $n$ and }r_n(s)+d^n_R\equiv_{q_n}r_n(\mcs_\beta(s))\}
\end{equation*}
We refer to $L_n$ and $R_n$ as the \emph{left lane} and \emph{right lane} respectively.

\

\noindent{\bf Notation}: Let $\beta^L_n, \beta^R_n$, be the measures of the left and right lanes at stage $n$.

\begin{lemma}\label{characterization of displacement}
Consider $(\mck, \nu,sh)$ and let $\iota_n$ be the measure of the collection of $s$ that are  not mature at stage $n$. Then:
    \begin{enumerate}
    \item $\ceil{\qn\beta}-q_n\beta\ge\beta_n^L\ge \ceil{\qn\beta}-q_n\beta-\iota_n$, 
    \item $\qn\beta-\floor{\qn\beta} \ge \beta_n^R\ge \qn\beta-\floor{\qn\beta}-\iota_n$
    \item $\beta_n^L+\beta_n^R+\iota_n=1$
    \end{enumerate}
In particular $\sum \beta_n^L<\infty$ if and only if $\sum(\ceil{\qn\beta}-\beta)<\infty$ and $\sum \beta_n^R<\infty$ if and only if $\sum(\qn\beta-\floor{\qn\beta})<\infty$.
\end{lemma}
\pf Let $M_n$ be the collection of $S$ that are mature at stage $n$. 
In the proof of Lemma \ref{two values}, we showed that  $L_n$ is 
{$\phi_0^{-1}(\bigcup_{j<q_n}[j/q_n, (j+1)/q_n-\gamma))\cap M_n$} and $R_n$ is 
{$\phi_0^{-1}(\bigcup_{j<q_n}[(j+1)/q_n-\gamma, (j+1)/q_n))\cap M_n$}, 
where $\gamma=(i+1)/q_n-\beta$ and $\beta\in \hoo{i}{q_n}$.  Since there are $q_n$ many 
levels and $q_n\gamma=\ceil{q_n\beta}-q_n\beta$ the inequalities in item 1 follow.  Item 2 is similar.  Item 3 follows since 
{$S=\phi_0^{-1}(\bigcup_{j<q_n}[j/q_n, (j+1)/q_n-\gamma)\cup(
\bigcup_{j<q_n}[(j+1)/q_n-\gamma, (j+1)/q_n)\cup M_n$.} The final assertion follows from Lemma \ref{bc1}.
\qed

\medskip

Restating the discussion:

\begin{lemma}\label{bc3}
For almost all $s\in S\subseteq \bk^c$ that are mature at stage $n$,  $\mcs_\beta(s)(0)=s(i)$ where $i\equiv_{q_n}d^n_L$ if $s\in L_n$  and $i\equiv_{q_n}d^n_R$ if $s\in R_n$.
\end{lemma}
\pf Assume that $n$ is mature for $s$. Then on its principal $n$-block, the  {projection of $s$ to 
$\mck_\alpha$} agrees with $w^\alpha_n$.\footnote{Recall $w_n^\alpha$ is the notation for the unique member of the $n^{th}$ element $\mcw^c_n$  of the construction sequence for 
$\mck_\alpha$.} The values $s(0)$ and $\mcs_\beta(s)(0)$ are the $r_n(s)^{th}$ and the $r_n(\mcs_\beta(s)))^{th}$ values of the word $w^\alpha_n$. From equation \ref{def of dn}, $r_n(\mcs_\beta(s)))=r_n(s)+d^n(s)$. Hence $\mcs_\beta(s)(0)=s(d^n(s))$, and the lemma follows.\qed

 The items in the following lemma are essentially Remark  \ref{interval coherence} and Lemma \ref{coherence of rn}  in a different context.
\begin{lemma}\label{coherence of ds}
For almost all $s$ and for $n<m$ that are mature for $s$ and $\mcs_\beta(s)$
\begin{enumerate} 
\item If $i\equiv_{q_n}r_n(s)+d^n(s)$ and $j\equiv_{q_m}r_m(s)+d^m(s)$ then the $j^{th}$ place {in} the principal $m$-block of $\mcs_\beta(s)$ is in the $i^{th}$ place of the principal $n$-block of $\mcs_\beta(s)$. 

\item Let $I$ be the ${r_n(s)+d^n(s)}^{th}$ interval of  $\mci^{q_n}$ and $J$ the ${r_m(s)+d^m(s)}^{th}$ interval of $\mci^{q_m}$ in the dynamical orderings. Then $J\subseteq I$.
\end{enumerate}
\end{lemma}

\pf This follows from Remark \ref{coherence remark} and Lemma \ref{coherence of rn}.
To see this note that $r_n(\mcs_\beta(s))\equiv_{q_n}r_n(s)+d^n(s)$; i.e. $\mcs_\beta(s)(0)$ is in the $i^{th}$ place of the principal $n$-block of $s$ where $i\equiv_{q_n}r_n(s)+d^n(s)$.
\qed

Thus typical points in $R_n$ and $L_n$ are those in which the $n$-block of $\mcs_\beta(s)$ containing $0$ 
is the shift of the block of $s$ containing $0$ by $d^n_R$ and $d^n_L$ respectively.

We now describe how $d^n(sh^k(s))$ changes.  As $k$ varies, $d^n(sh^k(s))$ measures the shift between $sh^k(s)(0)$ and $\mcs_\beta(sh^k(s))(0)$. In regions where the principal $n$-subwords of both $sh^s(s)$ and $\mcs_\beta(sh^k(s))$ exist and are repeating
$d^n(sh^k(s))$ is constant.  It is also constant as it crosses boundary regions of $sh^k(s)$ and $\mcs_\beta(sh^k(s))$ as long as those boundary regions have length $q_n$ and are lined up with adjacent $n$-subwords. However for $m\ge n+1$, if the boundary section of an $m$-word of $s$ or $\mcs_\beta(s)$  has length not divisible by $q_n$, the relative alignment between $s$ and $\mcs_\beta(s)$ changes. This happens on regions of $\bigcup_{m\ge n+1}\boundary_m(s)\cup \bigcup_{m\ge n+1}\boundary(\mcs_\beta(s))$.

If $n$ is mature for $s$, the principal $n$-word of $s$ repeats on both sides of $s(0)$ and thus we see:

\begin{lemma}\label{making shifts constant}
If $s$ is mature at stage $n$, then $d^n(s)$ is constant on the principal $n$-block of $s$. 
Moreover on $d^n(s)$ is constant on the even and odd overlaps of 
2-subsections of $n+1$ subwords of $s$ and {$\mcs_\beta(s)$.}

\end{lemma}

The next lemma is used for the ``nesting" arguments in Section \ref{red zones}. It says that the 
measure of the set of $s\in S$ with $d^n(s)=d^n_L$ or $d^n(s)=d^n_R$ can be closely computed as a 
density in every scale  bigger than $n$.

\paragraph{Remark}{The notation $d^n_L$ and $d^n_R$ are supposed to be suggestive of the left and right lanes. To a close approximation, if $s$ is mature and in a left lane then $d^n(s)=d^n_L$ and similarly for the right lanes.}
\medskip

\begin{lemma}\label{localize}

Let $n<m\in \nn$ be natural numbers. Then $\{0, 1, 2, \dots q_m-1\}=P^n_L\cup U\cup P^n_R$ such that  for almost every $s$ for which $n$ is mature:\footnote{Properly speaking the $P^n_R$ and $P^n_L$ notation should indicate $m$ as well. Without any contextual indication of what $m$ is we take $m=n+1$.	}
\begin{enumerate}
  \item If $r_m(s)\in P^n_L$, then $s\in L_n$,
  \item If $r_m(s)\in P^n_R$ then $s\in R_n$,
  \item $|U|\le 2q_n$,
  \item $\left| {|P^n_L|\over q_m}-\beta^L_m\right|<{2q_n\over q_m}$ and 
  \item $\left| {|P^n_R|\over q_m}-\beta_m^R\right|<{2q_n\over q_m}$.
\end{enumerate}
\end{lemma}

\pf As in Lemma \ref{two values}, let $\gamma={(i+1)\over q_n}-\beta$, where $i=p_nD_n(\beta)$ (See figure \ref{fig left right}). 
The partition $\mci_{q_m}$ splits each interval $I\in \mci_{q_n}$ into ${q_m\over q_n}$ subintervals. Let 
$U$ be the indices of the  $\mci_{q_m}$ intervals that lie over or under $\gamma$ and 
$\gamma+{1\over q_m}$.  Explicitly: suppose that $\gamma\in I^m_{i_0}$ and $\gamma+{1\over q_m}\in I^m_{i_1}$. 
Let 
\begin{eqnarray*}
U=&\{i:\mbox{ for some }0\le j<q_n, I^m_i=\mcr_{\alpha_n}^j I^m_{i_0}\}\cup\\
&\{i:\mbox{ for some }0\le j<q_n, I^m_i=\mcr_{\alpha_n}^j I^m_{i_1}\}.
\end{eqnarray*}
Then $|U|={2q_n}$, and if $i\notin U$, then either:
\begin{eqnarray}
I^m_i & \subseteq & \bigcup_{j<q_n}[j/q_n,(j+1)/q_n-\gamma) \mbox{ or } \label{L}\\
I^m_i & \subseteq & \bigcup_{j<q_n}[(j+1)/q_n-\gamma, (j+1)/q_n)\label{R}
\end{eqnarray}
For $i\notin U$, put $i\in P^n_L$ if it satisfies equation \ref{L} and $i\in P^n_R$ if it satisfies equation \ref{R}. 
It follows that for almost all $s$, if $n$ is mature for $s$ and $r_n(s)\in P^n_L$, then $d^n(s)=d^n_L$ and 
similarly for $P^n_R$. 
Since $P^n_R\cup P^n_L\cup U$ is a partition of $q_m$ and $|U|\le 2q_n$, the lemma follows.
\qed

\begin{lemma}\label{left or right doesnt matter}
Let $f\in \{0, 1\}^\nn$ and $s$ be a typical member of $S(\beta)$.
\begin{enumerate}
  \item Let $\beta^*_n={p_nD_n(\beta)+f(i)\over q_n}$. Then $\la \mcr_{\beta^*_n}:n\in\nn\ra$ converges  to 
  $\mcr_\beta$ in the $C^\infty$-topology.

  As a result, in the language of  symbolic systems:
  \item \label{An} Let $A_n=D_n(\beta+{f(i)\over q_n})$ and $T$ be the shift map on $\mck_\alpha$. Then $A_n$ is either $d^n_L$ or $d^n_R$, depending on the value of $f$ and  for almost every $s\in S, \lim_{n\to \infty}T^{A_n}s=\mcs_\beta(s)$.

  \item\label{pm infinity} With $A_n$ as in item \ref{An} and ${\bk^c}$ an arbitrary circular system with the given coefficient sequence $\la k_n, l_n:n\in\nn\ra$, define $a_n$ and $b_n$ to be the left and right endpoints of the principal $n$-block of $T^{A_n}(s)$.  Then for almost all $s$, $\lim_{n\to \infty}a_n=-\infty$ and $\lim_{n\to \infty}b_n=\infty$.
\end{enumerate}
\end{lemma}

\pf The first item follows because $|\beta^*_n-\beta|<2/q_n$. Hence $\beta^*_n$ converges rapidly to $\beta$. The second item follows from  the first via the isomorphism $\phi_0^{-1}$.  The third item follows because $\mcs_\beta(s)\in S$ and $T^{A_n}(s)$ converges to $\mcs_\beta(s)$ topologically. Hence for all $n$ there is an $N$ such that for all $m\ge N$,  the principal $n$-block of $T^{A_m}(s)$ is the same as the principal $n$-block of $\mcs_\beta(s)$. Since the principal $m$-block of $T^{A_m}$ contains the principal $n$-block of $\mcs_\beta(s)$ and $\mcs_\beta(s)\in S$, item three follows. \qed

If $a_n$ and $b_n$ are as in item \ref{pm infinity}, then:\begin{eqnarray}\label{to infinity}
a_n=-r_n(s)+A_n & \mbox{ and }& b_n=q_n-r_n(s)+A_n.
\end{eqnarray}

\section{The Displacement Function}\label{displacement function}

In this section we define a function $\Delta$ from $S^1$ to the extended positive real numbers that will eventually be 
 shown to have the properties that 
 \begin{itemize}
 \item $\Delta(\beta)<\infty$ implies that there is an element of the centralizer of $\bk^c$  having $\mcr_\beta$ as its associated rotation.
 \item if ${{\bk^c}}$ is built suitably randomly, then every element of the centralizer of ${{\bk^c}}$, or isomorphism from ${{\bk^c}}$ to $({\bk^c})^{-1}$ has rotation factor $\beta$ with $\Delta(\beta)<\infty$.
 \end{itemize}
 The idea behind the displacement function is simple:  the number $\beta$ determines $\mcs_\beta$ and hence a shift at each scale $n$.  The words in 
 $\mcw_{n+1}^c$ are of the from $\mcc(w_0,\dots w_{k_n-1})$.  If the shift at stage $n$ lines up most $n$-words with other $n$-words in the same argument of $\mcc$ then it is possible to build an element of the centralizer of any $\bk^c$ having rotation factor $\beta$.  If not, and we build $\bk^c$ suitably randomly, then we can arrange that $\beta$ is not a central value.
 
 Fix $\beta$ for the rest of this section, and let $T:\mck_\alpha\to \mck_\alpha$ be the shift map. The next lemma says that the principal  $n$-blocks of $T^{d^n(s)}(s)$ and $\mcs_\beta(s)$ are exactly aligned.

 \begin{lemma}\label{coherence of rn and rn+1} Let  $s,s^*\in \mck_\alpha$ be typical and  $n<m$ be mature for both. Define $t^*=T^{d^m(s)-d^n(s)}(s^*)$. 
 Then $t^*(0)$ is in the same position of its principal $n$-block as $s^*(0)$ is in $s^*$'s principal $n$-block.
   In particular, $T^{d^m(s)-d^n(s)}(s^*)$ has its 0 in a  position inside an $n$-word in the construction sequence for some copy of $w^\alpha_n$.
 \end{lemma}
 
 \pf {Since the $n$-blocks of $s^*$ repeat on either side of the principal $n$-block of $s^*$,  and these have length $q_n$,}  it suffices to show that $d^m(s)-d^n(s)\equiv_{q_n}0$.
 Let $t=T^{d^m(s)-d^n(s)}(s)$ and consider the  point $s'=T^{d^m(s)}(s)$. Then $s'(0)$ is in the $(r_m(s)+d^m(s))^{th}$ place in its principal $m
$-block. By Lemma \ref{coherence of ds}, $s'(0)$ is in the $(r_n(s)+d^{n}(s))^{th}$ place in its principal 
$n$-block. Since $t=T^{-d^n(s)}(s')$, the point $t$ has its $0$ in the 
$r_n(s)^{th}$ place of its principal $n$-block.  Hence $r_n(t)=r_n(s)$ and by so by remark \ref{life is easy in the interior} $d^m(s)-d^n(s)\equiv_{q_n}0$.
\qed

At first glance Lemma \ref{coherence of rn and rn+1} looks puzzling as  we are not assuming that any of $d^m(s)=d^m(s^*)$, $d^n(s)=d^n(s^*)$ or $r_n(s)=r_n(s^*)$. However the assertion is a statement about how the $n$-towers sit inside the $n+1$-towers.  For mature $s, s^*$ this nesting repeats on either side of the principal $n$-blocks and hence behaves as in the cyclical approximations.  Thus it is independent of the value of $d^n(s^*), d^m(s^*)$ or $r_n(s^*)$, and simply reflects the cyclical structure.
 \medskip

For a particular $s\in \mck$,  the sequence of shifts $T^{d^n(s)} (s)$  converges to $\mcs_\beta(s)$. Lemma \ref{coherence of rn and rn+1}  tells us that 
this happens promptly: for mature $n$, $T^{d^n(s)}(s)$ has its $0^{th}$ place in the same position of its principal $n$-block as 
$\mcs_\beta(s)$ does.
 \medskip
   
  We now consider the location of 0 in the principal $n+1$-block of the point $T^{d^{n+1}(s)-d^n(s)}(s)$
  relative to the position of $0$ in the principal $n+1$-block of $s$. 
 For some $j_0$ and $j_1$ 
 the principal $n$-block of $T^{d^{n+1}(s)-d^n(s)}(s)$ arises from the $j_0^{th}$ argument of 
 $\mcc(w_n^\alpha, \dots w_n^\alpha)$ and the principal $n$-block of  $s(0)$ is in a position coming from the  
 $j_1^{st}$ argument.

\hypertarget{bmatch}{\begin{definition}} Let $s\in \mck$.
With $j_0$ and $j_1$ as just described, the $j_0^{th}$ argument  of $\mcc(w_n^\alpha, \dots w_n^\alpha)$ \emph{$\beta$-matches} the $j_1^{st}$ argument. The point $s\in \mck$ is \emph{well-$\beta$-matched} at stage $n$ if $s$ is mature at $n$ and $j_0=j_1$. If $n$ is mature for $s$ and $j_0\ne j_1$, then $s$ is \emph{ill-$\beta$-matched}.
\end{definition}

\begin{lemma} \label{convergence to n} Let $\bk^c$ be a circular system and consider $S\subseteq \bk^c$. Let $s, s^*\in S$ and
suppose that $n$ is mature for $\pi(s), \pi(s^*),\mcs_\beta(\pi(s))$ and $\mcs_\beta(\pi(s^*))$ and that $\pi(s)$ is well-$\beta$-matched at stage $n$. Let $A_n=d^n(s)$ and
 $A_{n+1}=d^{n+1}(s)$.  Then:
 \begin{enumerate}
 \item $r_n(T^{A_n}s^*)=r_n(T^{A_{n+1}}s^*)$ and 
 \item if $I$ is the interval $[-r_n(T^{A_n}s^*), q_n-r_n(T^{A_n}s^*))\subseteq \poZ$, then 
 \begin{equation}\label{agreement}
 (T^{A_n}s^*\rest I)=(T^{A_{n+1}}s^*\rest I)
 \end{equation}
 \end{enumerate} 
\end{lemma}

\pf   Lemma \ref{coherence of rn and rn+1} asserts that $0$ is in the same 
place in the principal $n$-block of $T^{A_{n+1}-A_n}(\pi(s^*))(0)$
as $0$ is in the principal $n$-block of $\pi(s^*)$. Since $n$ is mature for $s^*$, the principal $n$-block of $s^*$ is 
repeated on either 
side of $s^*(0)$. Since $n$ is mature for $\mcs_\beta(\pi(s^*))$, the principal $n$-block of $T^{A_{n+1}}s^*$ is 
repeated at least twice on either side of $T^{A_{n+1}}(s^*)(0)$. It follows that $0$ is in the same place in the 
principal $n$-block of 
  $T^{A_n}(T^{A_{n+1}-A_n}(s^*))$ as $0$ is in the principal $n$-block of $T^{A_n}(s^*)(0)$. This proves the first assertion.
  
A repetition of this argument shows the second assertion as well, using the fact that $s$ is well-$\beta$-matched. Indeed 
the definition of \emph{well-$\beta$-matched} implies that the principal $n$-words of $T^{A_{n+1}-A_n}s$ and $s$ are identical. Applying $T^{A_n}$ to both, and using the fact that the principal $n$-words repeat one sees that the principal $n$-words $T^{A_{n+1}}s$ and $T^{A_n}s$ are identical. Since the issue of alignment only involves $\pi(s)$, item 2 holds for all $s^*$ with $\pi(s)=\pi(s^*)$. Moreover, arguing as in the last paragraph using the repetition of the principal $n$-blocks, shifting by an $l<q_n$ does not change this. \qed

\paragraph{Comment}  The terminology in this definition extends easily to general circular systems by saying that $j_0^{th}$ argument and $j_1^{st}$ arguments are $\beta$-matched in $s\in \bk^c$  if and only if this is true in $s^\pi$, where $s^\pi$ is the projection of $s$ to $\mck$. Similarly we write $d^n(s)$ for $d^n(\pi(s))$.


\subsection{The Definition of $\Delta$} \label{def of delta} 
Let $(X,\mcb,\mu, T)=(\bk^c, \mcb, \nu, sh)$ be {a} circular system.
Define 
\begin{equation}\label{delta n beta}
\Delta_n(\beta)=\nu(\{s:s\mbox{ is ill-$\beta$-matched at stage } n\}) 
\end{equation}
 and set\footnote{Since being well or ill-matched only depends on $\pi(s)$ in this section we will not carefully distinguish between $s$ and $\pi(s)$.} 
 \begin{equation}\label{delta beta}
 \Delta(\beta)=\sum_n\Delta_n(\beta).
 \end{equation}
\begin{definition}\label{definition of central}
The number $\beta\in S^1$ is a \emph{central value} iff $\Delta(\beta)<\infty.$
\end{definition}
Note that $\Delta(\beta)$ is  defined using the block structure of the $\mcw^c_n$ and hence is determined 
 by $\beta$ together with  the sequences $\la k_n\ra$ and $\la l_n\ra$. 
Thus for $\beta\in S^1$ the property of being \emph{central} depends only on the circular coefficient sequence $\la k_n, l_n:n\in\nn\ra$, rather than on the particular circular system  $\bk^c$. 

In section \ref{building central}, we show that if $\Delta(\beta)$ is finite then there is an element $T^*$ in the weak closure of $\{T^{n}:n\in \poZ\}$ such that $(T^*)^\pi=\mathcal S_\beta.$ In particular $\beta$ is the rotation factor of an element of the centralizer. That result does not use the results of the rest of this section.

\subsection{Deconstructing $\Delta(\beta)$}\label{deconstruction}

Fix a  $\beta$. 
{Recall that  $\la \varepsilon_n:n\in\nn\ra$ is the  sequence satisfying numerical requirement \ref{varepsilons}: $\varepsilon_N>6\sum_{n>N}\varepsilon_n$.}

{Suppose that $s$ is typical, $n$ is mature and $s$ is ill-$\beta$-matched.  Then there are 4 possibilities:}

\begin{enumerate}
\item  $d^n(s)=d^n_L$ or $d^n_R$ and
\item $d^{n+1}(s)=d^{n+1}_L$ or $d^{n+1}_R$
\end{enumerate}
Call these possibilities $P_{LL}, P_{LR}, P_{RL}, P_{RR}$.  

\begin{lemma}\label{misalignment}Let $n, m\in \nn$ with $n+1<m$.
There is a partition {$\{ P^{n,m}_{{{hd_1},{hd_2}}}:{{hd_1},{hd_2}}\in \{L,R\}\}\cup \{U\}$ of the set $\{0,1,\dots q_{m}-1\}$} such that 
for $s\in S$,
if $n$ is mature for $s$ then 
\begin{enumerate}
\item {$r_{m}(s)\in P^{n,m}_{{{hd_1},{hd_2}}}$ implies $(d^n(s), d^{n+1}(s))=(d^n_{hd_1}, d^{n+1}_{hd_2})$}

\item $|U|\le 2q_n+2q_{n+1}$.
\end{enumerate} 
\end{lemma}
\pf
This follows immediately from Lemma \ref{localize} by holding $m$ fixed and applying the lemma successively to $n$ and $n+1$. Except for a set $U=_{def}U_n\cup U_{n+1}$ that has at most $2q_n+2q_{n+1}$ elements, every point in $\{0, 1, \dots q_m-1\}$ belongs to some $P^n_i\cap P^{n+1}_j$. \qed

The levels of the $q_m$-tower reflect the construction of $w_m^\alpha$ from $n$-words with $n<m$. 
If $s$ and $\mcs_\beta(s)$ are mature at stage $n<m$, then the locations of $s(0)$ and 
$T^{d_{n+1}(s)-d_n(s)}(s)(0)$ in their principal $m$-block and  the pair $(d^n(s), d^{n+1}(s))$ determine whether $s$ is ill-$\beta$-matched or not.
For  particular choices of $hd_1, hd_2\in \{L,R\}$ either all typical $s$ in 
 {$P_{{{hd_1},{hd_2}}}$}  with $n$ mature for both $s$ and $\mcs_\beta(s)$ are well-$\beta$-matched \emph{or} none are.

In the next section we will fix a particular choice of $hd_1$ and $hd_2$. For now let $n, hd_1$ and $hd_2$ be such that  all $n$-mature $s$   in configuration $P_{{hd_1},{hd_2}}$ are ill-$\beta$-matched. We use the symbol 
$\misal_n$  {(In LaTeX: {\small \verb|\not\Downarrow|})} to indicate the \emph{misaligned} points  at stage $n$.
Let

\begin{equation}\label{newdef misal}
\misal_n=\{s:s \mbox{ is ill-$\beta$-matched at stage $n$ and in configuration {$P_{{{hd_1},{hd_2}}}$}}\}.
\end{equation}
We  need to localize the sets $\misal_n$. 
The next  lemma tells us that they are uniformly close to open sets:

\begin{prop}\label{cheating the Ergodic Theorem} Let 
 $n, m\in \nn$ with $n+1<m$.
 Then there is a set $d^{n,m}\subseteq \{0, 1, \dots q_m-1\}$ such that if $s\in S$,  $n$ is mature for $s$ and $r_m(s)+k\in d^{n,m}$, then
 \begin{enumerate}
 \item  $n$ is mature for $sh^k(s)$,
 \item $d^n(sh^k(s))=d^n_{hd_1}$ and $d^{n+1}(sh^k(s))=d^{n+1}_{hd_2}$.
 \item $sh^k(s)\in \misal_n$ and
 \end{enumerate}
 
 \begin{equation*}
\left|{|d^{n,m}|\over{q_m}}-\nu(\misal_n)\right|<
2\left({q_n+q_{n+1}\over q_m} \right)+{1\over l_{n-1}}+\varepsilon_{n-1}.
\end{equation*}
 \end{prop}
\pf Let $s$ be an arbitrary point in $S$ that is mature for $n$.
 Take $d^{n,m}$ to be those numbers of the form  $r_m(s)+k$ (where $k\in [-r_m(s), q_m-r_m(s))$) such that $sh^k(s)$ has its zero point in {$P^{n,m}_{{{hd_1},{hd_2}}}$} and $n$ is mature for $sh^k(s)$. Then $d^{n,m}$ is independent of the choice of $s$. By Lemma \ref {getting old}, the collection of $k$ such that $sh^k(s)$ is not mature for $n$ has density at most ${1\over l_{n-1}}+\varepsilon_{n-1}$.
 \qed

\subsection{Red Zones}\label{red zones}

Suppose that $\beta$ is not central, i.e. that  $\Delta(\beta)=\infty$.
Then for 
some fixed choice of {{$({{hd_1},{hd_2}})$, with $hd_i$ belonging to $ \{L, R\}$}}, 
    \begin{equation*}
    \sum_n\nu(\{s:s\mbox{ is ill-$\beta$-matched at stage $n$ and in configuration {$P_{{{hd_1},{hd_2}}}$}}\})
    \end{equation*}
is infinite. Fix such an {${{hd_1},{hd_2}}$}. Then {with this choice} for all $n, \misal_n$ is well-defined, and moreover there is a set $G\subseteq \nn$ such that if $n<m$ belong to $G$ then $n+2<m$ and 
\begin{equation}
\label{root cause} 
\sum_{n\in G}\nu(\misal_n)=\infty.
\end{equation}
Let $s$ be a  point in $\mck_\alpha$ such that all of the shifts of $s$ and $\mcs_\beta(s)$ are generic with respect to 
basic open sets, the $E^i_n$'s, $\misal_n$, {$P^{n,m}_{{{hd_1},{hd_2}}}$} and the sets $L_n$, $R_n$.
For large enough $M$, we describe how to use $s$ and $\bigcup_{n\in G}\misal_n$  
to identify a subset of the interval $[-r_{M}(s), q_M-r_M(s))$ consisting of misaligned points and having 
{density arbitrarily close to one.}

Assume that $s\in {\misal_n}$ and $n$ is mature for $s$ {and $\mcs_\beta(s)$}.
In defining $\misal_n$, the choice that $(d^n(s), d^{n+1}(s))=(hd1, hd2)$ together with $s(0)$, give us the relative locations of the overlap of the principal $n+1$-blocks of $s$ and $\mcs_\beta(s)$. 

Let $u$ be the principal $n+1$-block of $s$ and $v$ be the principal $n+1$-block of $\mcs_\beta(s)$. and assume that they are in the position determined by $d^{n+1}(s)$.   By Lemma \ref{numerology lemma}, on the overlap the 2-subsections of {$v$} split the 2-subsections of {$u$} into either one or two pieces, and  the positions of all of the even pieces are shifted by the same amount relative to the 2-subsections of $v$ and similarly for the odd pieces. 

We analyze the case where $s(0)$ occurs in an $n+1$-block where the 2-subsections are split into two pieces. If they 
are only split into one piece (i.e. they aren't split) the analysis is similar and easier.  
Without loss of generality we will assume that $s(0)$ occurs in an even overlap. 

Since neither $s(0)$, nor $\mcs_\beta(s)(0)$ occur in the first or last $\varepsilon_nk_n$ 1-subsections of the 
principal 2-subsection that contains them, we know that the overlaps of the principal 2-subsections of $s(0)$ and 
$\mcs_\beta(s)(0)$ contain at least $\varepsilon_nk_n$ 1-subsections.   The 0-subsections of the form $w_j^{{l_n}-1}$ of each 1-subsection of $s$ in this overlap 
are split into at most three pieces, powers of the form $w_i^{s^n_0}$, $w^r_i$ and $w_i^{s^n_1}$ where $0\le r\le 2$, 
$l_n-(s^n_0+s^n_1)\le 3$ and the middle power  $w^r_i$ crosses a boundary section of $\mcs_\beta(s)$.   The powers $s^n_0$ and $s^n_1$ are constant on the overlap of the 
2-subsections,  constant in all of the even pieces of the overlap of the 2-subsections of the principal $n+1$-block, 
and are determined by {$(hd_1, hd_2)$}. Moreover, $s^n_i>\varepsilon_nl_n$. Again, without loss of generality 
we assume that $s(0)$ is in the left overlap corresponding to the power $s^n_0$. 
\bigskip

{\bf \noindent Observation:} There is a number $j_0$ between $0$ and $k_n-1$  determined by the pair $(d^n(s), d^{n+1}(s))$ such that the even piece of a 2-subsection that 
contains $s(0)$ is of the form $\prod_{j<j_0}b^{q_n-j_i}w_j^{l-1}e^{j_i}$, except that the last 1-subsection 
may be truncated. Moreover, since $d^{n+1}(sh^k(s))$ is constant for $k$ in the principal $n+1$-block of $s$,  if 
\begin{equation}\label{counting the overlap}
t=k_n-j_0,
\end{equation} then $t\ne 0$ and  for all $j<j_0$ the powers {$w_j^{s^n_0}$} are \hyperlink{bmatch}{$\beta$-matched}   with $w^{s^n_0}_{j+t}$
except for portions of the first and last power.

In particular,  if $k$ is such that the $0$ position of  $sh^k(s)$ lies in the interior of initial power  $w_j^{s^n_0}$ in an even overlap and $j<j_0$, then $sh^k(s)\in \misal_n$ because it is lined up with $w_{j+t}$.

\bigskip

\begin{lemma}\label{blocking out}
Let $s\in \mck$ and suppose that $s$ and  $\mcs_\beta(s)$ 
 are generic, and that $s$ is mature at $n$.
Suppose that $m> n+2$. Then there is a set $B_n\subseteq \{0, \dots q_{m}-1\}$ such that if
$k\in [-r_m(s), q_m-r_m(s))$ and $r_m(s)+k\in B_n$, then:

\begin{enumerate}
\item $sh^{k}(s)$ has its zero located in $B_n$,
\item $n$ is mature for $sh^k(s)$,
\item $sh^k(s)\in \misal_n$,

\item \label{unions of blocks} There is a $j_0>\varepsilon_nk_n$  and a $t\ne 0$ such that $B_n$ is: 

\begin{enumerate} 
\item a union of sets, each  of  the form $\bigcup_{j<j_0}U_j$
\item each set $\bigcup_{j<j_0}U_j$ is a subset of a position of an occurrence in $s$ of an  $n+1$-subword 
$\mcc(u_0, u_1, \dots u_{k_n-1})$  of $w_m^\alpha$ (with $u_i=w_n^\alpha$), 

\item each $U_j$ is a collection of positions non-$n$-boundary positions in $u_j^{s^n_0}$ such that $u_j^{s^n_0}$ is $\beta$-matched with 
$u_{j+t}^{s^n_0}$, except perhaps for the first or last copy of $u_j$ in $u_j^{s^n_0}$  

and 

\item each set $\bigcup_{j<j_0}U_j$ is the collection of all non-$n$-boundary  positions in $u_j^{s^n_0}$ in a block of the form 
{$\prod_{j<j_0}b^{q_n-j_i}u_j^{l_n-1}e^{j_i}$}.

\end{enumerate}

\end{enumerate} 
and
\begin{equation*}
\left|{B_n\over q_m}-\nu(\misal_n)\right|<2\left({q_n+q_{n+1}\over q_m} \right)+{1\over l_{n-1}}+\varepsilon_{n-1}.
\end{equation*}
\end{lemma}

\pf The first statement is automatic since $B_n\subseteq \{0, 1 \dots q_m-1\}$. Let $d^{n,m}$ be as in Proposition \ref {cheating the Ergodic Theorem}. If $k\in d^{n,m}$ then, as in the 
discussion before the statement of Lemma \ref{blocking out},  $sh^k(s)(0)$ occurs in the position of a power 
$u^{s^n_0}$, where $u$ is the principal $n$-block of $sh^k(s)$ and  $u^{s^n_0}$ occurs on the left overlap 
of 1-subsections of the principal $n+1$-block of $sh^k(s)$.

As in the observation before this lemma,  to each $k\in d^{n,m}$ we  can associate a 
set $\bigcup_{j<j_0}U_j$ containing $k$ by taking all of the positions of the powers  $u_j^{s_0^n}$ in the even 
overlap determined $sh^k(s)(0)$, where $k$ is not in the boundary of a $u_j$. Let $B_n$ be the union of all of the collections $\bigcup_{j<j_0}U_j$  as $k$ 
ranges over $d^{n,m}$. 

Assertion 4(c) follows from the observation and the fact that  $d^n$ and $d^{n+1}$
 are constant (and equal to {$d^n_{hd_1}$ and $d^{n+1}_{hd_2}$}) on $d^{n,m}$. 

We show that if $k'\in B_n$ then $n$ is mature for $sh^{k'}(s)$ and that $sh^{k'}(s)\in \misal_n$.
The maturity of $n$ follows immediately from the maturity of $s$ and the fact  that the location of 0 in $sh^k(s)$ is in a non-boundary portion of an $n$-subword 
 of its principal $n+1$-block. That $sh^{k'}(s)\in \misal_n$ follows from the fact that
 $u_j^{s^n_0}$ is $\beta$-matched with $u_{j+t}^{s^n_0}$, and $t\ne 0$.
 
To finish, note that 
 \[d^{n,m}\subseteq \bigcup \bigcup_{j<j_0}U_j\subseteq \misal_n.\]
 Hence 
 \[{|d^{n,m}|\over q_n}\le {|\bigcup\bigcup_{j<j_0}U_j|\over q_n}\le \nu(\misal_n).\]
 Thus Lemma \ref{blocking out} follows from Lemma \ref{cheating the Ergodic Theorem}.
\qed

We now define the \emph{red zones} corresponding to $\beta$. 
Recall that if $n<m\in G$ then $n+2<m$ and $\sum_{n\in G}{\nu(\misal_n)}=\infty$. For $n<m$ consecutive elements of $G$, define 
\begin{equation*}
\delta_n=4\left({q_{n+1}\over q_m}\right)+{1\over l_{n-1}}+\varepsilon_{n-1}
\end{equation*}
Then we see that:
\begin{itemize}
\item $\sum_{n\in G}\delta_n<\infty$, so
\item $\sum_{n\in G}({\nu(\misal_n)}-\delta_n)=\infty$
\end{itemize}
and if $B_n$ is the set defined in Lemma \ref{blocking out}, then ${\nu(\misal_n)}-\delta_n\le {|B_n|\over q_m}\le{\nu(\misal_n)}$.

\begin{lemma}\label{new red zone} Let $N$ be a natural number and $\delta>0$. Suppose that $s$ and  $\mcs_\beta(s)$ 
are generic, and
that $s$ is mature at $N$. 
Then there is a sequence of natural numbers 
$\la n_i:1\le i\le i^*\ra$, an $M$ and sets $R_i\subseteq \{0, 1, 2 \dots q_M-1\}$, for $1\le i\le i^*$, such that 

\begin{enumerate}
\item $N<n_1$ and $n_i+2<n_{i+1}<M$,\label{basic condition}
\item $R_i$ is disjoint from $R_j$ for $i\ne j$, \label{disj cond}
\item $R_i$ is a union of blocks of the form $B_{n_i}$ described in condition \ref{unions of blocks}  in Lemma \ref{blocking out} inside $n_{i+1}$-subwords of $w^\alpha_M$\label{block condition}
\item if $k\in R_i$, then $sh^k(s)\in \misal_{n_i}$,\label{misal condition}
\item the density of $\bigcup_i R_i$ in $\{0, 1, \dots q_M-1\}$ is at least $1-\delta$. 
\end{enumerate}

\end{lemma}

\pf
We can assume that $N$ is so large that $\bigcup_{n\ge N}\partial_n$ has measure less than $\delta/100$ and $1/l_N+\varepsilon_N<\delta/100$. From the definition of $G$ we can find a collection $\la n_i:i\le i^*\ra$ of consecutive elements of $G$  so that  
\[\prod_{1\le i\le i^*}(1-\nu(\misal_{n_i})+\delta_{n_i})<\delta/100.\]
 Choose an $M>n_{i^*}+2$, and for notation purposes set $n_{i^*+1}=M$. 

Define  sets $R_i$  and $I_i$ by reverse induction from  $i=i^*$ to $i=1$ with the following properties:
\begin{description}

\item[i.] $\{0, 1\dots q_M-1\}\setminus ((\bigcup_{i^*\ge j\ge i}I_j)\cup (\bigcup_{i^*\ge j\ge i}R_j))$ consists of entire locations of  words $w^\alpha_{n_{i}}$ in $w^\alpha_M$,
\item[ii.] $R_i\subseteq \{0, 1\dots q_M-1\}\setminus ((\bigcup_{i^*\ge j>i}I_j)\cup (\bigcup_{i^*\ge j>i}R_j))$ and has relative density at least ${\nu(\misal_{n_i})}-\delta_{n_i}$,  
\item[iii.] the set $I_i\subseteq \bigcup_{j=n_i+1}^{n_{i+1}}\partial_j\cap \{0, 1\dots ,q_M-1\}$ and hence,

\item[iv.] $I_i$ has density less than or equal $1/l_{n_{i}}$ in $\{0, 1\dots q_M-1\}$

\end{description}

To start, apply Lemma \ref{blocking out} with $m=n_{i^*+1}$, to get a set $B_{n_{i^*}}\subseteq \{0, 1, \dots q_M-1\}$ of density at least $\nu(\misal_{n_{i^*}})-\delta_{n_{i^*}}$ 
satisfying conditions \ref{block condition}-\ref{misal condition} of the lemma we are proving.
Set $R_{i^*}=B_{n_{i^*}}$. Let $I_{i^*}=\bigcup^{M}_{j=n_{i^*}+1}\partial_j\cap \{0, 1\dots, q_M-1\}$.

Suppose that we have defined $\la R_j: i^*\ge j>i\ra$ and $\la I_j:i^*\ge j>i\ra$  satisfying the induction hypothesis (i-iv).

Apply Lemma \ref{blocking out}  again to get a set $B=B_{n_i}$ a subset of $ \{0, 1, \dots q_{n_{i+1}}-1\}$.
Inside each copy $\{k, k+1, \dots ,k+q_{n_{i+1}}-1\}$ corresponding to a location in $w^\alpha_M$ of a $w^\alpha_{n_{i+1}}$ in the complement of $((\bigcup_{i^*\ge j>i}I_j)\cup (\bigcup_{i^*\ge j>i}R_j))$, we have a translated copy of $B$,  $k+B$. Let $R_{i}$ be the union of the sets $k+B$ where $k$ runs over the locations the words $w^\alpha_{n_{i+1}}$ in the complement of $((\bigcup_{i^*\ge j>i}I_j)\cup (\bigcup_{i^*\ge j>i}R_j))$.

Then the density of $R_i$ relative to
\[\{0, 1, \dots q_M-1\}\setminus ((\bigcup_{i^*\ge j>i}I_j)\cup (\bigcup_{i^*\ge j>i}R_j))\]
is at least ${\nu(\misal_{n_i})}-\delta_{n_i}$. It follows from conclusion 3 of lemma \ref{blocking out} that  $R_i$ is a union of non-boundary portions of blocks of length $q_{n_i}^{s^0_{n_i}-1}$ corresponding to  positions of
$w^\alpha_{n_i}$ in $w^\alpha_M$ ,

Since $R_i$ consists of a union of the non-boundary portion of locations of words $w^\alpha_{n_i}$, 
	\[\{0, 1, \dots q_M-1\}\setminus \left((\bigcup_{i^*\ge j>i}I_j)\cup (\bigcup_{i^*\ge j>i}R_j))\cup R_i\right)\]
consists of the entire blocks of locations of $w^\alpha_{n_i}$ together with elements of $\bigcup_{j=n_i}^{n_{i+1}}\partial_j$. The latter set has density less than or equal to $1/l_{n_i-1}$. Let 
	\[I_i=\left(\{0, 1, \dots q_M-1\}\cap \bigcup_{j=n_i}^{n_{i+1}}\partial_j\right) \setminus \left((\bigcup_{i^*\ge j>i}I_j)\cup (\bigcup_{i^*\ge j>i}R_j))\cup R_i\right).\]

It remains is to calculate the density of $\bigcup_{1\le i\le i^*}R_i$. At each step in the induction, we remove a portion of density at least $\nu(\misal_{n_i})-\delta_{n_i}$ from $\{0, 1, \dots q_M-1\}\setminus ((\bigcup_{i^*\ge j>i}I_j)\cup (\bigcup_{i^*\ge j>i}R_j))$. Let $\partial=\bigcup_{1\le i\le M}\partial_{n_i}$. Then the density of 
the union of the $R_i$'s is 
at least
\begin{equation*}
1-\prod_{i^*\ge i\ge 1}(1-\misal_{n_i})-|\partial|/q_m,
\end{equation*}
which is at least  
$1-\delta.$    \qed

\section{The Centralizer and Central Values} \label{centralizer and central}

 In the first part of this section we show that every central value is rotation factor of an element of the closure of the powers of $T$ and hence an element of the centralizer.

The second part shows a converse: if $\bk^c$ is built sufficiently randomly then the rotation factor of every element of the centralizer is a rotation by a central value.

\medskip

We note in passing that every circular system is \emph{rigid}: if $s$ is mature for $n$, then $T^{q_n(l_n-2)}(s)$ has the same principal $n$-block as $s$ does. It follows that $\overline{\{T^n:n\in\poZ\}}$ is a perfect Polish monothetic group.

\subsection{Building Elements of the Centralizer} \label{building central}

If $\Delta(\beta)$ is finite, then the Borel-Cantelli lemma implies that for $\nu$-almost every 
$s$,  there is an $n_0$ such 
that for all $n\ge n_0$, $s$ is well-$\beta$-matched at stage $n$. As a consequence, certain sequences of translations  converge. Precisely:

\begin{theorem}\label{central values in closure}
Suppose that ${{\bk^c}}$ is a uniform circular system with coefficient sequence $\la k_n, l_n:n\in\nn\ra$. Let $T$ 
be the shift map on ${{\bk^c}}$ and $\beta\in \zoo$ be a number such that $\Delta(\beta)<\infty$. Then there is a 
sequence of integers $\la A_n:n\in\nn\ra$ such that $\la T^{A_n}:n\in\nn\ra$ converges  pointwise almost 
everywhere to a $T^*\in C(T)$ with $(T^*)^\pi=\mcs_\beta$. In particular there is a sequence $\la A_n:n\in\nn
\ra$ such that $\la T^{A_n}:n\in\nn\ra$  converges in the weak topology to a $T^*$ with $(T^*)^\pi=\mcs_\beta$.
\end{theorem}

\begin{corollary}\label{central are central}
If 
$\beta$ is central, then there is a $\phi\in \overline{\{T^n:n\in\poZ\}}$ such that $\phi^\pi=\mcs_\beta$. 
\end{corollary}

\pf Let $\mct$ be the tree of finite sequences $\sigma\in \{L,R\}^{<\infty}$. Choose an $n_0$ such that
 \[G=\{s:n_0\mbox{ is mature for $s$ and for all }m\ge n_0, s\mbox{ is well-$\beta$-matched at stage }m\}\]
  has positive measure. 
By the K\"onig Infinity Lemma there is a function $f:\{m:m\ge n_0\}\to \{L,R\}$ such that for all $m\ge n_0$, $\{s
\in G:d^n(s)=d^n_{f(n)}$ for all $n$ with $n_0\le n\le m\}$ has positive measure.  Let $A_n=d^n_{f(n)}$.

By Lemma \ref{left or right doesnt matter}, item \ref{pm infinity} it follows that for a typical $s$  the left and 
right endpoints of the principal $n$-blocks of $T^{A_n}s$ go to negative and positive infinity respectively. 
Let $s^*$ be a typical element of $S$; e.g. $\pi(s^*)$ and $\mcs_\beta(\pi(s^*))$ both belong to $S^\pi$,  
large  enough $n$ are mature for $s^*$ and for all large $n$, $\pi(s^*)$ is well-$\beta$-matched at stage $n$. Then for all large $n$,    the left and right endpoints of 
the principal $n$-block of $T^{A_n}s$ and $T^{A_{n+1}}s$ are the same. If $s^*$ is well-$\beta$-matched at stage $n$, then   the words 
constituting principal $n$-block of $T^{A_n}s$ and $T^{A_{n+1}}s$ are the same. It follows that for typical $s^*\in S$, the sequence $T^{A_n}s^*$ converges in the 
product topology on $(\Sigma\cup \{b, e\})^\poZ$.

 We now show that the map $s\mapsto \lim T^{A_n}s$ is one-to-one. If $s\ne 
s'$, then either $\pi(s)\ne\pi(s')$ or there is an $N$ such that for all $n\ge N$ the principal $n$-blocks of $s$ and $s'$ differ. We can assume that this $N$ is 
so large that $n$ is mature and well-$\beta$-matched for $\pi(s), \pi(s')$.

If $\pi(s)\ne \pi(s')$, then $\mcs_\beta(\pi(s))\ne \mcs_\beta(\pi(s'))$. Hence the limits of $T^{A_n}s$ 
and $T^{A_n}s'$ differ. So assume that $\pi(s)=\pi(s')$. Then, since $T^{A_n}$ is a translation by at most 
$q_n-1$ and $n$ is mature for all parties (so the principal $n$-blocks of $T^{A_n}s$  and $T^{A_{n}}s'$ 
repeat) we know that the principal $n$-blocks of $T^{A_n}s$ and $T^{A_{n}}s'$ differ. But for all $m>n$, 
the principal $n$-blocks of $T^{A_m}s$ agree with the principal $n$-blocks of $T^{A_n}s$ (and similarly for 
$s'$). {Hence for all $m>N$ the principal $N$-blocks of $T^{A_m}s$ and $T^{A_m}s'$ differ. It follows that the limit map is one-to-one.}

We need to see that for almost all $s, \lim_{n\to \infty}T^{A_n}s$ belongs to 
${{\bk^c}}$. By definition of ${{\bk^c}}$ this is equivalent to showing that for almost all $s$ if $I\subseteq \poZ$ is an 
interval, 
then $\lim_{n\to \infty}T^{A_n}s\rest I$ is a subword of some $w\in \mcw^c_m$ for some $m$.   However, by 
Lemma \ref{convergence to n},  for almost all $s$ we can find an $n$ so large that:
\begin{enumerate} 
\item $I\subseteq [-r_n(s), q_n-r_n(s))$,
\item $T^{A_n}s$ and  $\lim_{n\to \infty}T^{A_n}s$ agree on the location of the principal $n$-block of containing $I$, and
\item $T^{A_n}s$ and  $\lim_{n\to \infty}T^{A_n}s$ agree on what word lies on the principal $n$-block.  
\end{enumerate}Since the principal $n$-block of $T^{A_n}s$ belongs to $\mcw^c_n$, we are 
done.

Summarizing, if $T^*=\lim_{n\to \infty}T^{A_n}s$, then for almost all $s$, $T^*s$ is defined and belongs to $S
$. Moreover $T^*$ is one-to-one and commutes with the shift map.

Define a measure $\nu^*$ on $S$ by setting $\nu^*(A)=\nu((T^*)^{-1}A)$. Then $\nu^*$ is a non-atomic, shift 
invariant measure on $S$. By  Lemma \ref{dealing with S}, we must have $\nu^*=\nu$. In 
particular we have shown that $T^*:{{\bk^c}}\to {{\bk^c}}$ is an invertible measure preserving transformation belonging 
to $\overline{\{T^n:n\in \poZ\}}$, with $(T^*)^\pi=\mcs_\beta$.
\qed

We make the following remark without proof as it is not needed in the sequel:

\begin{remark}
Suppose that  ${{\bk^c}}$ satisfies the hypothesis of Theorem \ref{central values in closure} and $\beta$ is a central value. Then for any sequence of natural numbers $\la A_n:n\in\nn\ra$ such that $A_n\alpha$ converges to $\beta$ sufficiently fast, the sequence $\la T^{A_n}:n\in\nn\ra$ converges to a $T^*\in C(T)$ with $(T^*)^\pi=\mcs_\beta$.
\end{remark}

\subsection{Characterizing Central Values}

The main result of this section is a converse of Corollary \ref{central are central}. If  ${\bk^c}$ is a circular system built from sufficiently random collections of 
words and $\phi$ is an isomorphism between ${\bk^c}$ and ${\bk^c}$ then $\phi^\pi=\mcs_\beta$ for some central $\beta$. Moreover,  if $\phi$ is an isomorphism between ${\bk^c}$ and  $(\bk^c)^{-1}$ then $\phi^\pi$ is of the form 
$\rev{} \circ\natural\circ\mcs_\beta$ for some central $\beta$.  

{In this section we will return to considering $(\bk^c)^{-1}$ as $(\rev{\bk^c}, sh)$ with the forward shift, and hence can use $\natural$ instead of $\rev{}\circ \natural$.}

\subsubsection{The Timing Assumptions}\label{intro timing}

Randomness assumptions about the words in the $\mcw^c_n$'s  will 
allow us to assert that that the rotations associated with elements of the centralizer of ${\bk^c}$ or  isomorphisms 
between ${\bk^c}$ and $(\bk^c)^{-1}$ arise from central $\beta$'s.
The last part of the paper shows that these additional randomness assumptions are consistent with the randomness assumptions used in \cite{FRW} and describes how to build words with both collections of specifications.

Recall from Definition \ref{circular definition}, that to specify a circular system with coefficient sequence 
$\la k_n, l_n:n\in\nn\ra$ it suffices to inductively specify  collections of prewords 
$P_{n+1}\subseteq (\mcw^c_n)^{k_n}$, and define $\mcw^c_{n+1}$ as the collection 
of words:
\[\{\mcc(w_0, \dots w_{k_n-1}): w_0w_1\dots w_{k_n-1}\in P_{n+1}\}.\]

In the construction, there will be an  equivalence relation $\mcq^1_1$ on $\mcw^c_1$ that is lifted from an analogous equivalence relation on the first step of the odometer construction $\mcw_1$. It is built in section \ref{the specs};  we describe its properties here. Let
$\la \mcq_1^n:n\in\nn\ra$ be the
sequence of propagations of $\mcq^1_1$. As the construction progresses there are  groups $G_1^n$ acting freely on the 
set of $\mcq^n_1$ equivalence classes of words in $\mcw^c_n$. Each $G_1^n$ is a finite sum of copies of
$\poZ_2$. 
Inductively, $G_1^{n+1}=G_1^n$ or $G_1^{n+1}=G_1^{n}\oplus \poZ_2$.  The action of 
$G^n_1$ on $\mcw^c_{n+1}$ arising from the $G^{n+1}_1$ action via the inclusion map of $G^n_1$ into $G^{n+1}_1$ is the skew-diagonal action. We will write 
$[w]_1$ for the $\mcq_1^n$-equivalence class of a $w\in \mcw^c_n$ and $G^n_1[w]_1$ for the orbit of 
$[w]_1$ under $G^n_1$.  If $w\in\mcw^c_{n+1}$ and $C\in \mcw^c_n/\mcq^n_1$ then we say that $C$ 
\emph{occurs at $t$} if there is a $v\in \mcw^c_n$ sitting on the interval $[t, t+q_n)$ inside $w$ and $C=[v]_1$.

\medskip

\begin{numreq}\label{Gs and Qs}
$\sum{|G^n_1|\over |\mcq^n_1|}<\infty$. This can be satisfied by taking ${|G^n_1|\over |\mcq^n_1|}<2^{-n}$.
\end{numreq}

We note that $G^n_1$ is determined directly by the first $n$-nodes in  tree we are using in the domain of the reduction, and hence $|G^n_1|$ is determined by the tree.  So this requirement on $|\mcq^n_1|$ does not depend on any of the other variables being chosen during the construction. In what follows we call such requirements \emph{absolute} requirements.

\medskip
\bfni{Notation:} As an aid to tracking corresponding variables,  script letters are used for sets and non-script Roman letters for the corresponding cardinalities. For example we will use $\mcq_n$ for an equivalence relation and $Q_n$ for the number of classes in that equivalence relation.
\medskip

Here are the the assumptions used to prove the converse to Corollary \ref{central are central}.
The first three assumptions follow immediately from the definitions in section \ref{propagano}.

\begin{enumerate}
\item[T1] The equivalence relation $\mcq_1^{n+1}$ is the equivalence relation on $\mcw^c_{n+1}$ 
propagated from $\mcq_1^n$.
\item[T2] $G^n_1$ acts freely on $\wcuprev{n}{1}$
\item[T3] The canonical generators of $G^n_1$ send elements of $\mcw^c_n/\mcq^n_1$ to elements 
of 
$\rev{\mcw^c_n/\mcq^n_1}$ and vice versa.
\end{enumerate}

The next axiom states that the $\mcq^n_1$ classes are widely separated from each other.

\begin{enumerate}

\item[T4]\hypertarget{T4}{}  \label{T4}
There is a $\gamma$ such that $0<\gamma<1/4$ such that for each {$n$ and each} pair $w_0, w_1\in \mcw^c_n\cup\rev{\mcw^c_n}$
 and each $j\ge q_n/2$
  if
$[w_0]_1\ne [w_1]_1$, then:
\begin{eqnarray*}
 \dbar(w_0\rest[0, j), w_1\rest [0, j))&\ge& {\gamma},\\
\dbar(w_0\rest [q_n-j,q_n), w_1\rest [q_n-j,q_n))&\ge& {\gamma}\\
\mbox{ and }&&\\
\dbar(w_0\rest [0,j), w_1\rest[q_n-j, q_n))&\ge& {\gamma}.
\end{eqnarray*}

\end{enumerate}
\begin{remark}\label{what approx means}
In the axioms $T5-T7$ we write $|x_n|\approx {1\over y_n}$ to mean that $||x_n|-{1\over y_n}|<\mu_n$ 
where $\mu_n\ll \min(\varepsilon_n, 1/Q^n_1)$.
\end{remark}

\begin{numreq}\label{mu_n}
$\mu_n$ is chosen small  relative to $\min(\varepsilon_n, 1/Q^n_1)$. {Explicitly: if $t_n=\min(\varepsilon_n, 1/Q^n_1)$ then $0<\mu_n<t_n\min_{k\le n}2^{-n-2}{1\over t_k}$.}
\end{numreq}
In the next assumption we count  the occurrences of particular $n$-word $v$ that are lined up  in an $n+1$-preword $w_0$ with the occurrences of a particular $\mcq^n_1$-class in the shift of another $n+1$-preword $w_1$ or its reverse. The shift (by $t$ n-subwords), must be non-zero and be such that there is a non-trivial overlap after the shift.

\begin{enumerate}
\item[T5] \hypertarget{T5}{}Let $w_0, w_1$ be prewords in $P_{n+1}$, and $w_1'$ be either $w_1$ or $\rev{w_1}$. 
Write $w_0=v_0v_1\dots v_{k_n-1}$ and  $w'_1=u_0u_1\dots u_{k_n-1}$, with $u_i, v_j\in \mcw^c_n\cup\rev{\mcw^c_n}$.     Let $\mcc\in \mcw^c_{n}/\mcq^{n}_1$ or $\mcc\in \rev{\mcw^c_n}/\mcq^n_1$ according to whether $w_1'=w_1$ or $w_1'=\rev{w_1}$.   For all integers $t$ with
$1\le t\le(1-\varepsilon_{n}) (k_{n})$, $v\in \mcw^c_n$ :
\begin{enumerate}
\item[T5a] (This is comparing $w_0$ with $sh^{tq_n}(w_1')$.)  \hypertarget{T5a}{}
Let
\[J(v)=\{k<k_{n}-t:v=v_k\}.\]
Then
\[
{|\{k\in J(v):u_{t+k}\in \mcc|\over |J(v)|}\approx{1\over Q^n_1}.\]

\item[T5b] (This is comparing $sh^{tq_n}(w_0)$\hypertarget{T5b}{}
 with $w_1'$.)
Let
 \[J(v)=\{k: t\le k\le k_{n}-1\mbox{ and } v=v_k\}.\]
Then
\[{|\{k\in J(v):u_{t-k}\in \mcc
|\over |J(v)|}\approx{1\over Q^n_1}.\]

\end{enumerate}

\item[T6] \hypertarget{T6}{} Suppose that $w_0w_1\dots w_{k_n-1}, w'_0w'_1\dots w'_{k_n-1}\in P_{n+1}$ are  prewords, $1\le t\le (1-\varepsilon_n)k_n$ and $\varepsilon_nk_n\le j_0\le k_n-t$. 
Let 
\[S=\{k<j_0:\mbox{for some }g\in G_1^n, g[w_k]_1=[w'_{k+t}]_1\}.\] Then:
\[{|S|\over j_0}\approx |G^n_1|/Q^n_1.\]

\item[T7] \hypertarget{T7}{} Let $w_0, w_1$ be prewords in $P_{n+1}$, and 
$w_1'$ be either $w_1$ or  $\rev{w_1}$. Suppose that $[w_1']_1\notin G_1^n[w_0]_1$.
Write $w_0=v_0v_1\dots v_{k_n-1}$ and $w'_1=u_0u_1 \dots u_{k_n-1}$, with $u_i, v_j\in \mcw^c_n\cup\rev{\mcw^c_n}$.
  Let $\mcc\in \mcw_n^c/\mcq^n_1$ or
$\mcc\in \rev{\mcw_n^c}/\mcq^n_1$ according to whether $w_1'=w_1$ or $w_1=\rev{w_1}$. Then for all $v\in \mcw^c_n$ if 
\[J(v)=\{t:v_t=v\}\]
then 
\begin{equation}\label{estimate of difference}
{|\{t\in J(v):u_t\in\mcc\} |\over |J(v)|}\approx{1\over Q^n_1}
\end{equation}

\end{enumerate}

\begin{definition}
We will call the collection of  axioms T1-T7 the 
\emph{timing assumptions} for a  construction sequence and an equivalence relation $\mcq^1_1$.
\end{definition}

\subsubsection{Codes and $\dbar$-Distance}

We now prove some lemmas about $\dbar$.\footnote{Basic notation and facts about stationary codes are reviewed in section \ref{random label}.}

\begin{lemma}\label{no code}
Let $w_0\in \mcw_{n+1}^c, w_1\in \mcw_{n+1}^c\cup \rev{\mcw_{n+1}^c}$ and $[w_0]_1\notin G^n_1[w_1]_1$. Let $r>1000$ and $J_0, J_1$ be intervals in $\poZ$ of length $r*q_{n+1}$. Let $I$ be the intersection of the two intervals. Put $w_0^r$ on $J_0$ and $w_1^r$ on $J_1$ and suppose that all but (possibly) the first or last copies of $w_0$ are included in $I$. Let $\bar{\Lambda}$ be a stationary code such that the length of $\Lambda$ is less than $q_n/10000$. 
Then:
\begin{equation}\label{tres bizarre}
\dbar(\bar{\Lambda}[w_0^r\rest I], w_1^r\rest I)> {1\over 50}(1-{1\over Q^n_1}){\gamma}.
\end{equation}
\end{lemma}
\pf
Since the length of the code $\Lambda$ is much smaller than $q_n$ and $r>10000$, the end effects of $\Lambda$ are limited to the first and last copies of $w_0$ and thus affect at most $(1/5000)$ proportion of  $\dbar{(\bar{\Lambda}[w_0^r\rest I], w_1^r\rest I)}$. Removing the portion of $I$ across from the first or last copy of $w_0$ leaves a segment of $I$ of proportion at least $4999/5000$.

For all of the copies of $w_0$, except perhaps at most one at the end of $J_0$, there is a corresponding copy of $w_1$  that overlaps $w_0$ in a section of at least $q_{n+1}/2$. Discard the portions of $I$ arising from copies of $w_0$ not overlapping the corresponding copies of $w_1$.   After the first two removals we have a portion of $I$ of proportion at least $(1/2)(4999/5000)$. 
 
 Because $w_0$ and $w_1$ have the same lengths, the relative alignment between any two corresponding copies of $w_0$ and $w_1$ in the powers $w_0^r$ and $w_1^r$ are the same.  In particular, the ``even overlaps" and ``odd overlaps" are the same in each remaining portion of the corresponding copies of $w_0$ and $w_1$. 
 
 By Lemma \ref{numerology lemma}, there are $s, t<q_n$ such that on the even overlaps  all of the  
 $n$-subwords
 of $sh^s(w_0^r)$ are either lined up with an $n$-subword of $w_1^r$ or with a boundary 
 section of $w_1$, and  all of the $n$-subwords of $w_0$ in an odd overlap are lined up with an 
 $n$-subword or a boundary section of $w_1^r$ by $sh^t(w_0^r)$.
 
Either the even overlaps or the odd overlaps  contain at least $1/2$ of the $n$-subwords that are not across from boundary portions of $w_1$. Assume that $1/2$ of the $n$-subwords lie in even overlaps and discard the portion of $I$ on the odd overlaps. (If more than $1/2$ of the $n$-subwords are in odd overlaps we would focus on those.)

 Let $(w_0^*)^r=sh^s(w_0^r)$ on the even overlaps. 
 Denote any particular copy of $w_0$ in 
$(w_0^*)^r$ as $w_0^*$. Then, except for 
 $\mcw^c_n$-words that get lined up with a boundary section of $w_1$, every $n$-subword of $(w_0^*)^r$ coming from  an even 
overlap of $(w_0)^r$ gets lined up with an $n$-subword of $(w_1)^r$. 
Write $w_0=\mcc(v_1, v_2, \dots v_{k_n-1})$ and
 $w_1=\mcc(u_1, u_2, \dots u_{k_n-1})$ (or, respectively, $w_1=\mcc^r(\rev{u_1},\rev{u_2}, \dots \rev{u_{k_n-1}})$). Then each $n$-subword of 
 $w_0^*$  coming from an even overlap is of the form $v_i$ for some $i$. 
 There is a $t$ such that for all $i$ if $v_i$ occurs in any copy of $w_0^*$  and comes from an even overlap 
 then either:
\begin{enumerate}
\item[a.)] $v_i$ is lined up with $u_{i+t}$ (respectively $\rev{u_{k_n-(i+t)-1}}$) or
\item[b.)] $v_i$ is lined up with a boundary portion of $w_1$ or
\item[c.)] $v_i$ is lined up with $u_{i+t+1}$ (respectively $\rev{u_{k_n-(i+t+1)-1}}$).
\end{enumerate}

On copies of $v_i$ coming from even overlaps of 2-subsections the powers of $v_i$ in alternatives a.) and c.) are constant. Since the even overlaps of the 2-subwords has size at least half of the lengths of the 2-subwords, it follows that $0\le t\le k_n/2$.

Since all of $v_i^{l_{n}-1}$ satisfies a.), b.), or c.), after discarding the  $v_i$'s in case b.) half of the remaining $v_i$'s satisfy a.) or c.). Keep the larger alternative and discard the other. What is left after all of the trimming has size at least:
\[(4999/5000)(1/2)(1/2)(1-2|\boundary_{n+1}|)>1/10\]
proportion of $I$.

For some $t$ 
 what remains consists of $n$-subwords $v_i$ in even overlaps of $(w_0)^r$ that, 
after being shifted by $s$ to be subwords of $(w_0^*)^r$, are aligned with occurrences of $n$-subwords of 
$(w_1)^r$ of the form $u_{i+t}$ ($\rev{u_{k_n-(i+t)-1}}$ respectively). For the rest of this proof of 
Lemma \ref{no code} we will call these the \emph{good 
occurrences} of $n$-subwords.

 \bigskip
{\bf \noindent Claim:} Suppose that $v\in \mcw^c_n$
 and let
\[J^*(v)=\{y\in I:y\mbox{ is at the beginning of a good occurrence of $v$ in } (w_0^*)^r\}.\]
 Let $\mcc\in \mcw^c_n/\mcq^n_1$ or $\mcc\in \rev{\mcw^c_n}/\mcq^n_1$ depending on whether $w_1\in \mcw^c_{n+1}$ or $w_1\in \rev{\mcw^c_{n+1}}$. Then
\begin{equation}\label{claim equation}
\left|{|\{y\in J^*(v):\mbox{some element of }\mcc \mbox{ occurs at }y\mbox{ in $w_1$} \}|\over |J^*(v)|}-{1\over Q^n_1}\right|
\end{equation}
is bounded by  $2/q_n+2/l_n+\mu_n$.
\bigskip

We prove the claim. We have two cases:

\medskip

{\bf\noindent Case 1:} $t=0$.
\medskip

 In this case we have a \emph{trivial split} in the language of section \ref{understanding the words}.
The overlap of the 2-subsections contains the whole of the two subsections except for a portion of one 1-subsection. 
Since $[w_0]_1\notin G_1^n[w_1]_1$ we can apply axiom \hyperlink{T7}{T7} to the words $w_0$ and $w_1$.  
The claim follows from inequality \ref{estimate of difference}, which  is the preword version of formula \ref{claim equation}, after taking into account the boundary and  the words at the ends of the blocks of $(w_0^*)^r$ and the truncated 1-subsections.

\medskip

 \medskip
{\bf\noindent Case 2:} $t\ne 0$.
\medskip

In this case the split is non-trivial. Because the even overlaps are at least as big as the odd overlaps of 2-subsections, the 
even overlap looks like: 
\[\prod_{j=0}^{t^*}(b^{q-j_i}v_j^{l-1}e^{j_i})\]
but with a portion of its last 1-subsection possibly truncated. In particular it has an initial segment of the form 
\[\prod_{j=0}^{t^*-1}(b^{q-j_i}v_j^{l-1}e^{j_i})\]
where  $t^*\ge k_n/2$.

 It follows from  the timing assumption \hyperlink{T5}{T5} that if $J'= \{y\in J(v):$ some element of $\mcc $ occurs across from a word starting at $y$ in the first $t^*-1$ 1-subsections$\}$ then
\[\left|{|J'|\over |J(v)|}-{1\over Q^n_1}\right|<\mu_n.\]

 Any variation  between the quantity in formula \ref{claim equation} and the estimate in 
\hyperlink{T5}{T5} is due to the portion of the last 1-subsection of the even overlaps. This takes up a proportion of the remaining  even 
overlap less than or equal to $1/t^*\le 2/q_n$.  This proves the Claim.\footnote{The axiom \hyperlink{T5b}{T5b} takes care of the case where the relevant overlaps is odd.}
\bigskip

We now shift $(w_0^*)^r$ back to be $w_0^r$ and consider $s$. 
There is an $l'\ge l/2-1\ge l/3$ such that
all of the good occurrences of a $v\in \mcw^c_n$ in $(w_0^*)^r$ are in a power $v^{l'}$. Depending on whether $s\le q_n/2$ or $s>q_n/2$, for each good occurrence of a $v_j$ in $(w_0^*)^r$  either:
\begin{enumerate}
\item[a.)] there are at least $l'-1$ powers of $v_j$ in the corresponding occurrence in $w_0$ such that their left overlap with $u_{j+t}$ has length at least $q_n/2$ 

or 
\item[b.)] there are at least $l'-1$ powers of $v_j$ in the corresponding occurrence in $w_0$ such that their right overlap with $u_{j+t}$ has length at least $q_n/2$ 
\end{enumerate}

Without loss of generality we assume alternative a.). Suppose that the overlap has length $o$ in all of the good occurrences. Then the left side of $v_j$ overlaps the right side of $u_{j+t}$ by at least $q_n/2$. 

By axiom \hyperlink{T4}{T4}, if $v\in \mcw^c_n$, 
\[\dbar(\bar\Lambda[(v\rest[0,o)], u_{j+t}\rest[q_n-o-1, q_n))<{\gamma}/2\] and
 \[\dbar(\bar\Lambda(v\rest[0,o), u_{j'+t}\rest[q_n-o-1, q_n))<{\gamma}/2\]
 then $[u_{j+t}]_1=[u_{j'+t}]_1$.  It follows that if we fix a $v\in \mcw^c_n$ and let 
 \[J(v)=\{j:v_j=v\}\]
  then
\begin{equation*}
{|\{j\in J(v):\dbar(c(v_j\rest[0,o), u_{j+t}\rest[q_n-o-1, q_n-1))<{\gamma}/2\}|\over |J(v)|}
\end{equation*} 
is less than ${1\over Q^n_1}+\mu_n$.

Since at least $1/20$ proportion of $I$ consists of left halves of good occurrences of  the various $v$'s belonging to $\mcw^c_n$ it follows that 
\begin{equation}\label{lower bound for dbar}
\dbar(\bar{\Lambda}[w_0^r\rest I], w_1^r)\ge {1\over 20}(1-{1\over Q^n_1}-\mu_n){(\gamma/2)}.
\end{equation}
The lemma follows.\qed
\subsubsection{Elements of the Centralizer}

In this section we prove the theorem linking central values to elements of the centralizer of ${\bk^c}$.

\begin{theorem}\label{centralizer central values}
Suppose that $({\bk^c},\mcb, \nu, sh)$ is a  circular system built from a circular construction sequence satisfying the timing assumptions. Let $\phi:{\bk^c}\to {\bk^c}$ be an automorphism of $({\bk^c},\mcb, \nu, sh)$. Then $\phi^\pi=\mcs_\beta$ for some central value $\beta$.
\end{theorem}

{This is a condition that does not involve any of the other variables being chosen: at the moment when $Q^n_1$ is being chosen $|G^n_1|$ is already determined by the tree $\mct$ in the domain of the reduction.}

\pf Fix a $\phi$ and suppose  that $\phi^\pi=\mcs_\beta$. We must show that $\beta$ is central. Suppose not.
The idea of the proof is to choose a stationary code $\overline{\Lambda^*}$ well approximating $\phi$ and an $N$ such such for all $M>N$,  passing over the principal $M$-block of most $s\in {\bk^c}$ with $\overline{\Lambda^*}$ gives a  string 
very close to $\phi(s)$ in $\dbar$-distance. Consider  an $s$ where $\overline{\Lambda^*}$ codes well on this  principal $M$-block. 
  
  Use Lemma \ref{new red zone} to build a red zone corresponding to $M$.  
  Lemma \ref{no code} implies that 
 $\overline{\Lambda^*}$ cannot code well on the red zone. Since the red zone takes up the vast majority of the principal $M$-block, 
  $\overline{\Lambda^*}$ cannot code well on the principal $M$-block, yielding a contradiction. In more detail:

Let $\gamma$ be as in Axiom \hyperlink{T4}{T4}.
By Proposition \ref{codes exist} there is an code $\Lambda^*$ such that for almost all $s\in{\bk^c}$, 
\[\dbar(\bar\Lambda^*(s), \phi(s))<10^{-9}{\gamma}.\]
 
By the Ergodic Theorem there is a $N_0$ so large that for a set $E\subseteq {\bk^c}$ of measure $7/8$ for all 
$s\in E$ and all $N>N_0$, $s$ is mature for $N$ and if $B$ is the principal $N$-block of $s$ then 
\begin{equation}
\dbar(\overline{\Lambda^*}(s\rest B), \phi(s)\rest B)<10^{-9}{\gamma}.\label{code accuracy}
\end{equation}
 Let $s\in E$.
 Choose an $N>N_0$ such that the code length of $\Lambda^*$ is much smaller than $q_N$, ${1\over Q^N_1}<10^{-9}$ and $l_N>10^{12}$. Apply 
 Lemma \ref{new red zone}, with $\delta=10^{-9}$ to find an $M$ and $\la R_i:i<i^*\ra$  satisfying the 
 conclusions of Lemma \ref{new red zone}. Since $\bigcup_{i<i^*}R_i\subseteq q_M$ we view $\bigcup_{i<i^*}R_i$ as a subset of the principal $M$-block of $s$.
 
Each $R_i$ is a union of collections of locations of the form $\bigcup_{j<j_0} U_j$,  with each $U_j$ consisting of the locations of $u_j^{s^{n_i}_0}$ for $j\in [0,j_0)$  {(for some $j_0$)}.\footnote{$s_0^{n_i}$  is as in condition \ref{unions of blocks}.c) of Lemma \ref{blocking out}.}
Moreover there is a $t$ such that  each power $u_j^{s^{n_i}_0}$ 
is $\beta$-matched with a $v_{j+t}^{s^{n_i}_0}$ in $\phi(s)$ for some $t\ne 0$.

 Because $j_0>\varepsilon_nk_n$  axiom \hyperlink{T6}{T6} applies and thus for at least 
 $(1-{|G^{n_i}_1|\over Q^{n_i}_1}+\mu_{n_i})$ proportion of  $\{u_0, u_1, \dots u_{j_0-1}\}$, 
 $u_j$ and $v_{j+t}$ are in different $G^{n_i}_1$-orbits. In Lemma \ref{no code}, 
 inequality \ref{tres bizarre} implies that  if $u_i$ and $v_{i+t}$ are in different $G^{n_i}_1$ orbits then, restricted to  the overlaps of the locations of all of the $u_j^{s^{n_i}_0}$ and 
 $v_{j+t}^{s^{n_i}_0}$, the $\dbar$ distance between $\bar\Lambda^*(s)\rest U_j$ and 
 $\phi(s)\rest U_j$ is at least ${1\over 50}(1-{1\over Q^{n_i}_1}){\gamma}$. Since the first and last 
 powers of $u_j$ in  $u_j^{s^{n_i}_0}$'s take up $2/s^{n_i}_0$ of $u_j^{s^{n_i}_0}$ and 
 $s^{n_i}_0\ge l_{n_i}/2-2$, we know that 
 \begin{equation*}
 \dbar(\overline{\Lambda^*}(s)\rest U_j, \phi(s)\rest U_j)\ge (1-10^{-11}){1\over 50}(1-{1\over Q^{n_i}_1}){\gamma}
 \end{equation*}
 Because the proportion of $j$'s for which $u_j$ and $v_{j+t}$ are in different $G^{n_i}_1$-orbits is at least $(1-{|G^{n_i}_1|\over Q^{n_i}_1}+\mu_{n_i})$
it follows that 
 \[\dbar(\bar\Lambda^*(s)\rest\bigcup_{j<j_0}U_j, \phi(s)\rest\bigcup_{j<j_0}U_j )\]
  is at least 
 \begin{equation*}
(1-{|G^{n_i}_1|\over Q^{n_i}_1}+\mu_{n_i})(1-10^{-11}){1\over 500}(1-{1\over Q^{n_i}_1}){\gamma}.
 \end{equation*}
 This in turn is at least ${\gamma}/1000$. Since $R_i$ is a union of sets of the form $\bigcup_{j<j_0}U_j$:
 \begin{equation*}
 \dbar(\overline{\Lambda^*}(s)\rest R_i, \phi(s)\rest R_i)\ge {\gamma}/1000.
 \end{equation*}
 Since $\bigcup_{i<i^*}R_i$ has density at least $1-10^{-9}$ if $B$ is the principal $M$-block of $s$:
 \[\dbar(\overline{\Lambda^*}(s\rest B), \phi(s)\rest B)>{\gamma}/10^4.\]
  However this contradicts the inequality \ref{code accuracy}.\qed
  
  \begin{corollary}\label{exact central}
  Let ${\bk^c}$ be a  circular system built from a circular construction sequence satisfying the timing assumptions. Then $\beta$ is a central value if 
  and only if there is a $\phi\in \overline{\{T^n:n\in\nn\}}$ with $\phi^\pi=\mcs_\beta$. It follows that for each construction sequence $\la k_n, l_n:n\in\nn\ra$ satisfying the Numerical Requirements collected in Section \ref{numreqs},   the central values form a subgroup of the unit circle.
  \end{corollary}
  
\pf Theorem \ref{central values in closure} says that if $\beta$ is central, then there is a $\phi\in   \overline{\{T^n:n\in\nn\}}$ with $\phi^\pi=\mcs_\beta$. Theorem \ref{centralizer central values} is the converse. To see the last statement we prove in Section \ref{the specs} that for every coefficient sequence satisfying the Numerical Requirements, we can find a  circular construction sequence satisfying the timing assumptions. \qed

  \subsubsection{Isomorphisms Between ${\bk^c}$ and ${(\bk^c)}\inv$}
  We now prove a theorem closely related to Theorem \ref{centralizer central values}
  \begin{theorem}\label{conjugacy central values}
 Suppose that $({\bk^c},\mcb, \nu, sh)$ is a  circular system built from a circular construction sequence satisfying the timing assumptions. Suppose that 
 $\phi: ({\bk^c},\mcb, \nu, sh) \to ({(\bk^c)}\inv,\mcb, \nu, sh)$ is an isomorphism. Then 
 $\phi^\pi=\natural \circ\mcs_\beta$ for some central value $\beta$.
  \end{theorem}
  
  \pf  We 
  concentrate here on the differences with  the proof of Theorem \ref{centralizer central values}. The general outline is the same: 
  Fix a $\phi$. Then there is a unique $\beta$ such that $\phi^\pi=\natural \circ \mcs_\beta$.  Suppose that $\beta$ is not central. Choose a stationary code $\overline{\Lambda^*}$ that well approximates $\phi$ in 
  terms of $\dbar$ distance (say within ${\gamma}/10^{10}$), and derive a contradiction by choosing a large $M$ and getting lower bounds 
  for $\dbar$ distance along the principal $M$-block of a generic $s$.  
  
  This is done by first comparing a typical $s$ with $\mcs_\beta(s)$. As in Theorem \ref{centralizer central values}, a definite fraction of a large principal $M$-block of $s$ is misaligned with $\mcs_\beta(s)$. But most of the $n$-blocks of $\mcs_\beta(s)$ are aligned with reversed $n$-blocks of $\natural(\mcs_\beta(s))$ that have been shifted by a very small amount. This can be quantified by looking at the codes $\bar{\Lambda}_n$ for large $n$, which agree with $\natural$ on the $M$-block of $\mcs_\beta(s^\pi)$.

  Here are more details. Recall $\natural$ is the limit of  a particular sequence of stationary codes $\la\bar{\Lambda}_n:n\in\nn\ra$. The proof of  Theorem \ref{mr natural} showed that for almost all $s^\pi\in \mck$ for all large enough $n$ the principal $n$-blocks of $\bar{\Lambda}_n(s^\pi)$ and $\bar{\Lambda}_{n+1}(s^\pi)$ agree. 
Fix a generic $s\in{\bk^c}$ and a large $N$ such that: 
  \begin{enumerate} 
\item the code $\overline{\Lambda^*}$ codes $\phi$ well on the principal $n$-block of $s$ for all $n\ge N$, 
\item for all $n\ge N$  
the principal $n$-blocks of $\bar{\Lambda}_n(\mcs_\beta(\pi(s)))$ and $\bar{\Lambda}_{n+1}(\mcs_\beta(\pi(s)))$ 
agree,
\item $s$ is mature at $N$, 
\item the length of $\Lambda$ is very small relative to $N$ and  
\item $l_{N}$ is very large. 
\end{enumerate}
 Comparing $\pi(s)$ and $\mcs_\beta(\pi(s))$, Lemma
 \ref{new red zone} gives us an $M>N$ and a red zone in the principal $M$-block  $s$. We assume that the red 
 zones take up at least $1-10^{-9}$ proportion of the principal $M$-block and have the form given in Lemma 
 \ref{new red zone}. 

We will derive a contradiction by showing that $\overline{\Lambda^*}$ cannot code well. This is done by considering the blocks of 
$\phi(s)$  that are lined up with 
the red zones of the principal $M$-block of $s$ and using Lemma \ref{no code} to see that $\overline{\Lambda^*}$ cannot code well on these 
sections. This is possible because the mismatched $n$-blocks of $\mcs_\beta(\pi(s))$ are lined up closely with the $n$-blocks of 
$\natural(\mcs_\beta(\pi(s)))=\phi^\pi(s)$. Explicitly:
\smallskip

Use Lemma \ref{new red zone}, to choose red zones $\la R_i:i<i^*\ra$  that take up a $1-10^{-9}$ proportion of the principal $M$-block of $s$.\footnote{We use the notation in  Lemma \ref{new red zone} and
Theorem \ref{mr natural}.} 

 The boundary portions of $n$-words with $n<M+1$ take up at most 
$2/l_{M}$ proportion of the overlap of the principal $M$-blocks of $s$ and $\phi(s)$. Since this proportion is so small, Remark \ref{cheating on dbar} allows us to completely ignore 
blocks corresponding to $n_i$-words in $s$ that are lined up with 
boundary in $\phi(s)$ and vice versa.

We now examine the how $\natural(\mcs_\beta(\pi(s)))$ compares with $\mcs_\beta(\pi(s))$. 
Temporarily denote $\mcs_\beta(\pi(s))$ by $s'$. By the choice of $s$, for all $n\in [N, M]$ the 
alignments of the principal $n$-blocks of $\bar\Lambda_n(s')$ and $\bar\Lambda_M(s')$  agree.

The red zones of $s^\pi$ line up blocks of the form $u_j^{s^{n_i}_0}$ with blocks of the form 
$v_{j+t}^{s^{n_i}_0}$ occurring in $s'$  that are shifted by $d^{n_i}(s)$ (so $t\ne 0$). Except for those blocks that line up with boundary portions of $\natural(s')$
  these blocks are lined up with blocks of the 
form $sh^{A_{n_i}}(\rev{v_{k_{n_i}-(j+t)-1}})$ in $\natural(s')$.\footnote{See the qualititative discussion of $\natural$ that occurs after its definition in \cite{global_structure}.} Inequality \ref{An is small}, says that $A_{n_i}<2q_{n_i-1}$. In particular the blocks of powers of $v_{j+t}$  are lined up with a very small shift of $\rev{v_{k_{n_i}-(j+t)-1}}$ in $\natural(s')$.

Thus vast majority of  blocks $U_j$ that are positions of $u_j^{s^{n_i}_0}$ in $s^\pi$ are lined up with a  shift by less than $q_{n_i}$ of a
block of $\natural(s^\pi)$ in a  position of $v^{s^{n_i}_0}_{k_{n_i}-(j+t)-1}$ in $\natural(s')$ . 
Consider  $s$ and $\phi(s)$. Suppose that $u_j$ are the $n_i$-words of $s$ corresponding to the $U_j$ and $v_{k_{n_i}-(j+t)-1}$ are the $n_i$-words of $\phi(s)$ across from them.
By axiom \hyperlink{T5a}{T5a},  at most ${1\over Q^{n_i}_1}+\mu_{n_i}$ of the  $j<j_0$ happen to have 
$[u_j]_1\in G^{n_i}_1[v_{k_{n_i}-(j+t)-1}]_1$.
At least $1-{1\over Q^{n_i}_1}+\mu_{n_i}$ proportion of the powers of  $u_j$ the $\dbar$-distance between 
$\overline{\Lambda^*}$ and $\phi$ is at 
least ${1\over 500}(1-{1\over Q^n_1}){\gamma}$.

It follows that on $R_i$ the $\dbar$-distance is at least ${\gamma}/1000$. If we choose $\bigcup_{i<i^*} R_i$ to 
have density at least $1-10^{-9}$ and let $B$ be the principal $M$-block of $s$ then (as in Theorem 
\ref{centralizer central values}) 
\[\dbar(\overline{\Lambda^*}(s\rest B), \phi(s)\rest B)>{\gamma}/10^4,\]
a contradiction.\qed

\subsection{Synchronous and Anti-synchronous  Isomorphisms}\label{perfect matches section}
View a circular system $({\bk^c}, \mcb, \nu, sh)$  as an element $T$ of  the space $\MPT$ endowed with the weak topology.
\begin{theorem}\label{perfect match}
Suppose that ${\bk^c}$ is a circular system  satisfying the timing assumptions. Then:
\begin{enumerate}
\item If there is an isomorphism $\phi:{\bk^c}\to {\bk^c}$ such that $\phi\notin \overline{\{T^n:n\in \poZ\}}$, then there 
is an isomorphism $\psi:{\bk^c}\to {\bk^c}$ such that $\psi\notin \overline{\{T^n:n\in \poZ\}}$ and $\psi^\pi$ is the identity map.

\item If there is an isomorphism $\phi:{\bk^c}\to (\bk^c)^{-1}$ then there is an isomorphism $\psi:{\bk^c}\to (\bk^c)^{-1}$ 
such that $\psi^\pi=\natural$.
\end{enumerate}
\end{theorem}

\pf To see the first assertion, let $\phi:{\bk^c}\to {\bk^c}$ be an isomorphism with 
$\phi\notin \overline{\{T^n:n\in \poZ\}}$. Then by Theorem \ref{centralizer central values}, 
$\phi^\pi=\mcs_\beta$ for a central $\beta$.  Corollary \ref{exact central} implies that there is a 
$\theta\in \overline{\{T^n:n\in\nn\}}$ such that $\theta^\pi=\mcs_{-\beta}$. Then $\phi\circ\theta:{\bk^c}\to {\bk^c}$ is 
an isomorphism such that $(\phi\circ\theta)^\pi$ is the identity map. Since $\overline{\{T^n:n\in\nn\}}$ is a 
group, $\phi\circ\theta\notin \overline{\{T^n:n\in\nn\}}$.

The proof of the second assertion is very similar. Suppose that $\phi:{\bk^c}\to(\bk^c)^{-1}$ is an isomorphism. 
Then, by Theorem \ref{conjugacy central values},  $\phi^\pi=\natural\circ\mcs_\beta$ for a central $\beta$. 
Let $\theta\in \overline{\{T^n:n\in\nn\}}$  be such that $\theta^\pi=\mcs_{-{\beta}}$. Then $\phi\circ\theta$ is an isomorphism between ${\bk^c}$ and $(\bk^c)^{-1}$ with $(\phi\circ\theta)^\pi=\natural$.

 \section{The Proof of the Main Theorem}\label{pf of mt}

In this section we prove the main theorem of this paper, Theorem \ref{lele}. By  Fact \ref{comp anal}, it suffices to prove the following:

\begin{theorem}\label{reductio}
There is a continuous function $F^s:\trees\to \difflam$ such that for $\mct\in \trees$, if $T=F^s(\mct)$:
\begin{enumerate}
\item $\mct$ has an infinite branch if and only if $T\cong T^{-1}$,
\item $\mct$ has two distinct infinite branches if and only if 
\[C(T)\ne \overline{\{T^n:n\in\poZ\}}.\]

\end{enumerate}
\end{theorem}
We split the proof of this theorem into three parts. In the first we assume the timing assumptions hold,  define $F^s$ and show that it is a reduction. In the second part we show that 
$F^s$ is continuous.

The third part of the proof augments the specifications of \cite{FRW} with two additional randomness properties, shows that the additional properties imply the timing assumptions and describes how to perform the word construction from \cite{FRW} with these 
additional requirements. We present the third part of the proof separately in Section \ref{the specs}.

We begin by defining $F^s$.  The main result of \cite{FRW} relied on the construction of a  continuous function $F:\trees \to \MPT$ such that for all $\mct\in \trees$, if $S=F(\mct)$ then:

\begin{description}
\item[Fact 1] $\mct$ has an infinite branch if and only if $S\cong S^{-1}$,
\item[Fact 2] $\mct$ has two distinct infinite branches if and only if 
\[C(S)\ne {\{S^n:n\in\poZ\}}.\]
\item[Fact 3] The function $F$ took values in the strongly uniform odometer based transformations
 and for $S$ in the range of $F$, $S\cong S^{-1}$ if and only if there is an anti-synchronous 
 isomorphism $\phi$ between $S$ and $S^{-1}$.
\item[Fact 4]\label{4 I hope} (\cite{FRW}, Corollary 40, page 1565) If $S$ is in the range of $F$ and 
$C(S)\ne\{S^n:n\in \poZ\}$ then there is a synchronous $\phi\in C(S)$ such that for some $n$, 
non-identity element $g\in G^n_1$ and all generic $s\in \bk$ and all large enough $m$, if $u$ and $v$ 
are 
the principal $m$-subwords of $s$ and $\phi(s)$ respectively then:
\[[v]_1=g[u]_1.\]

\item[Fact 5] (Equations 1 and 2 on pages 1546 and 1547 of \cite{FRW}) For all $n_0$ there is an $M$ such that if $\mct$ and $\mct'$ are trees and\footnote{See Section \ref{dst basics} for notation.} 
\[\mct\cap \{\sigma_n:n\le M\}=\mct'\cap \{\sigma_n:n\le M\}\]
then the first $n_0$-steps of the construction sequences for $F(\mct)$ are equal to the first $n_0$-steps of the construction sequence for $F(\mct')$;
i.e. $\la\mcw_k(\mct):k<n_0\ra=\la \mcw_k(\mct'):k<n_0\ra$. 
\item[Fact 6] The construction sequence for $F(\mct)$ satisfies the \emph{specifications} given in \cite{FRW}.  In Section \ref{old specs}, these specifications are augmented by the addition of J10.1 and J11.1. In Section \ref{the TA} we argue that if $\la \mcw_n:n\in\nn\ra$ is a construction sequence for an odometer based system that satisfies the augmented specifications, then the associated circular construction sequence $\la \mcw_n^c:n\in\nn\ra$ satisfies the timing assumptions.

Moreover the construction sequence for $F(\mct)$ is strongly uniform and hence the construction sequence for $\mcf\circ F(\mct)$ is strongly uniform. 
\item[Fact 7] Construction sequences satisfying the augmented specifications are easily built using the techniques of \cite{FRW} with no essential changes; consequently we can assume that the the construction sequences for $F(\mct)$ satisfy the augmented specifications.
\end{description}

In \cite{prequel} (Theorem 60) it is shown  that if $\la W^c_n:n\in\nn\ra$ is a strongly uniform circular 
construction sequence with coefficients $\la k_n, l_n:n\in\nn\ra$,
where $\la l_n: n\in\nn\ra$ grows fast enough  and $|\mcw^c_n|$ goes to infinity then there is a smooth 
measure preserving diffeomorphism $T\in \difflam$ measure theoretically isomorphic to $\bk^c$. This 
gives a map $R$ from circular systems with fast growing coefficients to $\difflam$.
 
 If $\mcf$ is the canonical functor from odometer systems to circular systems we define
\[F^s=R\circ\mcf\circ F\].

\begin{figure}[h]
\centering
\includegraphics[height=.25\textheight]{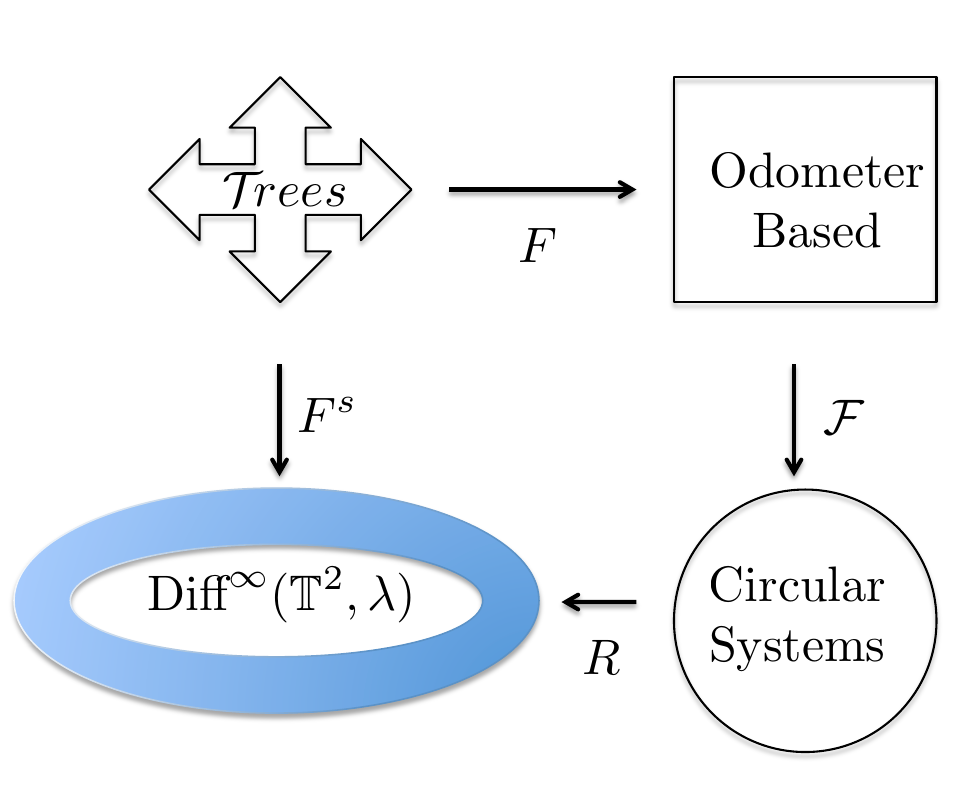}
\caption{The definition of  $F^s$.}
\label{def of F^s}
\end{figure}

\subsection{$F^s$ is a Reduction}
Because  $R$ preserves isomorphism, to show that $F^s=R\circ \mcf\circ F$ is a reduction, it is suffices to 
show that $\mcf\circ F$ is a reduction. Let $S$ be the transformation corresponding to the system $\bk=F(\mct)$ and $T$ the transformation corresponding 
to $\bk^c=\mcf\circ F(\mct)$.
\medskip

\bfni{Item 1 of Theorem \ref{reductio}:}  Suppose that $\mct$ is a tree and $\mct$ has an infinite branch. By Facts 1 and 3,   there is an anti-synchronous isomorphism $\phi:\bk\to \bk^{-1}$. By Theorem 105 of \cite{global_structure}, if $\bk^c=\mcf(\bk)$, there is an isomorphism $\phi^c:\bk^c\to (\bk^c)^{-1}$. 

Now suppose that $F^s(\mct)\cong (F^s(\mct))^{-1}$, then $\bk^c\cong (\bk^c)^{-1}$. By Fact 6,  the construction sequence $\la \mcw_n^c:n\in\nn\ra$   for $\mcf^s(\mct)$ satisfies the timing assumptions.  By Theorem \ref{perfect match}, there is an anti-synchronous isomorphism $\phi^c:\bk^c\to (\bk^c)^{-1}$. Again by Theorem 105 of \cite{global_structure}, there is an isomorphism between $\bk$ and $\bk^{-1}$. By \cite{FRW}, $\mct$ has an infinite branch.
\medskip

\bfni{Item 2 of Theorem \ref{reductio}:} Suppose that $\mct$ has at least two infinite branches.  Then the centralizer of $S=F(\mct)$ is not equal to the powers of $S$. By Fact 4, we can find a synchronous
$\phi\in C(S)\setminus \{S^n:n\in \poZ\}$. Let $\psi=\mcf(\phi)$, then $\psi$ is synchronous. We claim that $\psi\notin \overline{\{T^n:n\in\poz\}}$.
Using Fact 4, and lifting the group action of $G^n_1$  and the  equivalence relation $\mcq^n_1$, we see that for all generic $s^c \in \bk^c$, and all large enough $m$, if $u^c$ and $v^c$ are the principal $m$-subwords of $s^c$ and $\psi(s^c)$, respectively, then:
\[ [v^c]_1=g[u^c]_1\]
for some $g\ne e$. In particular, $[v^c]_1\ne[u^c]_1$.

{By the timing assumption {T4}, 
{there is a $\gamma>0$ such that for all large $m$ and all shifts $A$ with $|A|$ of size less than $q_m/2$, we have 
	\begin{equation}
	\dbar(T^A(u^c),v^c)>\gamma.\label{getting bigger than gamma}
	\end{equation} 
Suppose that $\psi\in \overline{\{T^n:n\in\poZ\}}$. Then, by Proposition \ref{weak and dbar}, we can find an $A\in \poZ$ and a generic $s^c$ such that 
	\begin{equation}
	\dbar(T^A(s^c),\psi(s^c))<\gamma/2.\label{smaller than gamma}
	\end{equation}
But inequality \ref{smaller than gamma} and the Ergodic Theorem imply that for large enough $m\gg A$ if $u^c$ and $v^c$ are the principal $m$-blocks of 
$s^c$ and $\psi(s^c)$ then 	
	\begin{equation*}
	\dbar(T^A(u^c),v^c)<\gamma,
	\end{equation*}
contradicting inequality \ref{getting bigger than gamma}.}

Now suppose that there is a $\psi\in C(T)$ such that $\psi\notin \overline{\{T^n:n\in\poz\}}$. Then by Theorem \ref{perfect match}, there is such a 
$\psi$ that is synchronous. In particular, for all $n$, $\psi\ne T^n$. Thus if $S$ is the transformation corresponding to $F(\mct)$, $\mcf^{-1}(\psi)$ belongs to the centralizer of $S$ and is not a power of $S$.

\subsection{$F^s$ is Continuous}

Fix a metric $d$ on $\difflam$ yielding the $C^\infty$-topology. For each circular system $T$, let $\la P^T_n:n\in\nn\ra$ be the sequence of collections of prewords used to construct $T$.
By Proposition 61  of \cite{prequel},  
 given $T=F^s(\mct)$ and a $C^\infty$-neighborhood $B$ of $T$,
  there is a large enough $M$, for all   $S\in \mbox{range}(R)$ if $\la P_n^S:n\le M\ra=\la P_n^T:n\le M\ra$, 
  then $S\in B$. For all odometer based transformations, the sequence $\la \mcw_n:n\le M\ra$ determines 
  $\la P_n:n\le M\ra$. {Hence} for all $\mct'$, if the first $M$ members of the construction sequence for  
  $F(\mct')$ are the same as the first $M$ members of the construction sequence for $F(\mct)$, then 
  $F(\mct')\in U$. By Fact 5,  there is a basic open interval $V\subseteq \trees$ that contains $\mct$ 
  and is such that the first $M$ members of the construction sequence  are the same for all 
$\mct'\in V$. It follows that for all $\mct'\in V, F^s(\mct')\in U$.

\subsection{Numerical Requirements Arising from Smooth Realizations}
The construction of $R$ depends on various estimates that put lower bounds on the growth of the coefficient sequences.    We now list these numerical requirements. The claims in this subsection presuppose a knowledge of \cite{prequel}.

{The map $R$ depends on various smoothed versions $h_n^s$ of the permutations $h_n$ of the unit interval arising from $\la \mcw_n:n\in\nn\ra$. To solve this problem, we fix in advance such 
 approximations, making sure that each approximation $h_n^s$ agrees sufficiently well with $h_n$ as to not 
 disturb the other estimates.}

This introduces various numerical constraints on the growth of the $l_n$'s. The diffeomorphism $T$ is built as a limit of periodic approximations $T_n$. To make the sequence of $T_n$'s converge at each stage, $l_n$ must be chosen sufficiently large. Thus the growth rate of $l_n$  depends on $\la k_m, s_m, h_m:m\le n\ra, \la l_m:m<n\ra$, $s_{n+1}, h_{n+1}$. Since there are only finitely many possibilities for  $\la h_m:m\le n\ra$'s corresponding to a given sequence $\la k_m:m\le n\ra$, $\la s_m:m\le n+1\ra$ we can find one growth rate that is sufficiently fast to work for all choices of $h_m$'s. {This is discussed in detail on page 34 of \cite{prequel}, where the lower bound is called $l_n^*$.}

 \begin{numreq}\label{inherent in smooth} 
 $l_n$ is big  enough relative to a lower   bound determined by $\la k_m, s_m:m\le n\ra$, $\la l_m:m<n\ra$  and $s_{n+1}$ to make the periodic approximations to the diffeomorphism converge. Moreover $k_n\le l_n$.
 \end{numreq}
 
 \begin{remark}
 {Choosing $\alpha_{n+1}$ close to $\alpha_n$ is a fundamental idea of the method of Approximation by Conjugacy, due to Anosov and Katok. By equations \ref{qn} and \ref{pn}, this is equivalent to taking $l_n$ large.  The magnitude of $l_n$ is not calculated, but instead it shown that as $l_n$ increases a sequence of periodic diffeomorphisms well approximates a given periodic diffeomorphism. Then in the original sources \cite{AK} and \cite{katoksbook}, one simply takes $l_n$ \emph{sufficiently large}. This is what requirement \ref{inherent in smooth} is repeating.  }
 \end{remark}
 
 The argument for the ergodicity of the diffeomorphism formally required that:
 \begin{numreq}\label{sn grows}
 $s_n$ goes to infinity as $n$ goes to infinity, $s_{n+1}$ is a multiple of $s_n$. 
 \end{numreq}
 \noindent The reader is referred to example \ref{to confuse the reader} for a discussion of $s(n)$ and its growth.

\medskip
The next requirement makes it possible to choose $s_{n+1}$ and then, by making $k_n$ sufficiently large, construct $s_{n+1}$ sufficiently random words using elements of $\mcw_n$.

\begin{numreq}\label{sn to kn}
$s_{n+1}\le s_n^{k_n}$
\end{numreq}

\section{The Specifications}\label{the specs}

 In this section we describe how the timing assumptions are related to the specifications given in \cite{FRW}, show that they are compatible and indicate how to construct odometer words so that both sets of assumptions hold. This completes the proof of Theorem \ref{reductio}, subject to the verification that all of the Numerical Requirements we have introduced are consistent with the numerical requirements of \cite{FRW}. We take this up in section \ref{numreqs}.
We will assume that the reader is familiar with sections 7 and 8 of \cite{FRW}.
\begin{figure}[h]
\centering
\includegraphics[height=.2\textheight]{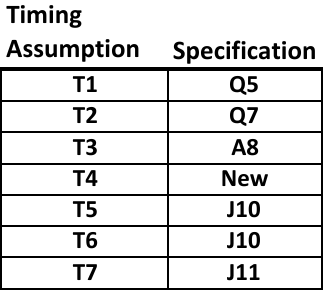}
\caption{The Specifications {in \cite{FRW}} related to the Timing Assumptions {in this paper.}}
\label{correspondence}
\end{figure}

\subsection{Corresponding Specifications}\label{transition info}

Figure \ref{correspondence} gives a table that links the Timing Assumptions we use in this paper to the corresponding 
Specification in \cite{FRW}. (We remind the reader that Appendix \ref{appendixA} has a table giving corresponding notation 
betwen \cite{FRW} and this paper.)

Specification T4 doesn't directly correspond to one of the Specifications, but (as we will show) holds naturally in the circular words lifted from an odometer construction satisfying the specifications.

\medskip

\begin{numreq}\label{two epsilons}
In the current construction we have two summable  sequences: $\la \epsilon_n:n\in\nn\ra$ and $\la \varepsilon_n:n\in \nn\ra$. We use the lunate ``$\epsilon_n$" notation for the specifications from \cite{FRW} and the classical ``$\varepsilon_n$" notation (``varepsilon" in LaTeX) for the numerical requirements relating to circular systems and their realizations as  diffeomorphisms. A requirement for the construction is that 
\[\epsilon_n<\varepsilon_n.\]
We also assume that the $\epsilon_n$'s are decreasing and $\epsilon_0<1/40$.
\end{numreq}

\subsection{Augmenting the  Specifications from \cite{FRW}}\label{old specs}

The paper \cite{FRW} constructs a reduction $F$ from the space of trees to the odometer based systems. The system 
$\bk=F(\mct)$ was built according to a list of specifications which we reproduce here in order to show how to strengthen them to imply the timing assumptions used in the proofs of Theorems \ref{perfect match} and \ref{reductio} and to verify that the strengthened assumptions are consistent. The specifications directly relevant to the timing assumptions are \emph{J10} and \emph{J11}. The others, which describe the scaffolding for the construction, are only relevant in that they set the stage for the application of the functor $\mcf$ defined in section  \ref{odo circ review}.

Here are some definitions from \cite{FRW} that are used in the specifications. {We advise the reader that a table giving the notational changes between \cite{FRW} and this paper is in Appendix \ref{appendixA}.}

Fix an 
enumeration of the finite sequences of natural numbers, $\la \sigma_n:n\in\nn\ra$, with the property that if $\sigma$ is an initial segment of $\tau$ then $\sigma$ is enumerated before $\tau$. Let $\mct$ be a tree whose elements are $\la \sigma_{n_i}:i\in \nn\ra$. Here are the specifications for the construction sequence 
 $\mcw=\mcw(\mct)$ used to build $\mcf(\mct)$. 

There is a sequence of groups $G^n_s$  built as follows.
For all $n$, $G^n_0$ is the trivial group $(e)$ and if we let 
\[X^n_s=\{\sigma_{n_i}:i\le n \mbox{ and }\sigma_{n_i} \mbox{ has length }s\}\]
 then 
\[G^n_s=\sum_{\sigma\in X^n_s}(\poZ_2)_\sigma\]
i.e. $G^n_s$ is a direct sum of copies of $\poZ_2$ indexed by elements of $X^n_s$. There are canonical homomorphisms from $G^n_{s+1}$ to $G^n_{s}$  that send a generator of $G^n_{s+1}$ corresponding to a sequence of the form $\tau\cat j$ to the generator of $G^n_s$ corresponding to $\tau$.

The sequence $\la \mcw_n:n\in\nn\ra$, equivalence relations $\mcq^n_s$ and the group actions of $G^n_s$  are constructed inductively. The words in 
$\mcw_n$ are sequences of elements of $\Sigma=\{0,1\}$.  To start
$\mcw_0=\{0,1\}$ and $\mcq^0_0$ is the trivial equivalence relation
with one class. The collection of words $\mcw_n$ is built when the $n^{th}$ element of $\mct$ is considered. We will say that words in $\mcw_n$ have \emph{even} parity and words in $\rev{\mcw_n}$ have \emph{odd} parity.

{We begin by restating the specifications from \cite{FRW} using the indexing conventions in this paper 
($n\mapsto n+1$ vs $m\mapsto n$). E1-A9 are exactly the same, however we modify the joining specifications J10, 
J11 slightly for the needs of this paper.}

\begin{enumerate}
\item[E1.] Any pair $w_1, w_2$ of words in $\mcw_n$ have the same
length.

\item[E2.]

Every word in
$\mcw_{n+1}$ is built by concatenating words in
$\mcw_n$. Every word in $\mcw_n$
occurs in each word of $\mcw_{n+1}$ exactly $p_n^2$ times, where $p_n$ is a large prime number chosen when the $n^{th}$ element of $\mct$ is considered.

{\item[E3.](Unique Readability)} If 
$w\in\mcw_{n+1}$ and 
\[w=pw_1\dots w_ks\]
where each $w_i\in \mcw_n$ and  $p$ or $s$ are sequences of $0$'s and $1$'s that have length less than that of any word in $\mcw_n$,
then both $p$ and $e$ are the empty word.
If $w,w'\in \mcw_{n+1}$ and $w=w_1w_2\dots w_{k_n}$ and
$w'=w'_1w'_2\dots w'_{k_n}$ with $w_i, w_i'\in\mcw_n$, and $k=[k_n/2]+1$ we have
$w_{k}w_{k+1}\dots w_{k_n}\ne w'_1w'_2\dots w'_{k_n-[k]-1}$; i.e. the first half of $w'$ is
not equal to the second half of $w$.
\end{enumerate}
\noindent Let $s(n)$ be the length of the longest sequence among the first $n$ sequences in $\mct$ and if $\mct=\la \sigma_{n_i}:i\in\nn\ra$ then $M(s)$ is the least $i$ such that $\sigma_{n_i}$ has length $s$.

The equivalence relations $\mcq^n_s$ on $\mcw_n$ are defined for all $s\le s(n)$. The equivalence relation $\mcq^0_0$ on $\mcw_0$ is the trivial equivalence relation with one class.

\begin{enumerate}

\item[Q4.] Suppose that $n=M(s)$. Then any two words in the same
$\mcq^n_s$ equivalence class agree on an initial segment of proportion least $(1-\epsilon_n)$.

\item[Q5.] For $n\ge M(s)+1$, $\mcq^{n}_s$ is the product equivalence
relation of
$\mcq^{M(s)}_s$. Hence we can view
$\mcw_n/\mcq^n_s$ as sequences of elements of $\mcw_{M(s)}/\mcq^{M(s)}_s$ and similarly for
$\rev{\mcw_n}/\mcq^n_s$.

\item[Q6.] ${\mcq}^n_{s+1}$ refines ${\mcq}^n_s$ and each ${\mcq}^n_s$
class contains
$2^{e(n)}$ many ${\mcq}^n_{s+1}$ classes, where $e$ is a strictly increasing function. {The speed of growth of $e$ is discussed in section \ref{numreqs}.}

\item[A7.] $G^n_s$ acts freely on $\mcw_n/{\mcq}^n_s\cup
\rev{\mcw_n/{\mcq}^n_s}$ and the $G^n_s$ action is subordinate to the
$G^n_{s-1}$ action via the natural homomorphism $\rho_{s,s-1}$ from
$G^n_s$ to $G^n_{s-1}$.

\item[A8.] The canonical generators of $G^{M(s)}_s$ send elements of
$\mcw_{M(s)}/\mcq^{M(s)}_s$ to elements of
$\rev{\mcw_{M(s)}}/\mcq^{M(s)}_s$ and vice versa. 

\item[A9.] If $M(s)\le n$ and we view
$G^{n+1}_s=G^{n}_s\oplus H$ then the action of $G^n_s$ on $\mcw_n/\mcq^n_s\cup
\rev{\mcw_n/\mcq^n_s}$ is extended to 
an action on 
$\mcw_{n+1}/\mcq^{n+1}_s\cup\rev{\mcw_{n+1}/\mcq^{n+1}_s}$ by the skew diagonal action. If $H$ is
non-trivial then $H=\poZ_2$ and its canonical generator maps $\mcw_{n+1}/\mcq^{n+1}_s$ to
$\rev{\mcw_{n+1}/\mcq^{n+1}_s}$.
\end{enumerate}
\paragraph{Note:} {While it is not explicitly given as a specification in \cite{FRW}, the construction sequence has the property that if $g\in G^n_s$ is a canonical generator, then for 
$m>n, \mcw_m/\mcq^m_s$ is closed under the skew diagonal action of $g$.}
\medskip

Suppose that $u$ and $v$ are elements of $\mcw_{n+1}\cup \rev{\mcw_{n+1}}$ and $(u', v')$ an ordered  pair from $\mcw_{n}\cup \rev{\mcw_{n}}$. Suppose that $u$ and $v$ are in positions shifted relative to each other by  $t$ units. Then an \emph{occurrence} of $(u', v')$ in $(sh^t(u), v)$ is a $t'$ such that $u'$ occurs in $u$ starting at $t+t'$ and in $v$ starting at $t'$. 
Let  $Q^n_s$ be the number of classes of $\mcq^n_s$ and $C^n_s$ be the number of elements of each $Q^n_s$ 
class.\footnote{We have changed the variables used in the statement of J10 in \cite{FRW} to conform to the notation 
described in the appendix \ref{appendixA}.}

\smallskip
To prove the timing assumptions we need to strengthen specifications J10 and J11 to deal with $\dbar$-distance on initial and tail segments and on words that are shifted. The spirit of specification J10 is that pairs of $n$-words $(u',v')$ occur randomly in the overlap of  $u$ and $v$ when $u$
is shifted by a suitable multiple $t$ of the lengths of $n$-words. J10.1 says the same thing relative to non-trivial initial segments of the overlap of the shift of $u$  and $v$. 

The specification $J11$ says that if $[u]_s$ is in the $G^n_s$-orbit of $[v]_s$ and $s$ is maximal with this property, then the occurrences of 
$(u',v')$ are approximately conditionally random. More explicitly, suppose that $g[u]_s=[v]_s$, and we are given $u'\in \mcw_n$. Then there are   $Q^n_s$ many pairs of 
$\mcq^n_s$-classes $([u^*]_s,[v^*]_s)$ with $g[u^*]_s=[v^*]_s$, and so $([u']_s,[v']_s)$ should occur 
randomly $1/Q^n_s$ proportion of the time. There are $C^n_s$ many elements of $\mcw_n$ in the 
$\mcq^n_s$-classes, and conditional on $g[u']_n=[v']_n$, the chances of such a pair $(u',v')$ randomly 
matching is $1/(C^n_s)^2$.  The specification J11.1 strengthens this (but only for $\mcq^n_0$, which 
is the trivial equivalence relation and $G^n_0=\la e\ra$) by asking that this holds over any non-trivial 
interval of length $j_0K_n$ at the beginning or end of an $n+1$-word. 
\smallskip

Here are the joining specifications as given in \cite{FRW}:

\begin{enumerate}

\item[J10.] Let $u$ and $v$ be 
elements of
$\mcw_{n+1}\cup\rev{\mcw_{n+1}}$.  Let $1\le t<(1-\epsilon_n)(k_n)$ be an integer. Then for each pair $u',
v'\in\mcw_{n}\cup\rev{\mcw_{n}}$ such that $u'$ has the same parity as $u$ and $v'$ has the same parity as
$v$, let
$r(u',v')$ be the number of occurrences of $(u',v')$ in $(sh^{tK_n}(u),v)$ on
their overlap. Then
\begin{equation*}
\left|{r(u',v')\over k_n-t}-{1\over s_{n}^2}\right|<\epsilon_n.
\end{equation*}

\hypertarget{J11}{\item[J11.] }  Suppose that $ u\in\mcw_{n+1}$ and $v\in \mcw_{n+1}\cup
\rev{\mcw_{n+1}}$. We let
$s=s(u,v)$ be the maximal $i$ such that there is a $g\in G^{n}_i$ such that
$g[u]_i=[v]_i$. Let $g=g(u,v)$ be the unique $g$ with this property and
$(u', v')\in
\mcw_{n}\times
(\mcw_{n}\cup\rev{\mcw_n})$ be such that
$g[u']_s=[v']_s$. Let $r(u',v')$ be the  number of occurrences of $(u',v')$ in $(u,v)$. Then:
\begin{equation*}
\left|{r(u',v')\over k_n}-{1\over Q^{n}_s}\left({1\over
C^{n}_s}\right)^2\right|<\epsilon_n.
\end{equation*}
\end{enumerate}

\noindent The strengthening of J10 is:

\begin{description}
\item[J10.1] \hypertarget{J10.1}{}
Let $u$ and $v$ be %
elements of
$\mcw_{n+1}\cup\rev{\mcw_{n+1}}$.  Let $1\le t<(1-\epsilon_n)(k_n)$. Let $j_0$ be a number between {$\epsilon_nk_n$}  and $k_n-t$. Then for each pair $u',
v'\in\mcw_{n}\cup\rev{\mcw_{n}}$ such that $u'$ has the same parity as $u$ and $v'$ has the same parity 
as
$v$, let
$r(u',v')$ be the number of $j<j_0$ such that   $(u',v')$ occurs in  $(sh^{tK_n}(u),v)$ in
the $(jK_n)^{th}$ position in their overlap. Then
\begin{equation*}
\left|{r(u',v')\over j_0}-{1\over
s_{n}^2}\right|<\epsilon_n.
\end{equation*}
\end{description}
The next assumption is a strengthening of a special case of J11.
\begin{description}
\item[J11.1] \hypertarget{J11.1}{}Suppose that $ u\in\mcw_{n+1}$ and $v\in \mcw_{n+1}\cup
\rev{\mcw_{n+1}}$ and $[u]_1\notin G_1^n[v]_1$.\footnote{In the language of J11: $s(u,v)=0$, $Q^n_0=1$ and $C^n_0=s_n$.}  Let $j_0$ be a number between $\epsilon_nk_n$ and $k_n$. Suppose that $I$ is either an initial or a tail segment of the interval $\{0, 1, \dots K_{n+1}-1\}$ having length $j_0K_n$. Then for each pair $u',
v'\in\mcw_{n}\cup\rev{\mcw_{n}}$ such that $u'$ has the same parity as $u$ and $v'$ has the same parity as
$v$,  let
$r(u',v')$ be the number of occurrences of $(u',v')$ in $(u\rest I,v\rest I)$. Then:
\begin{equation*}
\left|{r(u',v')\over j_0}-{1\over
s_n^2}\right|<\epsilon_n.
\end{equation*}
\end{description}

{We have augmented the specifications in \cite{FRW} with J10.1 and J11.1. Formally we must argue that it is 
possible to build  construction sequences satisfying the additional specifications. This means 
constructing $s_{n+1}$ many pseudo-random words. This is done using a routine modification of
the techniques of \cite{FRW}, where the collections of words $\mcw_n$ are built probabilistically.
 For $n\ge 1$ the words in  $\mcw_{n+1}$ are built by iteratively substituting  words into $K_{n+1}/K_{M(i)}$-sequences of classes 
$\mcq_i^n$, by  induction on $i\le i^*$ where $i^*$ is maximal with $M(i^*)\le n$. The classes of words $\mcw_{n+1}/\mcq_{n+1}^i$ are built by induction on $i$. A word $w\in \mcw_{n+1}/\mcq^{n+1}_{i+1}$ 
(or in $\mcw_{n+1}$ if $i=i^*$) can be viewed as a result  of substituting elements of $\mcw_n/\mcq^{n}_{i+1}$ (or $\mcw_n$) into a word in $\mcw_{n+1}/\mcq^{n+1}_i$.  
}

{Suppose that $[w]_{i}\in \mcw_{n+1}/\mcq^{n+1}_i$ has been built and is given by $K_{n+1}/K_{M(i)}$ many consecutive classes $C_1C_2\dots C_{K_{n+1}/K_{M(i)}}$. Then $[w]_{i+1}\in \prod_{j<K_{n+1}/K_{M(i)}}C_j$. Viewing these as independent trials and taking $k_n$ large enough (so that $K_{n+1}/K_{M(i)}$ is very large) the finitary \emph{Law of Large Numbers} shows that the vast majority of  choices of $2^{e(n)}$ words satisfy J10, J10.1, J11 and J11.1} 

\begin{remark}\label{befuddle the reader}
As noted in Example \ref{to confuse the reader}, given the number of substitutions to be made (which is one more than  the maximal $s$ such that $\mcq_s^n$ is defined) and the size of the groups $G^n_s$ one can give an explicit formula relating the sizes of $e(n+1)$ and $s_{n+1}$.  Given one of the two one can solve for the other. Moreover when one goes up the other does as well.  This co-determination means that the requirements can be stated in terms of either variable.  We state the requirements in terms of the $s_n$'s.   
\end{remark}

{In the construction, getting the additional .1 for J10 and J11  only involves taking $k_n$ larger than was necessary in \cite{FRW}. This is described in this notation in \cite{GD}.}

This leads to a numerical requirement:
\begin{numreq}\label{lln}
$k_n$ is chosen sufficiently large relative to a lower bound determined by $s_{n+1}$ for the Law of Large Numbers arguments to work.
\end{numreq}

\subsection{Verifying the Timing Assumptions}\label{the TA}

In this section we prove that the augmented specifications E1-J11.1 imply the {timing assumptions}, introduced in Section \ref{intro timing}. 
The first three timing assumptions \emph{T1-T3} follow easily from the results in section \ref{propagano} together with specifications Q5, Q7 and A8.

The following remark is easy and illustrates the idea behind the demonstrations of T4-T7.

\begin{remark}\label{warmed over}
Suppose that $\mcl$ is an alphabet with $s$ symbols in it and $\mcc\subset\mcl$ with $|\mcc|=C$. For $u, v$ words in $\mcl$ of the same length and 
$x,y\in \mcl$, set $r(x,y)$ to be the number of occurrences of $(x,y)$ in $(u,v)$, $r(x, \mcc)$ to be the number of occurrences of  some element of 
$\mcc$ 
opposite an occurrence of  $x$ in $u$ and $f(x)$ to be the number of occurrences of $x$ in $u$. Then for all $\mu>0$ there is a $\epsilon=\epsilon(\mu, s)$ such that whenever $u,v$ are two words in $\mcl$ of the same length $\ell$, if for all $x,y\in \mcl$,
    \begin{eqnarray*}
    \left|{r(x,y)\over \ell}-{1\over s^2}\right|<\epsilon\\
    \end{eqnarray*}

{then} for all $x$:
\begin{equation*}
\left|{r(x,\mcc)\over f(x)}-{C\over s}\right|<\mu
\end{equation*}
\end{remark}

\pf Because $f(x)=\sum_y r(x,y)$, by taking $\epsilon$ sufficiently small we can arrange that 
\[{f(x)\over \ell}\approx{1\over s},\]
and the approximation improves as $\epsilon$ gets smaller.
 Simplemindedly:
\begin{eqnarray*}
{r(x,y)\over f(x)}&=&{r(x,y)\over \ell}{\ell\over f(x)}\\
&\approx&{1\over s^2}s\\
&\approx&{1\over s}
\end{eqnarray*}

Since $r(x, \mcc)=\sum_{y\in \mcc}r(x,y)$ we see that 
\[{r(x, \mcc)\over f(x)}\approx{C\over s}.\]

As we take $\epsilon $ smaller the final approximation improves.\qed

We now establish the timing assumptions \hyperlink{T4}{T4}-\hyperlink{T7}{T7}. {Recall that in the context of the timing assumptions the notation $a\approx b$ means that $|a-b|<\mu_n$.}
\medskip

\bfni{Assumption T5:} 
Assume that specification J10 holds for sufficiently small $\epsilon_n$. 
To use remark \ref{warmed over} to see  \hyperlink{T5}{T5}, take 
$\mcl=\mcw_n$,  the number $f(x)$ to be $|J(v)|$ and $C$ to be the cardinality of any equivalence class of 
$\mcq^n_1$ and $s=s_n$. Since each class of $\mcq^n_1$ has the same number of elements, ${s\over C}$ 
is equal to the number of classes:   ${s\over C}=Q^n_1$. Thus ${C\over s}={1\over Q^n_1}$ and T5 follows. 
\medskip

\bfni{Assumption T6:}
 We can write the set $S$ as:
\[S=\bigcup_{v\in \mcw_n^c}\bigcup_{g\in G^n_1}\{k<j_0:v=w_k\mbox{ and }w'_{k+t}\in g[v]_1\}.\]
which can be written in turn as:
\[S=\bigcup_{v\in \mcw_n^c}\bigcup_{g\in G^n_1}\bigcup_{v'\in g[v]_1}\{k<j_0:v=w_k\mbox{ and }w'_{k+t}=v'\}.\]

Thus, using J10.1, we can estimate the size of $S$ as
\[|S|\approx s_n|G^n_1|C^n_1\left({j_0\over s_n^2}\right).\]
Since $C^n_1=s_n/Q^n_1$ we can simplify this to ${|G^n_1|\over Q^n_1}j_0$. The assumption T6 follows.
\medskip

\bfni{Assumption T7:} Under the assumption that 
$[w_1']_1\notin G^n_1[w_0]_1$, $s=0$ and $\mcq^n_0$ is the trivial equivalence relation.  The estimate in \hyperlink{J11}{J11} simplifies to:
\begin{equation}\label{getting the hype}
\left|{r(u',v')\over k_n}-{1\over s_n^2}\right|<\epsilon_n.
\end{equation}
To apply Remark \ref{warmed over}, we again set $\mcl=W_n$ and $x= v$ and $|J(v)|=f(x)$, in the language of the remark. With this notation, $l=k_n$ and  equation 
\ref{getting the hype} is the hypothesis of Remark \ref{warmed over}. 
The conclusion of the remark 
is that
\begin{equation}
\label{T7 arg}
{|\{t\in J(v):\mcc\mbox{ occurs at $t$ in } [u'_1]_1[u'_2]_1\dots [u'_{k_n-1}]_1\} |\over |J(v)|}\approx{C^1_n\over s_n}.
\end{equation}
Since ${C^1_n\over s_n}={1\over Q^1_n}$, assumption T7 follows.
\medskip

{We note that the verification of T5-T7 uses remark \ref{warmed over} for a small enough $\epsilon(\mu_n,s_n)$.  We 
make this a requirement on $\epsilon_n$.}

\begin{numreq}\label{eps and mu}
 \[ \epsilon_n\mbox{ is sufficiently small  relative to }\mu_n 
 \mbox{ that the timing assumptions T5-T7 hold}.\] 
\end{numreq}

\bfni{Assumption T4:} \hyperlink{T4}{T4} is the hardest timing assumption to verify. We motivate the proof by remarking that if $u, v$ are long mutually random words in a language $\mcl$ that has $s$ letters, then $\dbar(u,v)\approx 1- 1/s^2$. Thus $u$ and $v$ are far apart. Specifications \hyperlink{J10.1}{J10.1} and \hyperlink{J11.1}{J11.1} imply that most $(u,v)$ and their relative shifts are nearly mutually random. We use this to establish that $w_0$ and $w_1$ are distant  in $\dbar$. 

\begin{numreq}\label{yeah yeah yeah}
$\epsilon_0k_0>20$, the $\epsilon_nk_n$'s are increasing and $\sum 1/\epsilon_nk_n$ is finite.
\end{numreq}

Let
\begin{equation*}
\gamma_1=(1-1/4-\epsilon_0)(1-1/\epsilon_0k_0)(1-1/l_0).
\end{equation*}
For $n\ge 2$, set:
\begin{equation*}
\gamma_n=\gamma_1\prod_{0<m<n} (1-10(1/k_m\epsilon_m+1/q_m+1/l_m+1/Q^m_1+\epsilon_{m-1}))
\end{equation*}
and finally:
\begin{equation*}
\gamma=\gamma_1\prod_{0<m} (1-10(1/k_m\epsilon_m+1/q_m+1/l_m+1/Q^m_1+\epsilon_{m-1})).
\end{equation*} 

Assumption \hyperlink{T4}{T4} says that if $w_0, w_1\in \mcw^c_n\cup\rev{\mcw_n^c}$ are not $\mcq^n_1$-equivalent, then the overlaps of sufficiently long initial segments, or sufficiently long tail segments or of a sufficiently long initial segment with a tail segment of $w_0$ and $w_1$ are at least $\gamma$ distant in $\dbar$. In T4  \emph{sufficiently long}  means at least half of the length of the word.  We prove something  stronger by induction on $n$:

\hypertarget{prop 156}{\begin{prop}\label{verify T4}} %
Let $n\ge 0$ and $w_0, w_1\in \mcw^c_{n+1}\cup\rev{\mcw_{n+1}^c}$ {with $[w_0]_1\ne[w_1]_1$}. Let $I$ be an initial segment and $T$ be a tail segment of 
of $\{0, 1, \dots q_{n+1}-1\}$ of the same length $\ell>\epsilon_{n}q_{n+1}$. Then we have:
    \begin{eqnarray}
    \dbar(w_0\rest I, w_1\rest I)&\ge&\gamma_{n+1} \label{14}\\
    \dbar(w_0\rest T, w_1\rest T)&\ge&\gamma_{n+1}  \label{6}\\ 
    \dbar(w_0\rest I, w_1\rest T)&\ge&\gamma_{n+1}. \label{7}
    \end{eqnarray}
\end{prop}

\pf We will consider the situation where $w_0,w_1\in \mcw_{n+1}^c$. The situation where they both belong to $\rev{\mcw_{n+1}^c}$ follows, and the argument in the case where $w_0, w_1$ have different parities is a small variation of the basic argument.

The strategy for the proof is to consider  $n+1$-words $w_0$ and $w_1$ and gradually eliminate small portions of $I$ and $T$ so that we are left with 
only segments of $n$-words that align in $w_0$ and $w_1$ in such a way that they have large $\dbar$-distance. The remaining portions of the $w_0$ and $w_1$ are far apart and they constitute most of the segments of each word. By Remark \ref{cheating on dbar}, we get an estimate on the distance of $w_0$ and $w_1$.

Suppose that 
\begin{eqnarray*}
w_0&=&\mcc(u_0,u_1,\dots, u_{k_n-1})\\
w_1&=&\mcc(v_0,v_1,\dots, v_{k_n-1}),
\end{eqnarray*}
and let $u_i'=c_n^{-1}(u_i), v_i'=c_n^{-1}(v_i)$.

A general initial segment $w\rest I$ of a word $w\in \mcw^c_{n+1}$ has the following form with $q=q_n, k=k_n, l=l_n$. For some $0\le i_0\le q_n, 0\le j_0\le k_n$:
	\[\prod_{i<i_0}(\prod_{j<k}b^{q-j_i}w_j^{l-1}e^{j_i}) * (\prod_{j<j_0}b^{q-j_{i_0}}w_j^{l-1}e^{j_{i_0}}) * (b^{q^*}w_{j_0}^{l^*}w^*e^{j^*})\]
where $w^*$ is a possibly empty, possibly incomplete $n$-word,  $0\le j^*<j_{i_0}$, $0\le l^*\le l-1$, $0\le q^*\le q-j_{i_0}$. This is a block of complete 2-subsections, followed by a block of complete 1-subsections, followed by a possibly empty, incomplete 1-subsection. 

Similary a general tail segment  $w\rest T$ as the following form:
\[(b^{q^*}w^*w_{j_0}^{l^*}e^{j^*})* (\prod_{j_0\le j<k}b^{q-j_{i_0}}w_j^{l-1}e^{j_{i_0}})*\prod_{i_0\le i<q}(\prod_{j<k}b^{q-j_i}w_j^{l-1}e^{j_i})\]

\medskip

\bfni{Initial Segments:} We now argue for inequality \ref{14}. 
To start we take $n=0$. In this case $q_0=1$ and $q_1=k_0l_0$. The initial segment 
$w_i\rest I$ are of the form 
$\prod_{j<j_0}bw_j^{l_0-1}*u$ where $u$ is a proper initial segment of a word of the form $bw_{j_0}^{l_0-1}$ that has length $M$, for some $M<l_0$. 

If we throw away the tail segment $u$ we have thrown away proportion $M/\epsilon_0k_0l_0$. Since $M<l_0$ we have removed a portion of less than $\epsilon_0k_0$ and the segment $I_0$ that is left has proportion at least $1-({1/\epsilon_0k_0})$ and is made up of a product of $j_0$ many 1-subsections.

We now consider $n>0$. Since $\epsilon_{n}q_{n+1}=(\epsilon_nk_nl_nq_n)*q_n$, one of the following holds:
\begin{enumerate}
\item There are no complete 2-subsections, in which case we must have \\ $j_0+1>\epsilon_{n}k_nq_n$.
\item There is at least one complete 2-subsection and $j_0\ge\epsilon_nk_n$.
\item There is at least one complete 2-subsection and $j_0< \epsilon_nk_n$.
\end{enumerate}
In the first case, since $j_0+1>\epsilon_nk_nq_n$ we know that $j_0>\epsilon_nk_n$. Thus  eliminating the partial 1-subsection at the end we are left with a  concatenation of at least 
$\epsilon_nk_n$ complete 1-subsections and we have removed less than 
$1/\epsilon_nk_n$ portion of $I$. 
 Similarly in the second case we can eliminate the incomplete 1-subsection at the end by removing proportion less than $1/\epsilon_nk_n$ of $I$. In the final case  by removing  both the final incomplete 1-subsection and $(\prod_{j<j_0}b^{q-j_{i_0}}w_j^{l-1}e^{j_{i_0}})$ we 
eliminate {at most} 
{$1/q_n$ proportion of $I$.}

In all three cases, we are left an $I_0$ such that $w_0\rest I_0$ and $w_1\rest I_0$ are made up of a possibly empty initial segment of complete 2-subsections followed either by no complete 1-subsections or at least $\epsilon_nk_n$ complete 1-subsections. We now delete the boundary portions of $w_0\rest I_0$, which are aligned with the boundary portions of $w_1\rest I_0$. These have proportion $1/l_n$ in each complete 1-subsection--hence proportion $1/l_n$ of $I_0$. Let $I_1$ be the remaining portion of 
$I$. Then $I_1$ contains proportion at least {$(1-1/\epsilon_nk_n-1/q_n)(1-1/l_n)$} of $I$.

\medskip

{\bf\noindent Case 1:} $[w_0]_1\notin G^n_1[w_1]_1\footnote{We note that because $G^0_1=\la e \ra$, if $n=0$ we are in Case 1.}$.

Let $u'$ be the concatenation of $(u_0', u_1' \dots u'_{k_n-1})$, and $v'$ similarly the concatenation of the $v_i'$'s. Then $u', v'\in \mcw_{n+1}$ and $[u']_1\notin G^n_1[v']_1$. Let $u,v\in\mcw_{n}$ and $I^*$ be an initial or final segment of $\{ 0, 1, \dots , k_n-1\}$ of length at least $\epsilon_nk_n$. 
\begin{sublemma}\label{fortheref}
If $\epsilon_n$ is sufficiently small as a function of  $Q^n_1$, then 
	\[{\left|\{i\in I^*:[u_i']_1=[v_i']_1\}\right|\over |I^*|}\] 
is within ${1\over Q^n_1}$ of
	\[{1\over Q^n_1}.\]
\end{sublemma}

\pf Let 
$(u^*, v^*)$ be the concatenations of $\{u_i':i\in I^*\}$ and $\{v_i':i\in I^*\}$.
By \hyperlink{J11.1}{J11.1}, we see that  
the number $r(u,v)$ of occurrences of $(u,v)$ in $(u^*, {v^*})$  satisfies:
\begin{equation}
{r(u,v)\over |I^*|}\approx \left({1\over s_n}\right)^2 \label{8}
\end{equation}
Fix such an $I^*$ and let 
 $\mcc$ be a $\mcq^n_1$-class. Then $\mcc$ has $C^n_1$ elements. It follows from equation 
\ref{8} that the number of occurrences of a pair $(u,v)$ in $(u^*,v^*)$ with $u, v\in \mcc$ takes proportion of $|I^*|$ approximately 
\begin{equation*}
{(C^n_1)^2\over s_n^2}=\left({1\over Q^n_1}\right)^2
\end{equation*}
Since there are $Q^n_1$ many classes $\mcc$ that need to be considered we see that the number of pairs $u'_i$ and $v'_i$ with $[u'_i]_1=[v'_i]_1$ is approximately 
\begin{equation}\label{orphan}(1/Q^n_1)|I^*|.
\end{equation}
Hence for small enough $\epsilon_n$, we can see the conclusion of the sublemma. \qed
\begin{numreq}\label{eps and Qn1}
{The numbers $\epsilon_n$ should be small  enough as a function of $Q^n_1$ that estimate 
in the conclusion of sublemma \ref{fortheref} hold:
\begin{equation}\label{pia}
\left|{\left|\{i\in I^*:[u_i']_1=[v_i']_1\}\right|\over |I^*|}-{1\over Q^n_1}\right |<{1\over Q^n_1}.
\end{equation}
}
\end{numreq}

The locations in $w_0\rest I_1$ are made up of powers $u_i^{l-1}$. These fall into two categories, those locations occurring in whole 2-subsections and those occurring in the final product of 1-subsections. Applying the previous reasoning separately to the whole 2-subsections and the either-empty-or-relatively-long product of 1-subsections at the end of $I$, we see that the proportion of $u_i$ occurring in $w_0\rest I_1$ across from a $v_i$ in $w_1\rest I_1$ that is $\mcq^n_1$ equivalent is also extremely close to $1/Q^n_1$.  

If $n=0$, then specification \hyperlink{J11.1}{J11.1} implies that 
\[\left|\dbar(u^*, v^*)-{3\over 4}\right|<\epsilon_0.\]
So $\dbar(w_0\rest I_1,w_1\rest I_1)>(1-1/4-\epsilon_0)$ and hence 
\[\dbar(w_0\rest I, w_1\rest I)>\gamma_1.\]

In general, the induction hypothesis yields that 
$\mcq^n_1$-inequivalent words have $\dbar$-distance at least $\gamma_n$-apart.
Thus on $I_1$:
\begin{equation}\label{2 over Qn1}\dbar(w_0\rest I_1, w_1\rest I_1)>(1- 2/Q^n_1)\gamma_n.
\end{equation}

Allowing for agreement on boundary portions and applying Remark \ref{cheating on dbar}  we see that
{\[\dbar(w_0\rest I, w_1\rest I)\ge 
\left(1- 2\left({1\over Q^n_1} +{1\over \epsilon_nk_n}+{1\over q_n} +{1\over l_n}\right)\right)\gamma_n>\gamma_{n+1}.\]}

\medskip

{\bf\noindent Case 2:} $[w_0]_1\in G^n_1[w_1]_1$. 

In this case $n\ne 0$. Let $g\in G^n_1$ with $g[w_1]_1=[w_0]_1$. Since $[w_0]_1\ne[w_1]_1$, $g$ is not the identity. Since $G^n_1$ acts diagonally,   for all $i$ with $u_i$ intersecting the interval $I_1$, we have $[u_i]_1= g[v_i]_1$. In particular, $[u_i]_1\ne [v_i]_1$.

Hence  $\dbar(w_0\rest I_1, w_1\rest I_1)\ge \gamma_n$, and thus
 \[ \dbar(w_0\rest I, w_1\rest I)\ge 
\left(1- 2\left({1\over \epsilon_nk_n} +{1\over q_n}+{1\over l_n}\right)\right)\gamma_n>\gamma_{n+1}.\]

\medskip

\bfni{Tail Segments:} The argument for tail segments (inequality \ref{6}) follows the argument for initial segments, except that we delete small parts of the beginning of $T$, instead of the end of $I$.

\medskip
\bfni{Tail Segments compared to initial segments:} To show inequality \ref{7}, we proceed by induction, considering $w_0, w_1\in \mcw^c_{n+1}$.  In the comparing  two initial segments or two tail segments, not only did the 2 and 1-subsections  line up, but the $n$-subwords did as well. When comparing initial segments with tail segments, the $n$-subwords may be shifted, causing additional complications.
The proof  proceeds as in the easier cases, eliminating small sections of $I$ (or equivalently $T$) a bit at a time until we are left with $n$-words. The alignment of these $n$-words allows us to apply the induction  hypothesis  and conclude that the vast majority of $I$ and $T$ have $\dbar$-distance at least $\gamma_n$.
\smallskip

\begin{description}

\item{a.)} Of the 2-subsections of $w_0$ that intersect $I$, at most one is not a subset of $I$ (namely the last one), and similarly  except  for possibly the first 2-subsection intersecting $w_1\rest T$, $w_1\rest T$ is made up of whole 2-subsections.

\item{b.)}  Each 2-subsection of $w_0\rest I$ overlaps one or two 2-subsections of $w_1\rest T$.   An overlap of a 2-subsection of $w_0\rest I$ with a 2-subsection of $w_1\rest T$ that has  proportion bigger than $\epsilon_n$ of the 2-subsection implies that the overlap contains at least $\epsilon_nk_n$ complete 1-subsections.
    \begin{enumerate}
     \item Among the complete 2-subsections of $w_0\rest I$, delete overlaps of proportion less than $\epsilon_n$. 
     \item Delete the possible partial 2-subsection at the end of $w_0\rest I$ if it contains less than $\epsilon_nk_n$ complete 1-subsections.
    \end{enumerate}
The proportion of $I$ that has been deleted is less than $2\epsilon_n$.

\item{c.)} It could be that some of the portions of the remaining 2-subsections start or end with incomplete 1-subsections; i.e. not a whole word of the form $b^{q_n-j_i}v_j^{l_n}e^{j_i}$. Delete these incomplete sections. This leaves initial or tail segments of 2-subsections of the form $\prod_{j<k_n}b^{q_n-j_i}v_j^{l_n-1}e^{j_i}$ that consist of at least 
$\epsilon_nk_n$ whole 1-subsections. This trimming removes at most 
$1/k_n\epsilon_n$ proportion of $I$.

\item{d.)} We also remove the boundary sections of $w_0\rest I$. This removes at most $1/l_n$ of what remains of $I$ at this stage.

\item{e.)} We are left with a portion $I'\subset I$ such that   $w_0\rest I'$ consisting entirely of 0-subsections. These are blocks of the form $u_j^{l-1}$, where $u_j\in \mcw_n^c$.
 Each individual $n$-word $u_i$ can occur opposite a portion of $w_1\rest T$ in various ways. These 
 are:
\begin{description}
\item{i.} $u_i$ might occur exactly opposite a $v_{i+t}$\footnote{This is what happens in the case that $n=0$.} or
\item{ii.} $u_i$ might span portions of two copies of $v_{i+t}$ in  a power $v_{i+t}^{l-1}$. The two copies have the form $v_{i+t}v_{i+t}$, or
\item{iii.} $u_i$ might overlap a portion of the boundary of $w_1$. This can happen in two ways: boundary inside a 2-subsection (i.e. boundary of the form $e^{j_i}b^{q_n-j_i}$) and boundary between consecutive 2-subsections (i.e. boundary of the form $e^{j_i}b^{q_n-j_{i+1}}$). In each $u_i^{l_n-1}$ there are at most 3 copies of $u_i$ overlapping boundary portions of $w_1$.  

Hence by removing proportion at most {$4/l_n$} we are left with a portion  of $w_0\rest I$ consisting of powers of $u_j$'s that do not overlap any boundary in $w_1$.

\end{description}
\item{f.)} After the  deletions described in a.)-e.) the remaining portions of $w_0\rest I$ consists of blocks of powers of $u_i$'s in   initial segments of 2-subsections:
\begin{eqnarray*}u_0u_0\dots u_0*u_0\dots u_0\#u_1u_1\dots u_1*u_1 \dots u_1\# \dots\#\\ u_k\dots u_k*u_k\dots u_k
\end{eqnarray*}
and in tail segments of 2-subsections: 
\begin{eqnarray*}u_ju_j\dots u_j*u_j\dots u_j\#u_{j+1}u_{j+1}\dots u_{j+1}*u_{j+1} \dots u_{j+1}\# \dots\#\\ u_{k_n-1}\dots u_{k_n-1}*u_{k_n-1}\dots u_{k_n-1}
\end{eqnarray*}
where $*$'s stand for $u$'s deleted opposite boundary of $w_1$ and $\#$'s stand for the boundary of $w_0$ that has been deleted. An  important point for us is that in each block $k\ge \epsilon_nk_n$ and $k_n-j-1\ge \epsilon_nk_n$.

Consider the $u_j$'s in situation described in item e).ii. above. The $v_{i+t}$'s split $u_i$ into two pieces. By deleting a portion of the  individual $u_j$'s of size less than $\epsilon_{n-1}q_n$ we can assume that all of the overlap of $u_j$'s is in sections of length at least $\epsilon_{n-1}q_n$. By doing this for all $u_j$'s we remove a parts of the remaining elements of $w_0$ of proportion at most {$\epsilon_{n-1}$}.

\item{g.)} We now look more carefully at the two types of blocks  of  words described in item f.).  The case in item e.)i. is similar and easier than the case in item e.)ii. so we omit it.  Along the blocks described in f.) the initial portions of $u_i$ are lined up with $v_{i+t}$ and the second portions are lined up with $v_{i+t+1}$. Critically, the $t$ is \emph{constant} along the block.

According to whether $t=0$ or not, we apply specifications J11.1 (as in Case 1 of  the \emph{Initial Segments} argument) and J10
 to see that at most proportion {$2/Q^n_1$} of the $u_i$'s in a segment of the forms in f.) are lined up with $v_{i+t}$  are $\mcq^n_1$-equivalent. Hence we can make a final deletion of proportion at most $2/Q^n_1$ to get a portion $I^*\subseteq I$  consisting of relatively  long  pieces of $\mcw_n^c$-words in $w_0\rest I'$ overlapping $\mcw_n^c$-words in $w_1\rest T$ that lie in different $Q^n_1$ equivalence classes.

\end{description}

We now finish the argument using Remark \ref{cheating on dbar}. After all of the deletions we are left with $I^*$ having  at least
 $(1-(2\epsilon_n+1/\epsilon_nk_n+5/l_n+\epsilon_{n-1}+2/Q^n_1))$-proportion of $I$ and $w_0\rest I^*$ consists of relatively long pieces of $\mcw_n^c$ words that are overlapping portions of $\mcw_n^c$ words in $w_1\rest T$ that lie in different $\mcw^n_1$-classes. 

By the induction hypothesis each of the pieces of n-words in $w_0\rest I^*$ of $\dbar$-distance at least $\gamma_n$ from the corresponding portion of $w_1$.
 Consequently:
    \begin{eqnarray*}\dbar(w_0\rest I, w_1\rest T)&>&
    \gamma_n(1-(2\epsilon_n+1/\epsilon_nk_n+5/l_n+\epsilon_{n-1}+2/Q^n_1))\\
    &>&\gamma_{n+1}
    \end{eqnarray*}
thus finishing the proof of Proposition \ref{verify T4}.\qed
Since assumption \hyperlink{T4}{T4} is an immediate corollary of Proposition \ref{verify T4} we have finished verifying the timing assumptions.

We note in passing that inequality \ref{7} holds even if $w_0=w_1$ provided that the choice of initial and tail segment misalign corresponding 1-subsections.

We have proved: 
\begin{theorem}\label{timing is everything}
Suppose that $\bk^c$ is a system in the range of $F^s$ with construction sequence $\la \mcw_n^c:n\in\nn\ra$. Then $\la \mcw_n^c:n\in\nn\ra$ satisfies the timing assumptions.
\end{theorem}

\section{The Consistency of the Numerical Requirements}\label{numreqs}

During the course of this construction we have accumulated  numerical conditions about growth and decay rates of  several sequences. The majority of the numerical constants are not inductively determined--they are given immediately by knowing a small portion of the tree $\mct$. We call these \emph{exogenous} requirements. Other sequences of numbers depend on previous choices for the numbers--hence are determined recursively.  In this section we list the the recursive requirements, explicate their interdependencies and resolve their consistency.

Some of the conditions are easy to satisfy, as they don't refer to other sequences. For example, Numerical Requirement \ref{l_n summability} (that $\sum_n1/l_n<\infty$) can be satisfied once and for all by assuming that $l_n>20*2^{n}$. Others are trickier, in that they depend on the growth rates of the other sequences. For example, in defining the sequence of $k_n$'s we require that $k_n$ be large relative to the choice of $s_{n+1}$.  We call the former {type of }conditions \emph{Absolute} and the latter \emph{Dependent}.  The Dependent conditions introduce the risk of circular or inconsistent growth and decay rate conditions.  

Our approach here is to gather all of the conditions arising in this paper %
and its predecessors
and classify them as Absolute or Dependent.  We label them A or D accordingly. This process allows us to  make a diagram of the 
Dependent conditions to verify that there are no circularities. The lack of a cycle in the diagram gives a clear method of recursively 
satisfying all of the numerical conditions. 

Due to an overabundance of numerical parameters  we were forced into  some awkward notational 
choices. As noted before we have two types of epsilons:  the \emph{lunate} $\epsilon_n$, often used for set membership  and the \emph{classical} $\varepsilon_n$.  They play similar but slightly different 
roles. The lunate epsilons come from construction requirements arising in 
 \cite{FRW} {and their strengthenings. }
The  classical epsilons come from requirements related to circular systems and realizing them as smooth systems.
As is to be expected there is interaction between the two. This occurs via the intermediary numbers we called
 $\mu_n$'s in Numerical Requirements \ref{mu_n} and \ref{eps and mu}. 
 
 \subsection{The Numerical Requirements Collected.}\label{alltogethernow}
 In this section we collect the relevant numerical requirements. Specifically, in constructing $F^s(\mct)$ we are 
 presented with $\mct$ as a subsequence $\la \sigma_{n_i}:i\in\nn\ra$ of a fixed enumeration of $\nn^{<\nn}$.

In the formal statements of the specifications in \cite{FRW} for the construction sequence corresponding to $\mct$,  
$\mcw_n$ is built just in case $\sigma_n\in \mct$. This leads to a construction sequence of the form 
$\la \mcw_{n_i}:i\in\nn\ra$ with gaps corresponding to $m$'s where $\sigma_m\notin \mct$.  
To simplify notation we reindex 
$\la \mcw_{n_i}:i\in\nn\ra$ as $\la \mcw_i:i\in\nn\ra$ where  $\la \mcw_i:i<j\ra$ is determined by $\la \sigma_{n_i}:i<j\ra$.
In \cite{FRW}, the specifications
 discussed ``successive" (or ``consecutive") elements of $\mct$. These are $\sigma_m$ and $\sigma_n$ that belong to 
 $\mct$, but have no 
 $\sigma_j\in \mct$ with $j\in (m, n)$.  In our new notation successive elements $\sigma_m$ and $\sigma_n$ of $\mct$ 
 correspond to $\mcw_i$ and $\mcw_{i+1}$ where $m=n_i$. Having adopted this convention we don't distinguish between 
 $\la \mcw_i:i\in\nn\ra$ and 
 $\la \mcw_n:n\in\nn\ra$.
 To emphasize the dependence on $\mct$ we will occasionally write $\la \mcw_n(\mct):n\in\nn\ra$.

We begin with the  requirements inherited from \cite{FRW}.

\bigskip
\bfni{\underline{Inherited Numerical Requirements}}

 We have changed the notation from \cite{FRW} as described in the appendix \ref{appendixA}.   
 The number of elements of $\mcw_m$ is denoted $s_m$; 
 the numbers $Q^m_s$ and $C^m_s$ denote the
    number of classes and sizes of each class of $\mcq^m_s$ respectively. In \cite{FRW} we have sequences $\la \epsilon_n:n\in\nn\ra$, 
    $\la s_n,k_n, e(n), p_n:n\in\nn\ra$

\begin{description}
  
    \item[Inherited Requirement 1] $\la \epsilon_n:n\in\nn\ra$ is summable.
    \item[Inherited Requirement 2] $2^{e(n)}$ the number of $\mcq^n_{s+1}$ classes inside each
    $\mcq^n_s$ class. The numbers  $e(n)$ will be chosen to grow fast enough that 
    \begin{equation}
    2^n2^{-e(n+1)}<\epsilon_n \label{equ: crudite 2.1}
    \end{equation} If
    $s$ is the maximal length of an element of $\mct\cap \{\sigma_m:m\le n\}$ and 
    $|\mct\cap \{\sigma_m:m\le n\}|=i_0$ then 
     we set $C^{i_0}_s=2^{e(i_0)}$ as well. This forces $s_n, Q^n_s$ and $C^n_s$ all to be powers of 2 that are determined by $e(n)$. In particular let $\sigma_m$ and $\sigma_n$ be successive elements of $\mct$. Then $s_n$ is the number of words one gets by iteratively substituting $e(n)$ many elements into words in $\mcw_{n}^i/Q^n_i$ and closing under $G_i^m$ are successive for $i=0, 1, \dots s$.\footnote{It is possible to give a closed form formula for this, but it is complicated and uninformative.}        
    
    {By remark \ref{befuddle the reader} $s_n$ and $e(n)$ are monotonically co-determined.  Hence we can state this requirement as saying:
    \begin{quotation}
    $s_{n+1}$ is large enough in terms of $\epsilon_n$ that inequality \ref{equ: crudite 2.1} holds.
    \end{quotation}}
    
    \item[Inherited Requirement 3]     	
    	 If $\mct=\la
    	\sigma_{n_i}:i\in\nn\ra$ then 
    	\begin{equation}
    	2\epsilon_{i}s_{i}^2<\epsilon_{{i-1}}\label{equ: crudite 1.1}
    	\end{equation}
    	
  \item[Inherited Requirement 4]  	
  \begin{equation}\epsilon_{i}k_is_{{i-1}}^{-2}\to \infty\mbox{ as }i\to
    	\infty  \label{eqn: random prep2}
    	\end{equation}
    	
    \item[Inherited Requirement 5]	\item \begin{equation}
    	\prod_{n\in\nn}(1-\epsilon_n)>0 \label{equ: chisel2}
    	\end{equation} Since this  is equivalent to the summability of the $\epsilon_n$-sequence, it is   redundant  and we will ignore in the rest of this paper
    	
    	\item[Inherited Requirement 6] There will be  prime numbers $p_{n}$ such that $K_n=p_n^2s_{n-1}K_{n-1}$
    (i.e. $k_{n}=p_{n}^{2}s_{{n-1}}$). The $p_n$'s grow fast enough to allow the probabilistic arguments  in 
    	\cite{FRW} involving $k_n$ to go through.
	
	\item[Inherited Requirement 7] $s_n$ is a power of 2.
	
	\item[Inherited Requirement 8] The construction of $F(\mct)$ requires that if $\mct=\la \sigma_{i_n}:n\in\nn\ra$ then $\epsilon_n<2^{-i_n}$.

\end{description}

\bigskip
\bfni{\underline{Numerical Requirements introduced in this paper:}}
\medskip

\begin{description}
    \item[Numerical Requirement \ref{l_n summability}]
    $l_0>20$ and $\sum_{k=n} 1/l_k<1/l_{n-1}$.

    \item[Numerical Requirement \ref{varepsilons}]
    $\la \varepsilon_n:n\in\nn\ra$ is a  sequence of numbers in $\zoo$ such that {$6\sum_{n>N}\varepsilon_n<\varepsilon_N$.} 

    \item[Numerical Requirement  \ref{varepsilons and q's}] 
   {$k_n, l_n$ and $q_n$ grow fast enough that 
    $\varepsilon_nk_n\to \infty$, \\  $\varepsilon_nl_n\to \infty$, $\varepsilon_nq_n\to \infty$. }    
    
     \item[Numerical Requirement \ref{Gs and Qs}]
    $\sum{|G^n_1|\over Q^n_1}<\infty$ which is satisfied if ${|G^n_1|\over Q^n_1}<2^{-n}$.

    \item[Numerical Requirement \ref{mu_n}]
    $\mu_n$ is chosen small  relative to $\min(\varepsilon_n, 1/Q^n_1)$.

    \item[Numerical Requirement \ref{inherent in smooth}] The number 
     $l_n$ is big  enough relative to a lower   bound determined by $\la k_m, s_m:m\le n\ra$, $\la l_m:m<n\ra$ and $s_{n+1}$ to make the periodic approximations to the diffeomorphism converge.\footnote{{This is discussed in detail on pages 34-35 of \cite{prequel}, where the lower bound is called $l_n^*$.}}  Moreover $k_n\le l_n$.
    
    \item[Numerical Requirement \ref{sn grows}] 
    $s_n$ goes to infinity as $n$ goes to infinity and $s_{n+1}$ is a power of $s_n$. 
    \item[Numerical Requirement \ref{sn to kn}] $s_{n+1}\le s_n^{k_n}$.    
    \item[Numerical Requirement \ref{two epsilons}] The $\epsilon_n$'s are decreasing, $\epsilon_0<1/40$ and
    $\epsilon_n<\varepsilon_n.$

    \item[Numerical Requirement \ref{lln}]
    $k_n$ is chosen sufficiently large relative to a lower bound determined by $s_{n+1}, \epsilon_n$ that the Law of Large Numbers argument from \cite{FRW} works.

    \item[Numerical Requirement \ref{eps and mu}]
     $\epsilon_n$ is small  relative to $\mu_n.$

    \item[Numerical Requirement \ref{yeah yeah yeah}]
    $\epsilon_0k_0>20$, the $\epsilon_nk_n$'s are increasing and $\sum 1/\epsilon_nk_n<\infty$.

    \item[Numerical Requirement \ref{eps and Qn1}]
    {The numbers $\epsilon_n$ should be small  enough, as a function of $Q^n_1$, that estimate \ref{pia} holds.}
    
\end{description}

 \subsection{Resolution}

\label{resolution}

\begin{center}
{\bf \underline{A list of parameters, their first appearances and their constraints}}
\end{center}
We classify the constraints on a given sequence according to whether they refer to other sequences or not. Requirements that inductively refer to the same sequence are straightforwardly consistent. Those that refer to other sequences risk the possibility of being circular and thus inconsistent. As noted above refer to the former as \emph{Absolute} conditions and the latter as \emph{Dependent} conditions.

\begin{enumerate}
\item {\bf The sequence $\la k_n:n\in \nn\ra$.}  

\underline{Absolute conditions:} None for $\la k_n:n\in\nn\ra$.

\underline{Dependent conditions:}

    \begin{description}
    \item[D1] Numerical Requirement \ref{lln}, $k_n$ depends on $s_{n+1},\epsilon_n$.
   
    \item[D2] Inherited Requirement 6.  We can satisfy Inherited Requirement 6   by taking $k_n$ large enough to satisfy Numerical Requirement \ref{lln} and of the form $k_n=p_n^2s_{n-1}$.
    \item[D3] \label{you got it} From Inherited Requirement 4, equation \ref{eqn: random prep2} requires that $\epsilon_{n}k_ns_{{n-1}}^{-2}$ 
    goes to $\infty$ as $n$ goes to $\infty$. This can be satisfied by choosing $k_n$ large enough as a function of $\epsilon_n, s_{n-1}$.
    
    We note that equation \ref{eqn: random prep2} implies that $\sum 1/\epsilon_nk_n$ is finite. 
    \item[D4] Numerical Requirement 12 says that $\epsilon_0k_0>20$ and the $\epsilon_nk_n$'s are increasing and $\sum 1/\epsilon_nk_n$ is finite.  As noted the last condition follows from D3. The other parts of Numerical Requirement 12 are satisfied by taking $k_n$ large relative to $\epsilon_n$.
    \item[D5] Numerical Requirement \ref{sn to kn} implies that $k_n$ is large enough that $s_{n+1}\le s_n^{k_n}$. This implies that $k_n$ is large relative to $s_{n+1}$.
    \end{description}

{\bf From D1-D5, we see that $k_n$ is dependent on the choices of $\la k_m, l_m:m<n\ra, \la s_m:m\le n+1\ra$, and  
$\epsilon_n$.}

\item {\bf The sequence $\la l_n:n\in\nn\ra$.}

\underline{Absolute conditions}
    \begin{description}
    \item[A1] Numerical Requirement 1 says that  $1/l_n>\sum_{k=n+1}^\infty 1/l_k$. We also require that $l_n>20*2^n$, an exogenous requirement.
    \end{description}

\underline{Dependent conditions}

    \begin{description}
    \item[D6] By Numerical Requirement \ref{inherent in smooth}, $l_n$ is bigger than a number determined by  $\la k_m, s_m:m\le n\ra, \la l_m:m<n\ra$ and 
    $s_{n+1}$.

   \item[D7] {The sequence $\la l_n:n\in\nn\ra$ must grow fast  enough that $\varepsilon_{n+1}q_{n+1}\to \infty$. This can be arranged by making $\varepsilon_{n+1}q_{n+1}>n+1$. Since $q_{n+1}=k_nl_nq_n^2$, this puts lower bound on $l_n$ dependent on $\varepsilon_{n+1}$.}
    \end{description}
    
     {\bf Thus $l_n$ depends on $\la k_m, s_m: m\le n\ra$, $\la l_n:m<n\ra$, $\varepsilon_{n+1}$ and $s_{n+1}$. }
\item {\bf The sequences $\la s_n:n\in\nn\ra$ and $\la e(n):n\in\nn\ra$.}
We treat these sequences as equivalent since $s_n$ is a power of $2$ determined by $e(n)$ and the elements of the tree in the domain of the reduction. Moreover increasing one increases the other and vice versa. Since they are co-determined, they are chosen at the same time.

\underline{Absolute conditions} 
\begin{description}
\item[A2] Inherited Requirement 7 says that $s_n$ is a power of 2.
\end{description}
Numerical Requirement 7 says that:
\begin{description}
\item[A3] The sequence $s_n$ goes to infinity.
\item[A4] $s_{n+1}$ is a multiple of $s_n$.
\item[A5]  { Since $e(n)$ determines $Q^n_1$, Numerical Requirement \ref{Gs and Qs} puts an exogenous sequence of lower bounds on $e(n)$, for example that $|G^n_1|/Q^n_1<2^{-n}$. This requires that $e(n)$ be chosen large and, since $e(n)$ and $s_n$ are inter-determined,  can be satisfied by taking $s(n)$ large.}

\end{description}

\underline{Dependent conditions}

\begin{description}

\item[D8] Numerical Requirement 3 makes $s_n$ depend on $\epsilon_{n-1}$.
\end{description}

{\bf The result is that $s_{n+1}$ depends on the first $n+1$ elements of $\mct$, $\la k_m, s_m, l_m:m<n\ra$, $s_n$,  and $\epsilon_n$.}\footnote{It is important to observe that the choice of $s_{n+1}$ does \emph{not} depend on $k_n$ or $l_n$.} 

\item {\bf The sequence $\la \epsilon_n:n\in\nn\ra$.}

\underline{Absolute conditions}
    \begin{description}
    \item[A6] Numerical Requirement \ref{two epsilons} and Inherited Requirement 1  require that the $\la \epsilon_n:n\in\nn\ra$ is decreasing and summable and $\epsilon_0<1/40$.
    \item[A7]  Inherited Requirement 8 says that if $\mct=\la \sigma_{i_n}:n\in\nn\ra$ then $\epsilon_n<2^{-i_n}$
    \end{description}

\underline{Dependent conditions}

    \begin{description}
    \item[D9] Numerical Requirement \ref{two epsilons} requires that  $\epsilon_n<\varepsilon_n$. 
    \item[D10] Equation \ref{equ: crudite 1.1} of Inherited Requirement 3 says $2\epsilon_{n}s_{n}^2<\epsilon_{n-1}$
    \item[D11] Numerical Requirement \ref{eps and mu} says that $\epsilon_n$ must be small enough relative to $\mu_n$.
    \item[D12] Numerical Requirement \ref{eps and Qn1} says that $\epsilon_n$ is small as a function of $Q^n_1$.
    \end{description}

   The result is that $\epsilon_n$ depends exogenously on the first $n$ elements of $\mct$, and on $Q^n_1, s_n$, $\varepsilon_n$, $\epsilon_{n-1}$ and $\mu_n$.

\item {\bf The sequence $\la \varepsilon_n:n\in\nn\ra$. }

\underline{Absolute conditions} 
	\begin{description}
	\item[A8] {Numerical Requirement 2 says that $6\sum_{n>N}\varepsilon_n<\varepsilon_N$. This can be arranged by taking $\varepsilon_{n}<12^{-n}\varepsilon_{n-1}$.}
	\end{description}

\underline{Dependent conditions}

{Numerical Requirement \ref{varepsilons and q's} imposes three Dependent conditions on $\varepsilon_n$: 
$\varepsilon_nk_n\to \infty$, \\  $\varepsilon_nl_n\to \infty$, $\varepsilon_nq_n\to \infty$. We deal with these in turn.}

\begin{enumerate}

\item { The requirement that $\la \varepsilon_nk_n:n\in\nn\ra$ goes to infinity already follows from the fact that 
$\epsilon_n<\varepsilon_n$ and item D4.}
\item {$\la \varepsilon_nl_n:n\in\nn\ra$ goes to infinity.  This follows from $k_n\le l_n$, which is covered in Dependent condition D6.}
\item {$\la \varepsilon_nq_n:n\in\nn\ra$ goes to infinity. This follows from Dependent condition D7.}
\end{enumerate}

{\bf Thus there are no new  Dependent conditions.}

\item {\bf The sequence $\la Q^n_1:n\in\nn\ra$.}

\underline{Absolute conditions:} There are no new Absolute conditions. 

\underline{Dependent conditions}
	\begin{description}
	\item[D13] { Numerical Requirement \ref{Gs and Qs} says that ${|G^n_1|\over Q^n_1}<2^{-n}$. 
	But since  $Q^n_1$ is determined by $s_n$ and the  first $n$-elements of the tree, Numerical requirement \ref{Gs and Qs} is taken 
	care of by A5.}
	\end{description}

{\bf There are no new  Dependent conditions.}

\item {\bf The sequence $\la \mu_n:n\in\nn\ra$.} 

This sequence gives the required pseudo-randomness in the timing assumptions.

\underline{Absolute conditions:} There are no new Absolute conditions.

\underline{Dependent conditions}

\begin{description}
\item[D14] Numerical Requirement \ref{mu_n} requires that 
$\mu_n$ be very small relative to $\varepsilon_n$ and $1/Q^n_1$. 
\end{description}
{\bf $\mu_n$ is dependent on $\varepsilon_n$ and $Q^n_1$.}

\end{enumerate}

\noindent The recursive dependencies of the various coefficients are summarized in Figure \ref{ecosystem}, in which an 
arrow from a coefficient to another coefficient shows that the latter is dependent on the former. Here is the order the the 
coefficients can be chosen consistently.

\subsection{{The inductive order of choices}}
We begin by setting: $s_0=2, s_1=8, p_0=0, q_0=k_0=1, l_0=21$. $Q^0_1$ is not defined, but $Q^1_1$ is determined by $s_1$.
 $\mu_0=\epsilon_0=k_0=l_0=1$, $\varepsilon_0=1.1$, $\varepsilon_1=\varepsilon_0/12$, 
\medskip

\bfni{Assume:}
\begin{quotation}
\noindent The coefficient sequences $\la k_m, l_m, Q^m_1, \mu_m, \epsilon_m:m<n\ra $,  $\la \varepsilon_m:m\le n\ra$ and $s_{n}$ have been chosen. The first $n+1$ sequences on the tree $\mct$ are known.
\end{quotation}
\medskip

\bfni{To do:}
\begin{quotation}\noindent Choose $k_n, l_n, Q^n_1, \mu_n, \epsilon_n, \varepsilon_{n+1}$ and $s_{n+1}$. Each requirement is to choose the corresponding variable \emph{large enough} or \emph{small enough} where these adjectives are determined by the dependencies enumerated above. 
\end{quotation}

Figure \ref{ecosystem} gives an order to consistently choose the next elements on the sequences; Choose the successor coefficients  in the following order:

\[Q^n_1, \varepsilon_n, \mu_n,\epsilon_n, s_{n+1}, k_n, l_n.\]

\begin{figure}
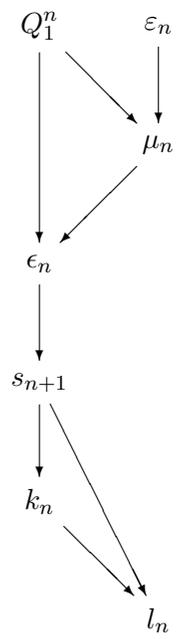

\centering

\[\begin{diagram}
\node{Q^n_1}\arrow[2]{s}\arrow{se}\node{}\\%
\node{}\node{\mu_n}\arrow{sw,r}{}\\
\node{\epsilon_n}\arrow{s}\node{}\\
\node{s_{n+1}}\arrow{s}\arrow{sse}\node{}\\
\node{{k_n}}\arrow{se}\node{{\varepsilon_{n+1}}}\arrow{s}{}\\
\node{}\node{l_n}
\end{diagram}\]

\caption{Order of choice of Numerical parameters dependency diagram.}
\label{ecosystem}
\end{figure}
\noindent We note that $Q^n_1$ is redundant in the diagram above since it is determined by $s_n$, but we include it as a bridge from stage $n-1$.

\appendix

\begin{center}
\bfni{Appendix} 
\end{center}
\section{Notation table}\label{appendixA}
 In this paper we have adopted the notation used in \cite{AK}, which conflicts with the notation in \cite{FRW}, accordingly we provide a table for translating between the two. In the table, \emph{NEW} means the notation used in this paper, \emph{OLD} means the notation used in \cite{FRW}.

\begin{eqnarray*}
\mbox{\underline{\bf NEW}} &  \mbox{\underline{\bf OLD}}& \mbox{\underline{\bf Description}}\\
s_n&  W_n & s_n  \mbox{ is the number of words in }\mcw^c_n\\\
k_n& l_{n+1}/l_n& \mbox{The number of words concatenated  } \\
& & \mbox{to make $\mcw_{n+1}$ from $\mcw_n$}\\
e(n)&  k(n)& \mbox{Controls the number of $\mcq_{s+1}$ classes}\\
& & \mbox{in each $\mcq_s$ class}\\
\gamma&  s_1& \mbox{The separation between $\mcq^n_1$ classes}\\
K_n&l_n&\mbox{$K_n$ is this paper's notation for the lengths of the}\\
&&\mbox{ odometer based words in $W_n,$}\\
&&\mbox{$l_n$ was the notation for the lengths of the words in \cite{FRW}}\\
q_n&l_n&\mbox{The lengths of the circular words in current paper}\\ &&\mbox{vs. odometer based words in   \cite{FRW}. The new $q_n$ }\\ 
&&\mbox{refers to the lengths of the words in $\mcw_n^c$.}\\
l_n& \mbox{no analogue}&\mbox{Coefficient needed to grow fast for smooth transformations}
\end{eqnarray*}
An equivalent description of the  numbers we are calling $k_n$ in this paper is that they are the number of words in 
$\mcw_n^c$ concatenated to form elements of $P_{n+1}$. The number $k_n$ is equal  to the number $K_{n+1}/K_n$ and $l_{n+1}/l_n$ in the old
notation of \cite{FRW}.

 \bibliography{citations}
\bibliographystyle{plain}

 \end{document}